\documentclass[12pt]{article}

\usepackage{amscd}
\usepackage{amsfonts}
\usepackage{amsmath}
\usepackage{amssymb}
\usepackage{amsthm}
\usepackage{comment}

\parskip=\medskipamount

\newdimen\Squaresize \Squaresize=20pt
\newdimen\thickness \thickness=1pt         

\font\cornerfont=cmr9

\def\Square#1{\hbox{\vrule width \thickness
   \vbox to \Squaresize{\hrule height \thickness\vss                                  \hbox to \Squaresize{\hss#1\hss}
   \vss\hrule height\thickness} 
\unskip\vrule width \thickness} 
\kern-\thickness}                                                            
                               
\def\vsquare#1{\vbox{\Square{$#1$}}\kern-\thickness}
\def\blank{\omit\hskip\Squaresize}

\def\young#1{\let\\=\cr	%added  let\\ =\cr 
\vbox{\smallskip\offinterlineskip
\halign{&\vsquare{##}\cr #1}}}

\def\smyoung#1{\let\\=\cr\Squaresize=15pt
\vbox{\smallskip\offinterlineskip
\halign{&\vsquare{##}\cr #1}}}

\def\fibyoung#1{\let\\=\cr		%added  let\\ =\cr 
\vbox{\smallskip\offinterlineskip
\halign{&\vsquare{##}\cr #1}}\,}

\newbox\vstrutbox
\setbox\vstrutbox=\hbox{\vrule height 0pt depth 0pt width 4pt}
\def\vstrut{\relax \ifmmode\copy\vstrutbox\else\unhcopy\vstrutbox\fi }
\newbox\smvstrutbox
\setbox\smvstrutbox=\hbox{\vrule height 0pt depth 0pt width 1pt}
\def\smvstrut{\relax \ifmmode\copy\smvstrutbox\else\unhcopy\smvstrutbox\fi }

\newbox\reghstrutbox
\setbox\reghstrutbox=\hbox{\vrule height 9.5pt depth 3.5pt width 0pt}
\def\reghstrut{\relax \ifmmode\copy\reghstrutbox\else\unhcopy\reghstrutbox\fi }
\newbox\smhstrutbox
\setbox\smhstrutbox=\hbox{\vrule height 4pt depth 2pt width 0pt}
\def\smhstrut{\relax \ifmmode\copy\smhstrutbox\else\unhcopy\smhstrutbox\fi }
\newdimen\fsquaresize \fsquaresize=40pt
\newdimen\fthickness 

\def\fsquare#1{\hbox{\vrule width \fthickness
   \vbox to \fsquaresize{\hrule height \fthickness
	\hbox to \fsquaresize{\hss\hstrut #1\smvstrut}
   \vss\hrule height\fthickness} 
\unskip\vrule width \fthickness} 
\kern-\fthickness}                                                            
                               
\def\fvsquare#1{\vbox{\fsquare{$#1$}}\kern-\fthickness}
\def\dfsquare#1{\hbox{\vrule width \fthickness
   \vbox to \fsquaresize{\hrule height \fthickness
  \vfill \hbox to \fsquaresize{\hfill\hstrut \hbox{\cornerfont #1}\smvstrut}
   \hrule height\fthickness} 
\unskip\vrule width \fthickness} 
\kern-\fthickness}                                                            
                               
\def\dfvsquare#1{\vbox{\dfsquare{#1}}\kern-\fthickness}
\def\fborder#1{\hbox{\vrule width 0pt
   \vbox to \fsquaresize{\hrule height 0pt
	\hbox to \fsquaresize{\hss \hstrut #1\vstrut }	
   \vss\hrule height 0pt} 		%had to mess with kerns;
\unskip\vrule width 0pt} 
\kern-4\fthickness} 
                                                         
\def\ftopborder#1{\hbox{\vrule width 0pt
   \vbox to \fsquaresize{\hrule  height 0pt\vss
	\hbox to \fsquaresize{ \hstrut #1\smvstrut \hss }	
   \hrule height 0pt} 		%had to mess with kerns;
\unskip\vrule width 0pt} 
\kern-\fthickness}

\def\fvborder#1{\vbox{\fborder{$#1$}}\kern-4\fthickness}

\def\fvtopborder#1{\vbox{\ftopborder{$#1$}}\kern-4\fthickness}
\def\fvcborder#1{\hbox{\vrule width 0pt
   \vbox to \fsquaresize{\hrule  height 0pt\vss
	\hbox to \fsquaresize{ \hstrut #1\vstrut \hss }	
   \vss\hrule height 0pt} 		%had to mess with kerns;
\unskip\vrule width 0pt} 
\kern-\fthickness}

\def\fvborder#1{\vbox{\fborder{$#1$}}\kern-4\fthickness}

\def\fvvcborder#1{\vbox{\fvcborder{$#1$}}\kern-4\fthickness}
\font\cellfontfive=cmssbx10 scaled \magstep5
\font\cellfontfour=cmssbx10 scaled \magstep4
\font\cellfontthree=cmssbx10 scaled \magstep3
\font\cellfonttwo=cmssbx10  scaled \magstep2
\font\cellfontone=cmssbx10  scaled \magstep1
\font\cellfontzero=cmssbx10  
\font\sevenrm=cmr7
\font\sixrm=cmr6
\font\fiverm=cmr5
\def\cellfont{\cellfontfive}

\def\xsquare#1#2{\hbox{\vrule width \fthickness
   \vbox to \fsquaresize{\hrule height \fthickness
	\vbox {
	   \hbox to \fsquaresize {\hss\hstrut {\cornerfont #2}\smvstrut }
	   \vss}
	\vss
	   \hbox to \fsquaresize{\hss\hbox{\cellfont #1}\hstrut \hss}
   	\vss  \hrule height\fthickness
	} 
\unskip\vrule width \fthickness} 
\kern-\fthickness}                                                            
                               
\def\xvsquare#1.#2.{\vbox{\xsquare{#1}{#2}}\kern-\fthickness} %?Why $$ ?
\def\Xsquare#1#2{\hbox{\vrule width \fthickness
   \vbox to \fsquaresize{\hrule height \fthickness \kern -5.5\fthickness
	\hbox to \fsquaresize{\hss 
	   \vbox to \fsquaresize{\vss\hbox{\cellfont #1}\vss} 
	  \hss
   	   \vbox to \fsquaresize{\vss \hbox{\hstrut \cornerfont #2}}
	}
    \hrule height\fthickness}
\unskip\vrule width \fthickness} 
\kern-\fthickness}

\def\Xvsquare#1.#2.{\vbox{\Xsquare{#1}{#2}}\kern-\fthickness} %?Why $$ ?

\def\fyoung#1{  
	\def\hstrut{\reghstrut}
	\fthickness=1pt  		%.025\fsquaresize  %%Don't need.
\def\>{\omit\fvborder}
\def\vc{\omit\fvvcborder}
\def\<{\omit\fvtopborder}
\def\x{\omit\xvsquare}	
\def\X{\omit\Xvsquare}	
\let\\=\cr
\def\blank{\omit\hskip\fsquaresize}
\vbox{\smallskip\offinterlineskip
\halign{&\fvsquare{##}\cr #1}}}	

\def\dfyoung#1{  
	\def\hstrut{\reghstrut}
	\fthickness=1pt  		%.025\fsquaresize  %%Don't need.
\def\>{\omit\fvborder}
\def\vc{\omit\fvvcborder}
\def\<{\omit\fvtopborder}
\def\x{\omit\xvsquare}	
\def\X{\omit\Xvsquare}	
\let\\=\cr
\def\blank{\omit\hskip\fsquaresize}
\vbox{\smallskip\offinterlineskip
\halign{&\dfvsquare{##}\cr #1}}}	

\def\smfyoung{\fsquaresize=30pt 
		\def\cellfont{\cellfontthree}  \def\cornerfont{\sevenrm}
				\fyoung}
\def\smdfyoung{\fsquaresize=19pt 
		\def\cellfont{\cellfontzero}
		\def\cornerfont{\sevenrm}\dfyoung} 
\def\textfyoung{\fsquaresize=12pt 
		\def\cellfont{\fiverm}  \fyoung}

\def\smxsquare#1#2{\hbox{\vrule width \fthickness
   \vbox to \fsquaresize{\hrule height \fthickness
	\vss	\hbox to \fsquaresize{\hss\hbox{\cellfont #1}\hstrut \hss}
   \vss  \hrule height\fthickness} 
\unskip\vrule width \fthickness} 
\kern-\fthickness}

\def\smXsquare#1#2{\hbox{\vrule width \fthickness
   \vbox to \fsquaresize{\hrule height \fthickness\kern %2\fthickness
	\vbox {\vss
	 	   \hbox to \fsquaresize{\hstrut\hss\hbox{\cellfont #1} \hss}
   	 \hbox to \fsquaresize {\hss\hstrut {\mit #2}\smvstrut }
	   }
	} 
\unskip\vrule width \fthickness} 
\kern-\fthickness}                                                            

\def\ssmfyoung{\fsquaresize=22pt \def\hstrut{\reghstrut}
		\def\cellfont{\cellfontone}
	\def\xsquare{\smxsquare}\def\Xsqaure{\smXsquare }
\fyoung}
\def\fomin{\fyoung}

\def\zz{\mathbb}
\def\mit{\rm}
\def\Cal#1{\mathcal{#1}}
\def\rs{\quad \stackrel{\rm R-S}{\longleftrightarrow}\quad }
\def\berele{\underset{\Cal B}{\longleftrightarrow}}%%  \quad seemed too big (IT)
\def\Bereleinsert{\underset{\mathcal B}{\leftarrow}}
\def\bet{\beta}
\def\cell#1#2#3#4#5#6#7#8#9{%
\vcenter{\baselineskip=0pt \lineskip=2pt \lineskiplimit=2pt \tabskip=2pt
\halign{&\hfil##\hfil\cr
&$#1$&$#2$&$#3$\cr
&\vbox to 3em{\vfil\hbox{$#4$}\vfil}%
&\vrule\vbox to 3em{\hrule\vfil\hbox to 3em{\hfil$#5$\hfil}\vfil\hrule}\vrule%
&\vbox to 3em{\vfil\hbox{$#6$}\vfil}\cr
&$#7$&$#8$&$#9$\cr}}}
\def\cover{\overset{\textstyle.}{\supset}}
\def\covered{\overset{\textstyle.}{\subset}}
\def\coveredrow#1{\overset{#1}{\subset}}
\def\coveredrowbelow#1{\underset{#1}{\subset}}
\def\coverrow#1{\overset{#1}{\supset}}
\def\coverrowbelow#1{\underset{#1}{\supset}}
\def\Cpx{\mathbb C}
\def\cross{{\times}}
\def\em#1{{\bf#1}}
\def\emp{\varnothing}
\def\gam{\gamma}
\def\Gam{\varGamma}%%%%%%%%% (IT)
\def\hdom{{\def\hr{\hrule width13.2pt}\def\vr{\vrule height6pt}%
\def\hs{\hskip6pt}%
\vcenter{\hr\hbox{\vr\hs\vr\hs\vr}\hr}}}
\def\inv{^{-1}}
\def\kap{\kappa}
\def\lam{\lambda}
\def\Lam{\varLambda}%%%%%%%% (IT)
\def\ord{\operatorname{ord}}
\def\packed{\offinterlineskip\tabskip=0pt}
\def\Par{\mathcal P}
\catcode`\@=11
\def\revddots{\mathinner{\mkern1mu\raise\p@\vbox{\kern7\p@\hbox{.}}\mkern2mu
\raise4\p@\hbox{.}\mkern2mu\raise7\p@\hbox{.}\mkern1mu}}
\def\sh{\operatorname{sh}}
\def\StandBereleinsert{\underset{\tilde{\Cal B}}{\leftarrow}}
\def\th#1{^{(#1)}}
\def\vcover{\overset{.}{\cup}}%%%%%%%%%%%%%%%%%%%% TO BE CHANGED %%%%%%%%%%
\def\vcovered{\overset{.}{\cap}}
\def\vcoveredrow#1{\overset{#1}{\cap}}
\def\vcoverrow#1{\overset{#1}{\cup}}
\def\vdom{{\def\hr{\hrule width6.8pt}\def\vr{\vrule height6pt}%
\def\hs{\hskip6pt}%
\vcenter{\hr\hbox{\vr\hs\vr}\hr\hbox{\vr\hs\vr}\hr}}}

\newtheorem{thm}{Theorem}[section]

\newtheorem{lemma}[thm]{Lemma}

\theoremstyle{definition}
\newtheorem{rem}[thm]{Remark}
\theoremstyle{definition}
\newtheorem{eg}[thm]{Example}	%This uses the def above
\newtheorem{defn}[thm]{Definition} %This also uses the def above
\def\rom#1{{\rm#1}}
\def\ssize{\scriptstyle}

\def\B#1{$\bar{\text{#1}}$}
\def\cyr#1{#1}

\title{A two-dimensional pictorial presentation\\
of Berele's insertion algorithm for symplectic tableaux}

\author{Tom Roby\\  
\small Department of Mathematics\\[-0.8ex] 
\small California State University\\[-0.8ex] 
\small Hayward, CA 94542, USA\\[-0.8ex]  
\small \texttt{troby@csuhayward.edu}\\[0.8ex]  
Itaru Terada\\ 
\small Graduate School of Mathematical Sciences\\[-0.8ex] 
\small University of Tokyo\\[-0.8ex] 
\small Komaba 3-8-1, Meguro-ku\\[-0.8ex] 
\small Tokyo 153-8914, Japan\\[-0.8ex] 
\small \texttt{terada@ms.u-tokyo.ac.jp}}
\date{\small May 17, 2004\\ 
\small MR Subject Classifications: 05E10, 05E15, 17B20, 20G05, 22E46}

\begin{document}

\maketitle

\begin{abstract}

We give the first two-dimensional pictorial presentation of Berele's
correspondence \cite{Berele}, an analogue of the Robinson-Schensted (R-S)
correspondence \cite{Robinson, Schensted} for the symplectic group
$Sp(2n, \Cpx )$.  From the standpoint of representation theory, the R-S
correspondence combinatorially describes the irreducible decomposition
of the tensor powers of the natural representation of $GL(n,\Cpx)$.
Berele's insertion algorithm gives the bijection that describes the 
irreducible decomposition of the tensor powers of the natural
representation of $Sp( 2n, \Cpx )$.  Two-dimensional
pictorial presentations of the R-S correspondence via local rules (first
given by S.~Fomin \cite{Fomin,FominGen}) and its many variants have
proven very useful in understanding their properties and creating new
generalizations.  We hope our new presentation will be similarly
useful.  
\end{abstract}

\section{Introduction}\label{sec:intro}

Our purpose is to give a new presentation of Berele's correspondence.

Berele's correspondence is a combinatorial construction devised by A.~Berele
in \cite{Berele}, as an $Sp(2n,\Cpx)$-analogue of one aspect of
the Robinson-Schensted correspondence, or the R-S correspondence for
short.  The R-S correspondence describes the irreducible decomposition
of the representation of the group $GL(n,\Cpx)$ on $(\Cpx^n)^{\otimes f}$
(where $f$ is a fixed positive integer)
derived from its natural action on the column vectors of $\Cpx^n$.

Similarly, Berele's correspondence describes the irreducible decomposition
of the representation of the group $Sp(2n,\Cpx)$ on $(\Cpx^{2n})^{\otimes f}$
also derived from its natural action on the column vectors in $\Cpx^{2n}$,
at least on the character level.
A further analysis of Berele's correspondence
was conducted by S.~Sundaram in her thesis \cite{Sundaram}.  (See also \cite{SundaramJCTA} and \cite[Theorem 3.10 and Appendix]{SundaramIMA}.)

While many interesting connections have been found between the R-S
correspondence and various algebraic and geometric objects, the
appearance of Berele's correspondence has been relatively limited.  We show
in this article that one more aspect of the R-S correspondence has its
counterpart for Berele's correspondence.

S.~Fomin \cite{Fomin,FominGen} showed
that the R-S correspondence can be presented
as a two-dimensional inductive application of ``local rules'',
which are based on the properties of Young's lattice $\Par$
(the poset of all partitions, ordered by containment of diagrams)
as a ``Y-graph'' or ``differential poset'' \cite{StanleyDif}.
T.~Roby \cite{Robyone} generalized this interpretation
to several variants of the R-S correspondence.  The local rules can be
derived directly from the original procedural definition of the R-S
correspondence given in \cite{Schensted} or \cite{Knuth} (first
appearance in \cite{Robinson}), as lucidly explained in
\cite{vanLeeuwen} by M. van~Leeuwen.
It is this type of analysis that we apply to Berele's correspondence
in this article.

Another important ingredient of Berele's correspondence
is Sch\"utzenberger's jeu de taquin or sliding algorithm \cite{Schut}.
Fomin, and later van Leeuwen, gave a local rules presentation
of this algorithm. 
A widely available treatment of Fomin's local rules approach to the R-S
correspondence and jeu de taquin can be found in  Section 7.13 of
\cite[Section 7.13, Appendix~1]{StanleyECII}.
We have been inspired by their work to extend the set of local rules
and create a ``stratification'' of the diagram that allows Berele's
correspondence to be presented pictorially.

The local rules thus extended turn out to have interesting
symmetries.  We think that the procedure defined by these local rules
are of intrinsic interest and deserves more investigation.  It would
also be interesting to connect our
algorithm with a poset invariant like the Greene-Kleitman correspondence
(\cite{Greene} or \cite{Greenetwo}), with geometric or Lie group theoretic
objects like flags or their generalizations, or with a precise
interpretation in terms of a quantum analogue of $Sp(2n,\Cpx)$; all
these still remain to be explored.

In Section~\ref{sec:berins} we review Berele's original approach to his
correspondence via bumping and jeu de taquin.  In
Section~\ref{sec:berlr} we describe the extended set of local rules and
the stratification of the diagram necessary to present Berele's
algorithm pictorially.  In Section~\ref{sec:rev} we describe the
procedures to handle the reverse correspondence.  This is more
complicated than the original R-S case, where one of the most satisfying
aspects of the pictorial description is the transparentness of
bijectivity.  Finally in Section~\ref{sec:oqr} we make some remarks and
mention directions for future research.

The authors are grateful to the Japan Society for Promotion of
Science for supporting the first author's postdoc at the University of
Tokyo.  We thank the departments at the University of Tokyo and MIT for
their hospitality.  We particularly benefited from conversations with
Sergey Fomin, Kazuhiko Koike, and Marc van Leeuwen.   

\section{Review of Berele's Correspondence by Insertion}\label{sec:berins}

Throughout this article, an interval $[i,j]$ will
be taken inside the ordered set $\zz Z$ of integers.
\subsection{Partitions}\label{subsec:part}
A {\bf partition} $\lam$ is a weakly decreasing sequence of nonnegative
integers $\lam = ( \lam_1, \lam_2, \lam_3, \ldots )$,
$\lam_1 \ge \lam_2 \ge \lam_3 \ge \cdots$ with only a finite number of
nonzero terms (called the {\bf parts} of $\lam$).
The number of parts of $\lam$ is called the {\bf length} of $\lambda $ and is
written $l (\lambda )$.  The sum of all the parts of $\lam$ is called the
{\bf weight} of $\lam$, and denoted by $|\lambda|$.
In writing concrete partitions, we generally suppress trailing zeros.
Moreover, in the figures below, we sometimes omit parentheses and
commas.  %%%}(IT)
Since no parts greater than 9 occur in any
of the examples, no confusion should result.
The unique partition of weight 0 is denoted by
$\varnothing$ or %%%%%%%%%%%%%%%%%%%%%%%%%%%%%%%%%%%%%%%%%%%%%%% (IT)
$0 $.
The set of all partitions will be denoted by $\Par$.

The {\bf (Young) diagram} of a partition $\lam$ is formally the set
$D_\lam = \{\, ( i, j ) \in \mathbb N^2 \mid 1 \le i \le l( \lam ), \;
1 \le j \le \lam_i \,\}$, which is sometimes identified with $\lam$ itself.
Each $( i, j ) \in D_\lam$ is called its {\bf square} or {\bf cell},
and we may visualize $D_\lam$ as a cluster of contiguous square boxes,
each representing a ``square'' $( i, j )$,
arranged in a matrix-like order.  (See Figure~\ref{fig:ydtab}(a).)

Define a partial order $\subseteq $ on partitions by $\mu
\subseteq \lambda $ if and only if $D_{\mu }\subseteq D_{\lambda }$.
This turns $\Par$ into a distributive lattice, called {\bf Young's lattice}.
We say that ``$\lambda $ \em{covers} $\mu
$'' and write $\lambda \cover \mu $
or $\mu\covered\lambda$ %%%%%%%%%%%%%%%%%%%%%%%%%%%%%%%%%%%%%%%%%%%%% (IT)
if
$\mu \subseteq \lambda $ and they differ by exactly one square.  We call
such a square a {\bf corner} of $\lambda $ and a {\bf cocorner} of $\mu
$ (following van Leeuwen).
If the difference lies in the $k$th row,
then we also write $\lam\coverrow{k}\mu$ or $\mu\coveredrow{k}\lam$.

\begin{figure}
\caption{A Diagram and a Tableau}
\label{fig:ydtab}
$$
\gathered
D_{\lam} =
\vcenter{\packed
\def\cell#1{\vbox to 1.5em{\vss\hbox to 1.5em{\hss$#1$\hss}\vss}}
\halign{\vrule#&&\cell{#}\vrule\cr\noalign{\hrule}
& & &\bullet& \cr\noalign{\hrule}
& & & \cr\multispan4\hrulefill\cr
& & & \cr\multispan4\hrulefill\cr
& \cr\multispan2\hrulefill\cr
}
} ( = \lam ) \\
\vspace{\medskipamount}
\bullet=\text{square }(1,3)\\
\vspace{\medskipamount}
\text{(a) The Young diagram of } \lam=(4,3,3,1)
\endgathered
\qquad
\gathered
T=
\vcenter{\packed
\def\cell#1{\vbox to 1.5em{\vss\hbox to 1.5em{\hss$#1$\hss}\vss}}
\halign{\vrule#&&\cell{#}\vrule\cr\noalign{\hrule}
& 3& 1& 8& 4\cr\noalign{\hrule}
& 1& 9& 9\cr\multispan4\hrulefill\cr
& 6& 5& 2\cr\multispan4\hrulefill\cr
& 5\cr\multispan2\hrulefill\cr
}
}
=\begin{matrix} 3& 1& 8& 4\\ 1& 9& 9\\ 6& 5& 2\\ 5\end{matrix}
\\
\vspace{\medskipamount}
T(1,3)=8\\
\vspace{\medskipamount}
\text{(b) A tableau of shape $(4,3,3,1)$}
\endgathered
$$
\end{figure}
\begin{eg}
\label{eg:cov}
The Young diagram of $\lam=(4,3,3,1)$, shown in
Figure~\ref{fig:ydtab}(a), 
has 3 corners and 4 cocorners, corresponding to the following
covering relations.

$$
\vtop{\packed
\def\pad{\hphantom{,0}}\def\vprop{\vphantom{\dfrac11}}
\halign{\vrule$#\vprop$&&\quad\hfil$#$\hfil\quad\vrule\cr\noalign{\hrule}
&\text{corner}&\text{covering relation}\cr\noalign{\hrule}
&(1,4)&\lam\coverrow{1}(3,3,3,1)\cr\noalign{\hrule}
&(3,3)&\lam\coverrow{3}(4,3,2,1)\cr\noalign{\hrule}
&(4,1)&\lam\coverrow{4}(4,3,3)\pad\cr\noalign{\hrule}
}
}
\qquad\qquad
\vtop{\offinterlineskip\tabskip=0pt
\def\pad{\hphantom{,0}}\def\vprop{\vphantom{\dfrac11}}
\halign{\vrule$#\vprop$&&\quad\hfil$#$\hfil\quad\vrule\cr\noalign{\hrule}
&\text{cocorner}&\text{covering relation}\cr\noalign{\hrule}
&(1,5)&\lam\coveredrow{1}(5,3,3,1)\pad\cr\noalign{\hrule}
&(2,4)&\lam\coveredrow{2}(4,4,3,1)\pad\cr\noalign{\hrule}
&(4,2)&\lam\coveredrow{4}(4,3,3,2)\pad\cr\noalign{\hrule}
&(5,1)&\lam\coveredrow{5}(4,3,3,1,1)\cr\noalign{\hrule}
}
}
$$
\end{eg}

\subsection{Tableaux}\label{subsec:tab}
There are many different conventions
for defining ``tableaux''.
In this article, a \em{tableau} of shape $\lambda$, with $\lambda\in\Par$,
formally means an arbitrary map from $D_{\lambda}$ to a fixed set $\Gam$,
which we call the \em{alphabet}, and whose elements are the \em{letters}.
A tableau $T$ of shape $\lambda$ is visualized as the same cluster of
square boxes as $D_\lam$
with each box $( i, j )$ containing the value $T( i, j )$ of the map $T$
at $( i, j )$.
Thus $T( i, j )$ is also called
the \em{entry} or \em{content} of the square $( i, j )$.
(See Figure~\ref{fig:ydtab}(b).)
The shape of a tableau $T$ will sometimes be denoted by $\sh(T)$.
For each $\gam\in\Gam$, let $m_T(\gam)$ denote
the number of occurrences (``multiplicity'') of the letter $\gam$ in
the tableau $T$.

If $\Gam$ is a totally ordered set, a tableau $T$ is called
\em{semistandard} or \em{column strict} if it satisfies
the following two conditions:
\begin{enumerate}
\item $T(i,1)\le T(i,2)\le\cdots\le T(i,\lam_i)$ for $1\le i\le l$,
where $l=l(\lam)$,
\item $T(1,j)<T(2,j)<\cdots<T(\lam'_j,j)$ for $1\le j\le\lam_1$,
where $\lam'_j$ denotes the length of the $j$th column of $\lam$.
\end{enumerate}

Fix a positive integer $n$,
and let $\Gam_n$ denote the totally ordered set
$\{1<\bar 1<2<\bar 2<\cdots<n<\bar n\}$.
A semistandard tableau $T$ of shape $\lam$, with entries from $\Gam_n$,
is called an $\boldsymbol{S}\boldsymbol{p}\boldsymbol{(}\boldsymbol{2}\boldsymbol{n}
\boldsymbol{)}${\bf-tableau}
or an {\bf $\boldsymbol n$-symplectic tableau}
if it satisfies an additional condition, called \em{the symplectic condition}:
\begin{enumerate}
\item[(3)] $T(i,j)\ge i$ for $1\le i\le l$, $1\le j\le\lam_i$.
\end{enumerate}
The sum of the weight monomials of $T$, defined by
$$
x_1^{m_T(1)-m_T(\bar 1)}x_2^{m_T(2)-m_T(\bar 2)}\cdots x_n^{m_T(n)
-m_T(\bar n)},
$$
for all $Sp(2n)$-tableaux $T$ of a given shape $\lam$,
equals the character $\lam_{Sp(2n)}$ of the irreducible representation
of $Sp(2n,\Cpx)$ labeled by $\lambda$
(see \cite{KoikeTerada}\ and a survey in \cite{SundaramIMA}).

\subsection{(Ordinary) row insertion}\label{subsec:rowins}
To define Berele insertion we first need to define ``(ordinary) row
insertion'' in the sense of Schensted and Knuth.  The description we
give here will be somewhat informal; a more formal version can be
found in \cite{Knuth}.

Given a semistandard tableau $T$ of shape $\lambda $ and a letter
$\gam$, we determine
a new tableau denoted by $T\leftarrow\gam$ as follows.
First ``insert''
$\gam$ into the first row of $T$, which means to replace by $\gam$
the leftmost letter $\gam'$ in the first row
which is strictly larger than $\gam$ (if such $\gam'$ exists),
in which case $\gam'$ is said to get ``bumped'' by $\gam$;
or put $\gam$ at the end of the first row if no such $\gam'$ exists.
As long as a letter gets bumped from one row,
we similarly insert that letter into the next row.
At some iteration, the bumped letter will come to rest at the end
of the next row (possibly creating a new row at the bottom).
The resulting object is a semistandard
tableau (which is denoted by $T\leftarrow\gam$),
whose shape covers $\lam$.
An example of this procedure is contained in Example~\ref{eg:berins}
of Berele insertion below.

\subsection{Jeu de taquin (sliding algorithm)}\label{subsec:jeu}
To define Berele insertion we also need the notion of a jeu de taquin
slide, due originally to Sch\"utzenberger.  %%%%%%%%%%%%%%%%%%%%%%%%%%%cite
Define a {\bf punctured shape} to be a pair $(\lam,h)$,
where $\lam$ is a partition and $h\in D_{\lam}$ (called the \em{hole}),
and its diagram to be $D_{\lam}\setminus\{h\}$.
Define a {\bf punctured tableau} of shape $(\lam,h)$ to be a pair $(T,h)$
where $T$ is a map $D_{\lam}\setminus\{h\}\to\Gam$.
It represents a filling of the squares of $\lam$ except for the ``hole'' $h$,
which is left blank.
It is called semistandard if it satisfies the inequalities (1) and (2) given
in the above definition of semistandard tableau,
in which the hole is to be skipped.

A {\bf (backward) slide} is a transformation $\xi: (T,h)\mapsto (T',h')$
between punctured tableaux.
For a fixed $\lam$, it is a bijection from the set of semistandard punctured
tableaux $(T,h)$ such that $h$ is not a corner of $\lam$,
to the set of those with $h\neq(1,1)$.
It is defined as follows.
Compare the contents of the two squares of $T$
that are below and to the right of $h=(i,j)$.
If $T(i+1,j)\leq T(i,j+1)$ (or if $(i,j+1)\not\in D_{\lam}$),
then set $T'(i,j)=T(i+1,j)$, $ h'=(i+1,j)$, and set
$T'$ to be identical to $T$ elsewhere.
Otherwise set $T'(i,j)=T(i,j+1)$, $h'=(i,j+1)$,
and set $T'$ to be identical to $T$ elsewhere.
Informally, we simply slide the smaller of these two letters
(or the one below if they are equal)
into the hole $h$ and make the vacated square the new hole.
$T'$ is again a semistandard punctured tableau.

Given a semistandard punctured tableau $(T,h)$,
one can repeat slides until the hole comes to rest at a corner
of the shape $\lam$.
At this point one can just forget the hole and consider $T'$ to be a
semistandard tableau of shape $\sh( T' ) \setminus h'$.
We use this procedure below.

\subsection{Berele insertion and Berele's correspondence}\label{subsec:ber}
\em{Berele insertion} is an explicitly given bijection from the set of pairs
$(T,\gam)$, where $T$ is an $Sp(2n)$-tableau of a given shape $\lam$,
and $\gam\in\Gam_n$, to the set of $Sp(2n)$-tableau whose shape
either covers $\lam$ or is covered by $\lam$ (in the poset $\Par$).
If the ordinary row insertion of $\gam$ into $T$
yields a valid $Sp(2n)$-tableau,
then it is also the result of the Berele insertion of $\gam$ into $T$
by definition.
In this case the resulting shape covers $\lam$.
On the other hand,
if the result of the row insertion violates condition (3),
then it must be that, for some $k$, a letter $\bar k$
that was in row $k$ in $T$ was bumped by a letter $k$.
Find the earliest such occurrence, and at this point erase both
the $k$ and $\bar k$ that are involved in this bumping,
leaving the position formerly occupied by the $\bar k$ as a hole.
After this, apply the sliding algorithm until the hole moves to a corner,
and then forget the hole.
This is by definition the result of the Berele insertion, and in this case
the resulting shape is covered by $\lam$.
Let $T\Bereleinsert\gam$ denote the result of the Berele insertion of $\gam$
into $T$.

\begin{eg}
\label{eg:berins}
{\bf Example of Berele insertion}
Berele insertion of $\bar 1$ into the following $Sp(2n)$-tableau $T$ proceeds
as follows, producing $T\Bereleinsert\bar1$ at the end.
In the bumping phase, the caption on the arrow means:
$\xrightarrow[\text{bump:}]{\text{insert:}}$.

\begin{align*}
&T=
\begin{matrix} 1& 1& \bar2& \bar2\\ 2& \bar2& 3& \bar4\\ 3& \bar3& 4\\
4& 4& \bar4\\ 5& \bar5\end{matrix}
\xrightarrow[\text{$\bar 2$ at $(1,3)$}]{\text{$\bar1$ into row $1$}}
\begin{matrix} 1& 1& \bar1& \bar2\\ 2& \bar2& 3& \bar4\\ 3& \bar3& 4\\
4& 4& \bar4\\ 5& \bar5\end{matrix}
\xrightarrow[\text{$3$ at $(2,3)$}]{\text{$\bar2$ into row $2$}}
\begin{matrix} 1& 1& \bar1& \bar2\\ 2& \bar2& \bar2& \bar4\\ 3& \bar3& 4\\
4& 4& \bar4\\ 5& \bar5\end{matrix}
\xrightarrow[\text{$\bar3$ at $(3,2)$}]{\text{$3$ into row $3$}}
\\
\intertext{Placing $\bar3$ in row $4$ would cause a violation,
so cancel $3$ and $\bar3$ and proceed to the sliding phase.}
&\begin{matrix} 1& 1& \bar1& \bar2\\ 2& \bar2& \bar2& \bar4\\ 3& & 4\\
4& 4& \bar4\\ 5& \bar5\end{matrix}
\xrightarrow[\text{move lower $4$ up}]{\text{lower $4\le$ right $4$}}
\begin{matrix} 1& 1& \bar1& \bar2\\ 2& \bar2& \bar2& \bar4\\ 3&  4& 4\\
4&  & \bar4\\ 5& \bar5\end{matrix}
\xrightarrow[\text{move $\bar4$ left}]{\text{lower $\bar5>$ right $\bar4$}}
\begin{matrix} 1& 1& \bar1& \bar2\\ 2& \bar2& \bar2& \bar4\\ 3&  4& 4\\
4& \bar4\\ 5& \bar5\end{matrix}
=T\Bereleinsert\bar1.
\end{align*}

\end{eg}

The weighted enumerative identity following from this bijection
represents the decomposition of the tensor product of the irreducible
representation $\lam_{Sp(2n)}$ of $Sp( 2n, \Cpx )$ labeled by $\lam$
and the natural representation.

\em{Berele's correspondence}, as we call it in this article,
is a bijection from the set of words $w=w_1w_2\dots w_f$
in the alphabet $\Gam_n$ of fixed length $f$
to the set of pairs $(P,Q)$,
where $P$ is an $Sp(2n)$-tableau of some shape $\lam$,
and $Q$ is an $\boldsymbol n$-\em{symplectic up-down tableau} of \em{degree} $f$
with initial shape $\emp$ and final shape $\lam$;
namely $Q=(\emp=\kap\th 0,\kap\th 1,\dots,\kap\th f=\lam)$,
$\kap\th i\in\Par$, $l (\kap\th i)\le n$, %%%%%%%%%%%%%%%% \ell-> l (IT)
and for each $i$ either $\kap\th{i-1}\covered\kap\th i$ or $\kap\th{i-1}
\cover\kap\th i$ holds.
(In the literature, $f$ is generally called the \em{length} of $Q$.
In this article, we call it the \em{degree} in order to avoid any association
with the length of each $\kap \th i$.)
If $w$ is such a word, then for $0\le i\le f$ put
$P_i=(\cdots((\emp\Bereleinsert w_1)\Bereleinsert w_2)\Bereleinsert\cdots
)\Bereleinsert w_i$, and let $\kap\th i$ be the shape of $P_i$.
Put $P=P_f$ and $Q=(\kap\th 0,\kap\th 1,\dots,\kap\th f)$.
Then, by definition, Berele's correspondence takes $w$ to this pair $(P,Q)$.
Following the convention for the Robinson-Schensted correspondence,
we call $P$ and $Q$ the (Berele) $\boldsymbol P$\em{-symbol}
and $\boldsymbol Q$\em{-symbol} of $w$
respectively.

\begin{eg}
\label{eg:ber}
Applying Berele insertion to the word
$w=\bar 3 1 \bar 2 \bar{3} 3
\bar{1} 1 2 \bar 3 \bar 1
\bar 2 2 3 \bar 2 \bar 1
2 2 \bar 3 1 2$
yields the following sequence $P_{i}$ of symplectic tableaux:

%%%%%%%%%%%%  IF YOU DON'T LIKE IT, PLEASE REMOVE %%%%%%%%%%%%%%%%
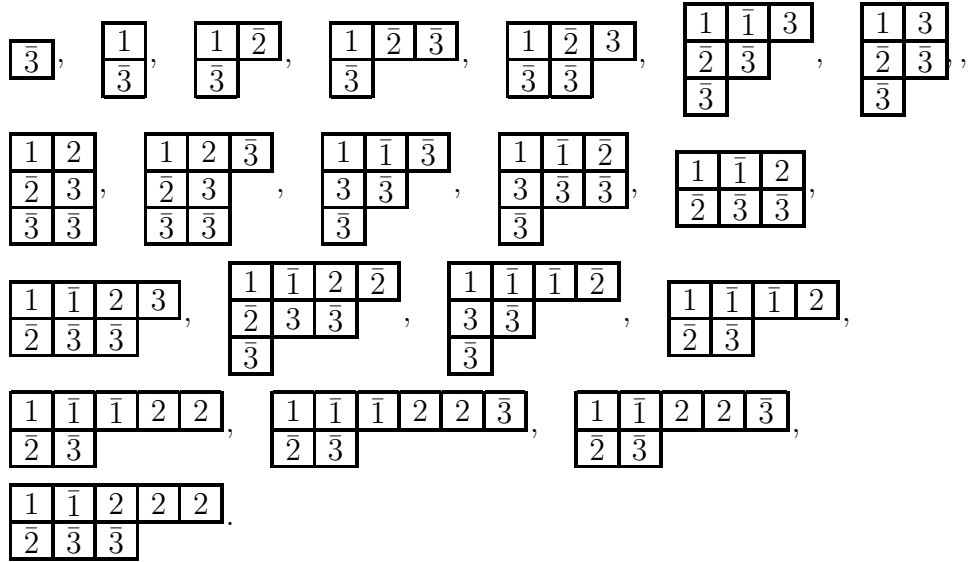
\def\smyoung#1{\Squaresize=15pt
\vcenter{\smallskip\offinterlineskip
\halign{&\vsquare{##}\cr #1}}}
%%%%%%%%%%%%%%%%%%%%%%%%%%%%%%%%%%%%%%%%%%%%%%%%%%%%%%%%%%%%%%%%%%
\begin{align*}
&
\smyoung{\bar{3}\cr }\,,\quad
\smyoung{1\cr
\bar{3}\cr }\,,\quad
\smyoung{1&\bar{2}\cr
\bar{3}\cr }\,,\quad
\smyoung{1&\bar{2}&\bar{3}\cr
\bar{3}\cr }\,,\quad
\smyoung{1&\bar{2}&3\cr
\bar{3}&\bar{3}\cr }\,,\quad
%five
\smyoung{1&\bar{1}&3\cr
\bar{2}&\bar{3}\cr
\bar{3}\cr }\,,\quad
\smyoung{1&3\cr
\bar{2}&\bar{3}\cr
\bar{3}\cr },,\quad
\\
%linebreak
&
\smyoung{1&2\cr
\bar{2}&3\cr
\bar{3}&\bar{3}\cr }\,,\quad
\smyoung{1&2&\bar{3}\cr
\bar{2}&3\cr
\bar{3}&\bar{3}\cr }\,,\quad
\smyoung{1&\bar{1}&\bar{3}\cr
3&\bar{3}\cr
\bar{3}\cr }\,,\quad
%ten
\smyoung{1&\bar{1}&\bar{2}\cr
3&\bar{3}&\bar{3}\cr
\bar{3}\cr }\,,\quad
\smyoung{1&\bar{1}&{2}\cr
\bar{2}&\bar{3}&\bar{3}\cr }\,,\quad
\\
%linebreak
&
\smyoung{1&\bar{1}&{2}&3\cr
\bar{2}&\bar{3}&\bar{3}\cr }\,,\quad
\smyoung{1&\bar{1}&{2}&\bar{2}\cr
\bar{2}&3&\bar{3}\cr
\bar{3}\cr }\,,\quad
\smyoung{1&\bar{1}&\bar{1}&\bar{2}\cr
3&\bar{3}\cr
\bar{3}\cr }\,,\quad
\smyoung{1&\bar{1}&\bar{1}&{2}\cr
\bar{2}&\bar{3}\cr}\,,\quad
\\
%linebreak
&
\smyoung{1&\bar{1}&\bar{1}&{2}&2\cr
\bar{2}&\bar{3}\cr}\,,\quad
\smyoung{1&\bar{1}&\bar{1}&{2}&2&\bar{3}\cr
\bar{2}&\bar{3}\cr}\,,\quad
\smyoung{1&\bar{1}&{2}&2&\bar{3}\cr
\bar{2}&\bar{3}\cr}\,,\quad
\\
&
\smyoung{1&\bar{1}&{2}&2&2\cr
\bar{2}&\bar{3}&\bar{3}\cr}\,.
\end{align*}
Berele's correspondence takes the word $w$ to the pair $(P,Q)$, where
$P$ is the tableau $
%%%%%%%%%%%%  IF YOU DON'T LIKE IT, PLEASE REMOVE %%%%%%%%%%%%%%%%
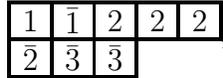
\def\smyoung#1{\Squaresize=15pt
\vcenter{\smallskip\offinterlineskip
\halign{&\vsquare{##}\cr #1}}}
%%%%%%%%%%%%%%%%%%%%%%%%%%%%%%%%%%%%%%%%%%%%%%%%%%%%%%%%%%%%%%%%%%
\smyoung{1&\bar{1}&{2}&2&2\cr
\bar{2}&\bar{3}&\bar{3}\cr}$, and $Q$ is the following:
$$
(0,1,11,21,31,32,321,221,222,322,321,331,33,43,431,421,42,52,62,52,53).
$$
\end{eg}

The enumerative identity following from the whole Berele correspondence
represents the decomposition of the $f$-fold tensor product of the
natural representation of $Sp(2n)$. %%%%%%%%%%%%%%%%%%%%%%%%%%%%%%%%%%
%%%%%%%%%%%%%%%% To be put into INTRODUCTION!!!
For more information about related matters,
we refer the interested reader to \cite{SundaramIMA},
which includes a nice survey and an interesting connection between
up-down tableaux and standard tableaux.
%%%%%%%%%%%%%%%%%%%%%%%%%%%%%%%%%%%%%%%%%%%%%%%

\subsection{Standardization of R-S Correspondence}\label{subsec:strs}

In order to give a pictorial interpretation,
we introduce a ``standardized'' version of the Berele correspondence.
Before discussing standardization of the Berele correspondence,
let us include a brief summary of the situation for the R-S
correspondence.

In Schensted's original paper, he is interested in enumerating the
number of permutations with a certain fixed length of longest increasing
subsequence.  To generalize this to words with repeated entries, in Part
II of his paper, he mapped such a word to a permutation (in a natural
way), applied his insertion algorithm to this permutation, and then
mapped the resulting entries of the $P$ symbol back. However, he
provides neither a formal definition of \em{standardization} nor a
proof that it commutes with insertion.   

Sch\"utzenberger \cite{SchutzenbergerLNM} not only defined
standardization of semistandard tableaux, but also showed the validity
of a commutative diagram like Figure~\ref{fig:strs} below
for semistandard tableaux by using the sliding algorithm to explicate
Schensted insertion.  To generalize to the symplectic case, we prefer to
have a lemma and a proof that directly compare semistandard and standardized
insertion.  Standardization of shifted tableaux was given by B. Sagan
in \cite{Sagan}. Our approach is most similar to his.  
%% his treatment is most similar to ours.

\begin{defn} \label{def:stand}
Let $w=w_1w_2\dots w_f$ be a word in the alphabet $\Gam = [ 1, n ]$ of
length $f$.  
The R-S correspondence for multiset permutations takes $w$ to a pair $( P, Q )$
where $P$ is a semistandard tableau of the same weight as $w$,
and $Q$ is a standard tableau of the same shape as $P$.
The \em{standardization} $\tilde{w}=\tilde{w_{1}}\tilde{w_{2}}\cdots
\tilde{w_{f}}$ is the word obtained from $w$ by replacing, for each
$\gamma \in \Gam $, the occurrences of the letter $\gamma $ in $w$ by
the symbols $\gamma _{1}, \gamma _{2},\ldots,\gamma _{m_{w}(\gamma )}$
from left to right, where $m_{w(\gamma )}$ is the number of such
occurrences.  Let $\tilde{\Gam }_{w}$ denote the totally ordered set
$1_{1}<1_{2}<\cdots <1_{m_{w}(1)}<2_{1}<2_{2}<\cdots <2_{m_{w}(2)}<\cdots
< n_{1}<n_{2}<\cdots <n_{m_{w}(n)}$.  
By the \em{standardization} $\tilde P$ of $P$ we mean a standard tableau
with entries from $\tilde \Gam_w$, instead of $[ 1, f ]$, obtained from $P$
by replacing
the occurrences of each letter $\gam$ in $P$ (which form a horizontal strip)
by $\gam_1$, $\gam_2$, \dots, $\gam_{ m_P( \gam ) }$ from left to right.
\end{defn}
We now give a direct proof that standardization commutes with Schensted
insertion that we will later generalize to the symplectic case.  

%(Schensted \cite{Schensted} did not write any proof of this,
%and Sch\"utzenberger \cite{SchutzenbergerLNM} proved this
%proved this in terms of the sliding algorithm rather than insertion:
%unfortunately we have been unable to find a reference for a proof
%directly comparing the insertion processes,
%so let us write it down in a form we use below).

\begin{lemma}
\label{lem:LemRSStand}
Let $w = w_1 w_2 \cdots w_f$ be a word in $\Gam = [ 1, n ]$ of length $f$,
and let $\tilde w = \tilde w_1 \tilde w_2 \cdots \tilde w_f$ be its
standardization.
Let $P \th i$ denote the tableau obtained by inserting
$w_1$, $w_2$, \dots, $w_i$ into the empty tableau,
and let $\tilde P \th i$ be the standardization of $P \th i$.
Then, for each $i$, the insertion of $\tilde w_i$ into $\tilde P \th{ i - 1 }$
follows exactly the same route as that of $w_i$ into $P \th{ i - 1 }$,
and the resulting tableau coincides with $\tilde P \th i$.
\end{lemma}

\begin{proof}
We compare the insertion of $\tilde w_i$ into $\tilde P \th{ i - 1 }$
(the standardized case) with that of $w_i$ into $P \th{ i - 1 }$
(the unstandardized case).
We will show the following claim holds row by row along with the insertion;
then the Lemma follows immediately.

We define one technical notion.
Let $T$ be a semistandard tableau of shape $\lam$ in the alphabet $\Gam =
[ 1, n ]$, and let $k \in \Gam$ be a letter.
For $r \ge 1$, let $c_{+}( k, r )$ and $c_{-}( k, r )$ be defined by:
\begin{align*}
c_-( k, r ) & = \max \{ 0 \} \cup \{\, j \mid T( r, j ) \le k \,\}, \\
c_+( k, r ) & = \begin{cases} \min \{ \lam_{ r - 1 } + 1 \} \cup \{\, j
\mid T( r - 1, j ) \ge k \,\} & \text{if $r \ge 2$,} \\
\infty & \text{if $r = 1$.} \end{cases}
\end{align*}
Roughly $c_-( k, r ) $ gives the rightmost column of row $r$ 
containing entries $\leq k$, while $c_+( k, r )$ gives the
leftmost column of the {\it previous} row with entries $\geq k$.  
The semistandardness guarantees that $c_+( k, r ) - c_-( k, r ) \ge 1$.
Let us say that $T$ has a \em{{\boldmath$k$}-gap} between rows $r-1$ and
$r$ if in fact 
$c_+( k, r ) - c_-( k, r ) \geq 2$.

%%Do we want a claim environment?
{\bf Claim.}
Suppose the bumping is about to reach row $r$ in both the standardized
and unstandardized cases.
Suppose the intermediate tableau $\tilde T$ of the standard case
at this point is obtained from the intermediate tableau $T$
of the unstandardized case by modified standardization in the following
sense (inductive hypothesis).
\begin{enumerate}
\item Let $k \in [ 1, n ]$ be the letter bumped from row $r - 1$
\rom(or $k = w_i$ if $r = 1$\rom) in the unstandardized case.
Then the letter bumped from row $r - 1$ in the standardized case
is $k_s$ with some index $s$.
\item For each $k' \neq k$, the $k'$ with various indices in $\tilde T$
occupy the same positions as the $k'$ in $T$, which form a horizontal strip,
and their indices increase from left to right.
\item The $k$ with various indices in $\tilde T$ occupy
the same positions as the $k$ in $T$,
which form a horizontal strip, and are indexed as follows.
The $k$ in rows $r$ and below are indexed from $1$ to $s - 1$
from left to right, and those in rows $r - 1$ and above are indexed
from left to right starting with $s + 1$.
This together with (2) assures that $T$ is semistandard,
and we further assume that $T$ has a $k$-gap
between rows $r-1$ and $r$.
\end{enumerate}
Then the following hold.
\begin{enumerate}
\item[(a)] The insertion terminates at row $r$ in the standardized case
if and only if it terminates at row $r$ in the unstandardized case.
\item[(b)] If the bumping continues, then the bumping at row $r$ occurs
at the same position for both cases, and the intermediate tableaux
after bumping from row $r$ satisfy \rom(1\rom)--\rom(3\rom) above
with $r$ replaced by $r + 1$.
\end{enumerate}
\smallskip
%%End of claim
First note that the insertion terminates at row $r$ in the unstandardized case
if and only if all entries in row $r$ of $T$ are at most $k$.
Since $k$ with indices greater than $s$ cannot exist in row $r$
by assumption, this is equivalent to saying that all entries in row $r$
of $\tilde T$ 
are less than $k_s$, precisely in which case the insertion terminates here
in the standard case.
Hence (a).

Now suppose the bumping continues.
Let the conditions (1)--(3) claimed in (b)
for the new intermediate tableaux be written as (1)$^*$--%
(3)$^*$, as opposed to the conditions (1)--%
(3) for $T$ and $\tilde T$ in the assumption.
Let $k\sp*$ be the letter bumped by $k$
from row $r$ of $T$ in the unstandardized case.
It is the leftmost letter greater than $k$ in this row.
Since again by assumption $\tilde T$ contains no $k$
with indices greater than $s$ in row $r$ of $\tilde T$,
the bumped letter in the standardized case is also a $k\sp*$,
more precisely $k\sp*$ with the smallest index in this row.
Since the indices of $k\sp*$ increase from left to right in a row
by assumption,
it is also the leftmost $k\sp*$ in this row of $\tilde T$.
So the bumping occurs at the same position in both cases,
and (1)$^*$ also follows.
Let $t$ be the index of this $k\sp*$.
The only difference to the $k\sp*$ in $\tilde T$ (resp.\ $T$)
caused by this bumping is that it loses $k\sp*_t$
(resp.\ the $k\sp*$ in the same position), so that (3)$^*$ follows
from the assumption (2) applied to $k' = k\sp*$%
%(modulo (2)$^*$, which we show in the following paragraph)%
.

Since no letters other than $k$ or $k\sp*$ move
during this bumping in row $r$,
(2)$^*$ for those other letters follows from the assumption
(2).
Now let us concentrate on the letters $k$.
We know that the letters $k$ form a horizontal strip in $T$
and $\tilde T$,
and the only change caused during this step was an addition of $k_s$
into row $r$, immediately to the right of column $c_-( r, k )$.
Because of the $k$-gap between rows $r-1$ and $r$ in $T$,
this is still to the
left of the column $c_+( r, k )$, so that the letters $k$
still form a horizontal strip after this addition.
All other $k$ in row $r$ have smaller indices,
and so do those in rows below.
Those in row $r - 1$ and higher have indices larger than $s$ by assumption
(3), so (2)$^*$ also holds for $k$.
\end{proof}

This lemma shows the validity of the commutative diagram in
Figure~\ref{fig:strs}.  
\begin{figure}
\caption{Standardization commutes with ordinary R-S correspondence}
\label{fig:strs}
\begin{equation*}
\begin{CD}
%w \xrightarrow[]{\text{R-S}} (P,Q)\\
w @>\text{R-S}>> (P,Q)\\
@V\text{standardization}VV @VV{
\begin{aligned}
& \ssize \tilde P = \text{standardization of $P$} \\
& \ssize \tilde Q = Q
\end{aligned}
}V \\
%\tilde w \xrightarrow[\text{R-S}]{} (\tilde P,\tilde Q)
\tilde w @>>\text{R-S}> (\tilde P,\tilde Q)
\end{CD}
\end{equation*}
\end{figure}

\subsection{Standardized Berele's correspondence}\label{subsec:stber}

Let $w$ be a word in $\Gam_n=\{1<\bar 1<2<\bar 2<\cdots<n<\bar n\}$.
Let $\tilde{\Gam }_{w}$ be defined as in Def.~\ref{def:stand}, but
with $\Gam=[1,n] $ replaced by $\Gam_{n} $, and $\ord\:\tilde\Gam_w\to[1,f]$ be the unique order-preserving bijection.   

Now we can define \em{standardized Berele's correspondence} for
standardized words.
In \em{standardized Berele insertion},
all our bumping and slides occur according to
the usual rules (though in this case all letters are distinct).
Violations of the symplectic condition are determined by
ignoring the subscripts of the symbols $\gam_t$.
The only point that needs careful consideration
is the handling of cancellation.

A violation occurs exactly when $k_s$ \label{1shiftpath}
($1\le k\le n$, $s$ being any index) tries to bump $\bar k_t$
($t$ again being any index) out of the $k$th row,
say from the cell $(k,c)$.
First put the $k_s$ at $(k,c)$, which action does not yet cause a violation,
and throw the $\bar k_t$ away instead of inserting it into the next row.
Note that now the tableau contains $k$'s in cells $(k,1)$--$(k,c)$,
because letters smaller than $k$ cannot appear in this row
due to the symplectic condition,
and $\bar k_t$ must have been the smallest $\bar k$ in this row
(therefore the leftmost) in order to be bumped.
Now remove the $k$ in the cell $(k,1)$, which 
is the smallest $k$ in the tableau, and is $k_s$ if and only if $c=1$, 
and move this hole by the sliding algorithm.\footnote{The authors are
grateful to K.~Koike for raising the question of which subscripted $k$
should be considered cancelled by $\bar{k}$ in this situation.}
Note that if $c>1$, then the hole continues to move to the right
up to $(k,c)$.
Therefore, if we discard the subscripts, this amounts to the same thing
as cancelling $k_s$ with $\bar k_t$ and
making $(k,c)$ the initial hole for the sliding algorithm.
It turns out that removing the smallest $k$ enables
easier consistent handling.

This insertion will be denoted by $\StandBereleinsert$.

\begin{eg}
\label{eg:cancel}
{\bf(Example of cancellation)}\enspace
Suppose $3_2$ has been bumped from the 2nd row, and is about to be inserted
into the 3rd row.
In this example it bumps $\bar3_1$, which cannot be placed in row 4.
So $\bar3_1$ is removed, and in cancellation the smallest $3$,
which in this case is $3_1$, is removed.
Note that all letters in row 4 must be $\ge 4$,
assuring that the sliding proceeds sideways until the hole comes
to the position previously occupied by $\bar3_1$.
Compare this to the $3$-$\bar3$ cancellation between the third and
fourth tableaux in Example~\ref{eg:berins}.

\begin{align*}
&
\begin{matrix} 1_1& 1_2& \bar1_1& \bar2_3\\ 2_1& \bar2_1& \bar2_2& \bar4_2\\
3_1& \bar3_1& 4_3\\ 4_1& 4_2& \bar4_1\\ 5_1& \bar5_1\end{matrix}
\xrightarrow[\text{remove $\bar 3_1$}]{\text{insert $3_2$ into row 3}}
\begin{matrix} 1_1& 1_2& \bar1_1& \bar2_3\\ 2_1& \bar2_1& \bar2_2& \bar4_2\\
3_1& 3_2& 4_3\\ 4_1& 4_2& \bar4_1\\ 5_1& \bar5_1\end{matrix}
\xrightarrow[]{\text{remove $3_1$}}
\begin{matrix} 1_1& 1_2& \bar1_1& \bar2_3\\ 2_1& \bar2_1& \bar2_2& \bar4_2\\
& 3_2& 4_3\\ 4_1& 4_2& \bar4_1\\ 5_1& \bar5_1\end{matrix}
\\
&
\xrightarrow[\text{move $3_2$}]{{4_1>3_2}}
\begin{matrix} 1_1& 1_2& \bar1_1& \bar2_3\\ 2_1& \bar2_1& \bar2_2& \bar4_2\\
3_2&  & 4_3\\ 4_1& 4_2& \bar4_1\\ 5_1& \bar5_1\end{matrix}
\end{align*}
\end{eg}

Given a standardized word $\tilde w=\tilde w_1\tilde w_2\cdots\tilde w_f$,
put $\tilde P_i=(\cdots((\emp\StandBereleinsert\tilde w_1)
\StandBereleinsert\tilde w_2)\StandBereleinsert\cdots)\StandBereleinsert
\tilde w_i$, and $\kap\th i=\sh(\tilde P_i)$ for $0\le i\le f$.
Put $\tilde P=\tilde P_f$ and $\tilde Q=(\kap\th 0,\kap\th 1,%\allowmathbreak
\dots,\kap\th f)$.
Standardized Berele's correspondence takes $\tilde w$ to the pair
$(\tilde P,\tilde Q$),
and the terms $P$-symbol and $Q$-symbol will be used
as in the original Berele's correspondence.

Then it is possible to define standardization of $Sp(2n)$-tableaux,
in a limited sense:

\begin{lemma}
\label{LemStand}

Let $w=w_1w_2\cdots w_f$ be a word of length $f$ in the alphabet $\Gam_n$,
and $\tilde w=\tilde w_1\tilde w_2\cdots\tilde w_f$ be its standardization.
Suppose $w$ corresponds to $(P,Q)$ by Berele's correspondence,
and $\tilde w$ to $(\tilde P,\tilde Q)$ by standardized
Berele's correspondence.
For each $\gam\in\Gam_n$, let $c_w(\gam)$ be the number of letters $\gam$
removed in cancellation during the process of Berele's correspondence
applied to the word $w$.
Note that $c_w(k)=c_w(\bar k)$ for any $k\in[1,n]$,
and that $\sum_{\gam\in\Gam_n}c_w(\gam)=2\sum_{k=1}^n c_w(k)=f-|\sh(P)|$.

Then we have $\tilde Q=Q$, and $\tilde P$ is obtained from $P$
by replacing, for each $\gam\in\Gam_n$, the occurrences of the letter $\gam$
in $P$ by the letters $\gam_{c_w(\gam)+1}$, $\gam_{c_w(\gam)+2}$, \dots,
$\gam_{m_w(\gam)}$ in this order from left to right.
\rom(Note that this makes sense since they form a horizontal strip in $P$.\rom)
\end{lemma}

\begin{proof}
One proves this by induction on $f$,
starting with the trivial case where $f=0$.
Now suppose $f>0$; then by the induction hypothesis
the lemma holds for $\bar w=w_1w_2\cdots w_{f-1}$.
The standardization of $\bar w$ is $\tilde w_1\tilde w_2\cdots\tilde w_{f-1}$.
For simplicity put $\bar P=P(\bar w)$ and $\tilde{\bar P}=
\tilde P(\tilde{\bar w})$.
A similar result for ordinary row insertion (see Lemma~\ref{LemStand}) 
assures that the bumping phase of $\tilde{\bar P}\StandBereleinsert\tilde w_f$
proceeds along exactly the same route
as that of $\bar P\Bereleinsert w_f$;
moreover, if cancellation is not involved, the lemma holds for $w$ as well,
and if cancellation is involved, it occurs at exactly the same timing
as it occurs in $\bar P\Bereleinsert w_f$.
In the latter case suppose the cancellation is for the pair $k$-$\bar k$.
Then the offending $\bar k$ that is bumped and removed
must be the smallest (leftmost) $\bar k$ in $\tilde{\bar P}$,
since it is the leftmost $\bar k$ in the $k$th row due to the rule of bumping,
and because of semistandardness any $\bar k$ to the left of this $\bar k$
must be in a row below, which is prohibited by the symplectic condition.
The next instruction by the standardized Berele insertion
is to remove the smallest $k$, whose subscript matches that of the $\bar k$
just removed.
So the requirement for the subscripts of $k$'s and $\bar k$'s
remaining in $\tilde P$ is fulfilled.
As stated above, the sliding in the standardized version
continues to move to the right
until it moves the $k$ that has just bumped the $\bar k$,
and after this point the sliding follows exactly the same route
as in the original version.
The left-to-right increasing order of the subscripts of each letter
is preserved under each step of sliding.
Therefore, in this case also, the lemma holds for $w$.
\end{proof}

\begin{eg}
\label{eg:stand}
For $w$ in Example~\ref{eg:ber}, we have
$$
\tilde w=\bar 3_1      1_1 \bar 2_1 \bar 3_2      3_1
\bar 1_1      1_2      2_1 \bar 3_3 \bar 1_2
\bar 2_2      2_2      3_2 \bar 2_3 \bar 1_3
2_3      2_4 \bar 3_4      1_3      2_5.
$$

The set $\tilde\Gam_w$ is the set of subscripted letters in the upper row
of the following table.
The table describes the ordinal function for this $w$.

$$
\vbox{\packed
\halign{\vrule#\vphantom{$\dfrac11$}&&\hfil$\;#\;$\hfil\vrule\cr
\noalign{\hrule}
&\gam_t&1_1&1_2&1_3&\bar 1_1&\bar 1_2&\bar 1_3&2_1&2_2&2_3&2_4&2_5&
\bar 2_1&\bar 2_2&\bar 2_3&3_1&3_2&\bar 3_1&\bar 3_2&\bar 3_3&\bar 3_4\cr
\noalign{\hrule}
&\ord(\gam_t)&1&2&3&4&5&6&7&8&9&10&11&12&13&14&15&16&17&18&19&20\cr
\noalign{\hrule}
}
}
$$

The table in Figure~\ref{fig:sbi} describes the proceedure
of standardized Berele's correspondence for the word $\tilde w$
in a step-by-step manner.
The whole procedure starts with an empty tableau,
which is omitted from the table.
Each line describes Berele insertion of one letter.
The field (A) lists the letters involved in the bumping phase,
excluding the inserted letter $\tilde w_i$, which is written
in the leftmost field.
If the symplectic condition is violated,
the offending letter, which is at the end of the list and is underlined,
gets removed and sliding starts.
The field (B) lists the letters involved in the sliding phase, if any.
The first letter, also underlined, gets removed in cancellation,
and the rest get moved.

\begin{figure}
\caption{A detailed example of standardized Berele insertion}
\label{fig:sbi}
$$
\vbox{\offinterlineskip
\hrule
\def\inter{height1pt&\omit&&\omit&&\omit&&\omit&&\omit&\cr
\noalign{\hrule}
height1pt&\omit&&\omit&&\omit&&\omit&&\omit&\cr}
\def\ul{\underline}
\def\tab{\smallmatrix}\def\etab{\endsmallmatrix}
\halign{&\vrule#&\strut\quad$#$\hfil\quad\cr
height2pt&\omit&&\omit&&\omit&&\omit&&\omit&\cr
&\hfil\tilde w_i&&\hfil\text{(A)}&&\hfil\text{(B)}&&\hfil\tilde P_i&
&\hfil\sh(\tilde P_i)&\cr
height2pt&\omit&&\omit&&\omit&&\omit&&\omit&\cr
\noalign{\hrule}
height1pt&\omit&&\omit&&\omit&&\omit&&\omit&\cr
&\bar 3_1&& && &&\tab \bar 3_1\etab&&(1)&\cr\inter
&1_1&&\bar 3_1&& &&\tab 1_1\\ \bar 3_1\etab&&(1,1)&\cr\inter
&\bar 2_1&& && &&\tab 1_1&\bar 2_1\\ \bar 3_1\etab&&(2,1)&\cr\inter
&\bar 3_2&& && &&\tab 1_1&\bar 2_1&\bar 3_2\\ \bar 3_1\etab&&(3,1)&
\cr\inter
&3_1&&\bar 3_2&& &&\tab 1_1&\bar 2_1&3_1\\ \bar 3_1&\bar 3_2\etab&&
(3,2)&\cr\inter
&\bar 1_1&&\bar 2_1,\bar 3_1&& &&\tab 1_1&\bar 1_1&3_1\\ \bar 2_1&\bar 3_2
\\ \bar 3_1\etab&&(3,2,1)&\cr\inter
&1_2&&\ul{\bar 1_1}&&\ul{1_1},1_2,3_1&&\tab 1_2&3_1\\ \bar 2_1&\bar 3_2\\
\bar 3_1\etab&&(2,2,1)&\cr\inter
&2_1&&3_1,\bar 3_2&& &&\tab 1_2&2_1\\ \bar 2_1&3_1\\ \bar 3_1&\bar 3_2
\etab&&(2,2,2)&\cr\inter
&\bar 3_3&& && &&\tab 1_2&2_1&\bar 3_3\\ \bar 2_1&3_1\\ \bar 3_1&\bar 3_2
\etab&&(3,2,2)&\cr\inter
&\bar 1_2&&2_1,\ul{\bar 2_1}&&\ul{2_1},3_1,\bar 3_2&&\tab 1_2&\bar 1_2&
\bar 3_3\\ 3_1&\bar 3_2\\ \bar 3_1\etab&&(3,2,1)&\cr\inter
&\bar 2_2&&\bar 3_3&& &&\tab 1_2&\bar 1_2&\bar 2_2\\ 3_1&\bar 3_2&\bar 3_3\\
\bar 3_1\etab&&(3,3,1)&\cr\inter
&2_2&&\bar 2_2,3_1,\ul{\bar 3_1}&&\ul{3_1}&&\tab 1_2&\bar 1_2&2_2\\ \bar 2_2&
\bar 3_2&\bar 3_3\etab&&(3,3)&\cr\inter
&3_2&& && &&\tab 1_2&\bar 1_2&2_2&3_2\\ \bar 2_2&\bar 3_2&\bar 3_3\etab&&
(4,3)&\cr\inter
&\bar 2_3&&3_2,\bar 3_2&& &&\tab 1_2&\bar 1_2&2_2&\bar 2_3\\ \bar 2_2&3_2&
\bar 3_3\\ \bar 3_2\etab&&(4,3,1)&\cr\inter
&\bar 1_3&&2_2,\ul{\bar 2_2}&&\ul{2_2},3_2,\bar 3_3&&\tab 1_2&\bar 1_2&
\bar 1_3&\bar 2_3\\ 3_2&\bar 3_3\\ \bar 3_2\etab&&(4,2,1)&\cr\inter
&2_3&&\bar 2_3,3_2,\ul{\bar 3_2}&&\ul{3_2}&&\tab 1_2&\bar 1_2&\bar 1_3&2_3\\
\bar 2_3&\bar 3_3\etab&&(4,2)&\cr\inter
&2_4&& && &&\tab 1_2&\bar 1_2&\bar 1_3&2_3&2_4\\ \bar 2_3&\bar 3_3\etab&&
(5,2)&\cr\inter
&\bar 3_4&& && &&\tab 1_2&\bar 1_2&\bar 1_3&2_3&2_4&\bar 3_4\\ \bar 2_3&
\bar 3_3\etab&&(6,2)&\cr\inter
&1_3&&\ul{\bar 1_2}&&\ul{1_2},1_3,\bar 1_3,2_3,2_4,\bar 3_4&&\tab 1_3&\bar 1_3&2_3&
2_4&\bar 3_4\\ \bar 2_3&\bar 3_3\etab&&(5,2)&\cr\inter
&2_5&&\bar 3_4&& &&\tab 1_3&\bar 1_3&2_3&2_4&2_5\\ \bar 2_3&\bar 3_3&
\bar 3_4\etab&&(5,3)&\cr
height1pt&\omit&&\omit&&\omit&&\omit&&\omit&\cr}
\hrule}
$$
\end{figure}

\end{eg}

\begin{rem}
(1) Let the $2n$-tuple of integers
$(m_w(\gam))_{\gam\in\Gam_n}$ be called 
the \em{literal weight} of $w$.
If we fix the literal weight of $w$ for example,
then Lemma~\ref{LemStand} gives the operation $(P,Q)\mapsto(\tilde P,\tilde Q)$
which makes the following diagram commute,
since we can determine the $c_w(\gam)$
by comparing the given $m_w(\gam)$ and the number of symbols $\gam$
remaining in $P$.

\begin{figure}
\caption{Standardization commutes with Berele's correspondence}
\label{fig:stber}
\begin{equation*}
\begin{CD}
%w \xrightarrow[]{\phantom{xxxxx}\text{Berele}\phantom{xxxxx}} (P,Q)\\
w @>\phantom{xxxxx}\text{Berele}\phantom{xxxxx}>> (P,Q) \\
@V\text{standardization}VV @VV\text{Lemma~\ref{LemStand}}V \\
%\tilde w \xrightarrow[\text{standardized Berele}]{} (\tilde P,\tilde Q)
\tilde w @>>\text{standardized Berele}> (\tilde P,\tilde Q)
\end{CD}
\end{equation*}
\end{figure}

It should be possible to determine the $c_w(\gam)$ from the pair $(P,Q)$ alone.
The sum $\sum_{k=1}^n c_w(k)$ equals the number of shrinks
($\kap\th{i-1}\cover\kap\th{i}$) in the sequence $Q$,
but the problem is how to ``distribute'' this sum among various $k$'s.
In principle it is possible since we can run Berele's correspondence
backwards to find $w$.
Unfortunately we have not yet found a direct method,
which would be extremely useful.
The same problem occurs when we try to reverse our pictorial presentation,
as discussed in \S3.

(2) Unlike the ordinary case, we cannot construct standardized Berele's
correspondence (i.e., the map along the bottom row of 
Figure~\ref{fig:stber} 
by replacing $\tilde w$ by a permutation
$$\begin{pmatrix} 1&2&\cdots&f \\ \ord(\tilde w_1)&\ord(\tilde w_2)&\cdots&
\ord(\tilde w_f) \end{pmatrix}\,.$$ For it is essential to know what each letter was before
standardization in order to detect violations of the symplectic condition
correctly.

We could however regard $\tilde w$ as a weighted permutation, as in
\cite{SaganStanley}.  
\end{rem}

\section{Berele's Correspondence by Local Rules}\label{sec:berlr}

\subsection{Picture of $\boldsymbol w$}\label{subsec:pic}
Next we explain our pictorial approach, which is a two-dimensional
presentation of Berele's algorithm based on a modified set of local
rules in the spirit of Fomin.
We draw an $f\times f$ lattice as in Example~\ref{eg:picw}.   
We employ the matrix coordinate system, and the vertices are labelled
 $(i,j)$ with $0\le i\le f$, $0\le j\le f$.
In this section, we use the letters $A$, $B$, $C$, and $D$
 to denote lattice points.
When we use these names together,
 we generally assume that they have coordinates
 $A=(i-1,j-1)$, $B=(i-1,j)$, $C=(i,j-1)$, and $D=(i,j)$ respectively,
 for some $i$ and $j$.
\begin{figure}
\caption{A cell in the pictorial grid}
\label{fig:abcd}
$$
\gathered
\cell{A}{}{B}{}{}{}{C}{}{D} \\
\endgathered
$$
\end{figure}

The square region with vertices $A$, $B$, $C$, and $D$ will be called
the \em{cell} at $(i,j)$. 
For each $\gam\in\Gam_n$, the cells $(i,j)$ with $i\in[\ord(\gam_1),
 \ord(\gam_{m_w(\gam)})]$ will be said to constitute the {\boldmath $\gam$}-\em{stratum}.
We will refer to this partitioning of the lattice
 as its {\bf stratification}.
The \em{picture of} $\boldsymbol w$ is obtained from this stratified grid
 by writing $\cross$
 inside the cells at $(\ord(\tilde w_j),j)$ for $1\le j\le f$.
We say that the $\cross$ at $(\ord(\tilde w_j),j)$ \em{represents}
 the letter $\tilde w_j$, and define the \em{contents} of
$(\ord(\tilde w_j),j)$ to be $\cross $; (the contents of) any cells not
marked with an $\cross $ is said to be empty.

\begin{eg}
\label{eg:picw}

In our Example~\ref{eg:ber}, the $\ord(\tilde w_j)$ are as follows:
\begin{figure}
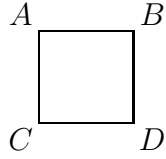

\caption{$\ord(\tilde w_j)$ for Example~\ref{eg:ber}}
\label{fig:ord}
$$
\vbox{\packed
 \halign{\vrule#\vphantom{$\dfrac11$}&&\hfil\;$#$\;\hfil\vrule\cr
  \noalign{\hrule}
  &j&1&2&3&4&5&6&7&8&9&10&11&12&13&14&15&16&17&18&19&20\cr
  \noalign{\hrule}
  &\tilde w_j&\bar 3_1&1_1&\bar 2_1&\bar 3_2&3_1&\bar 1_1&1_2&2_1&\bar 3_3&
   \bar 1_2&\bar 2_2&2_2&3_2&\bar 2_3&\bar 1_3&2_3&2_4&\bar 3_4&1_3&2_5\cr
  \noalign{\hrule}
  &\ord(\tilde w_j)&17&1&12&18&15&4&2&7&19&5&13&8&16&14&6&9&10&20&3&11\cr
  \noalign{\hrule}
 }
}
$$
\end{figure}
The picture of this $w$ is shown in Figure~\ref{fig:picw}.   

On the left edge, the corresponding letters in $\tilde\Gam_w$ are shown.
Thicker horizontal grid lines separate the strata.
On the top are the column numbers of the cells.

\begin{figure}
\caption{The Picture of $w$ for Example~\ref{eg:ber}}
\label{fig:picw}
$$
\vbox{\packed
 \def\hcell#1{\vbox{\hbox to 1.5em{\hss$#1$\hss}}}
 \def\vcell#1{\vbox to 1.5em{\vss\hbox{$#1$}\vss}}
 \def\cell#1{\vbox to 1.5em{\vss\hbox to 1.5em{\hss$#1$\hss}\vss}}
 \def\nr{\omit&\multispan{21}\hrulefill\cr}
 \def\tr{\omit&\multispan{21}\leaders\hrule height1.0pt\hfill\cr}
 \def\x{\cross}
 \halign{\phantom{$#$}&\hskip 0.4pt#&&\hcell{#}\hskip 0.4pt\cr
      1_1&& 1& 2& 3& 4& 5& 6& 7& 8& 9&10&11&12&13&14&15&16&17&18&19&20\cr
 }
 \halign{\vcell{#}&\vrule#&&\cell{#}\vrule\cr
  \tr
      1_1&&  &\x&  &  &  &  &  &  &  &  &  &  &  &  &  &  &  &  &  &  \cr%1
  \nr
      1_2&&  &  &  &  &  &  &\x&  &  &  &  &  &  &  &  &  &  &  &  &  \cr%2
  \nr
      1_3&&  &  &  &  &  &  &  &  &  &  &  &  &  &  &  &  &  &  &\x&  \cr%3
  \tr
  \bar1_1&&  &  &  &  &  &\x&  &  &  &  &  &  &  &  &  &  &  &  &  &  \cr%4
  \nr
  \bar1_2&&  &  &  &  &  &  &  &  &  &\x&  &  &  &  &  &  &  &  &  &  \cr%5
  \nr
  \bar1_3&&  &  &  &  &  &  &  &  &  &  &  &  &  &  &\x&  &  &  &  &  \cr%6
  \tr
      2_1&&  &  &  &  &  &  &  &\x&  &  &  &  &  &  &  &  &  &  &  &  \cr%7
  \nr
      2_2&&  &  &  &  &  &  &  &  &  &  &  &\x&  &  &  &  &  &  &  &  \cr%8
  \nr
      2_3&&  &  &  &  &  &  &  &  &  &  &  &  &  &  &  &\x&  &  &  &  \cr%9
  \nr
%           1  2  3  4  5  6  7  8  9 10 11 12 13 14 15 16 17 18 19 20
      2_4&&  &  &  &  &  &  &  &  &  &  &  &  &  &  &  &  &\x&  &  &  \cr%10
  \nr
      2_5&&  &  &  &  &  &  &  &  &  &  &  &  &  &  &  &  &  &  &  &\x\cr%11
  \tr
  \bar2_1&&  &  &\x&  &  &  &  &  &  &  &  &  &  &  &  &  &  &  &  &  \cr%12
  \nr
  \bar2_2&&  &  &  &  &  &  &  &  &  &  &\x&  &  &  &  &  &  &  &  &  \cr%13
  \nr
  \bar2_3&&  &  &  &  &  &  &  &  &  &  &  &  &  &\x&  &  &  &  &  &  \cr%14
  \tr
      3_1&&  &  &  &  &\x&  &  &  &  &  &  &  &  &  &  &  &  &  &  &  \cr%15
  \nr
      3_2&&  &  &  &  &  &  &  &  &  &  &  &  &\x&  &  &  &  &  &  &  \cr%16
  \tr
  \bar3_1&&\x&  &  &  &  &  &  &  &  &  &  &  &  &  &  &  &  &  &  &  \cr%17
  \nr
  \bar3_2&&  &  &  &\x&  &  &  &  &  &  &  &  &  &  &  &  &  &  &  &  \cr%18
  \nr
  \bar3_3&&  &  &  &  &  &  &  &  &\x&  &  &  &  &  &  &  &  &  &  &  \cr%19
  \nr
  \bar3_4&&  &  &  &  &  &  &  &  &  &  &  &  &  &  &  &  &  &\x&  &  \cr%20
  \tr
 }
}
$$
\end{figure}
\end{eg}

Note the following simple facts, which follow directly from the definitions:

\begin{rem}\label{rem:picw}

(1)
The picture of $w$ contains exactly one $\cross$ in each row
 and in each column of cells.

(2)
If two $\cross$'s  are in the same stratum,
 then the one on the right is in a row below than the one on the left.
(We say that the $\cross$'s occur \em{in increasing order} within a stratum.)
\end{rem}

\subsection{Shape array for $\boldsymbol w$---local rules}\label{subsec:lr}
Now if $A=(i,j)$ is any lattice point,
 let $\tilde w(A)=\tilde w(i,j)$ denote the word in $\tilde\Gam_w$
 obtained from the rectangular section of the
 picture to the left of and above the vertex $A$, i.e.,
 $\tilde w(A)$ is the subword of $\tilde w_1\tilde w_2\cdots\tilde
 w_j$ consisting of all letters with ordinals $\le i$.
Let $w(A)=w(i,j)$ denote the word in $\Gam_n$ obtained from $\tilde w(A)$
 by discarding the subscripts of the letters.
By the above Remark (2), $\tilde w(A)$ equals the standardization
 of $w(A)$.
Let $\Lam(A)=\Lam(i,j)$ denote the shape of the $Sp(2n)$-tableau obtained
 by applying Berele's correspondence to $w(A)$.
By Lemma~\ref{LemStand}, it is also the final shape obtained by applying
 the standardized Berele correspondence to $\tilde w(A)$.

\begin{eg}
\label{eg:wordsec}
Let $w$ as in our previous examples.
For $A=(7,8)$, the relevant region is
$$
\vcenter{\packed
 \def\hcell#1{\vbox{\hbox to 1.5em{\hss$#1$\hss}}}
 \def\vcell#1{\vbox to 1.5em{\vss\hbox{$#1$}\vss}}
 \def\cell#1{\vbox to 1.5em{\vss\hbox to 1.5em{\hss$#1$\hss}\vss}}
 \def\nr{\omit&\multispan9\hrulefill\cr}
 \def\tr{\omit&\multispan9\leaders\hrule height1.0pt\hfill\cr}
 \def\x{\cross}
 \halign{\phantom{$#$}&\hskip 0.4pt#&&\hcell{#}\hskip 0.4pt\cr
      1_1&& 1& 2& 3& 4& 5& 6& 7& 8\cr
 }
 \halign{\vcell{#}&\vrule#&&\cell{#}\vrule\cr
  \tr
      1_1&&  &\x&  &  &  &  &  &  \cr%1
  \nr
      1_2&&  &  &  &  &  &  &\x&  \cr%2
  \nr
      1_3&&  &  &  &  &  &  &  &  \cr%3
  \tr
  \bar1_1&&  &  &  &  &  &\x&  &  \cr%4
  \nr
  \bar1_2&&  &  &  &  &  &  &  &  \cr%5
  \nr
  \bar1_3&&  &  &  &  &  &  &  &  \cr%6
  \tr
      2_1&&  &  &  &  &  &  &  &\x\cr%7
  \nr
 }
}
\quad
\text{so that}
\quad
\tilde w(A)=1_1 \bar1_1 1_2 2_1,
\quad
w(A)=1 \bar1 1 2.
$$
\end{eg}

Then we have the following:

\begin{thm}\label{ThmLocRules}
\rom{(1)}
Consider any cell, located at $(i,j)$, and let $A=(i-1,j-1)$, $B=(i-1,j)$,
 $C=(i,j-1)$, and $D=(i,j)$ be the four vertices surrounding the cell.
Then the quadruple of shapes $(\Lam(A),\Lam(B),\Lam(C),\Lam(D))$
 falls into exactly one of the following cases.
Note that, only the case marked as $(\cross)$
 has an $\cross$ written in the cell.

\noindent\rom(The carry-over group\rom)
\begin{enumerate}
\item[(\phantom{$\cross$})] $\Lam(A)=\Lam(B)=\Lam(C)=\Lam(D)$
\item[(\phantom{$\cross$})] $\Lam(A)=\Lam(B)\neq\Lam(C)=\Lam(D)$
\item[(\phantom{$\cross$})] $\Lam(A)=\Lam(C)\neq\Lam(B)=\Lam(D)$
\end{enumerate}

These three cases are visualized as follows.

$$
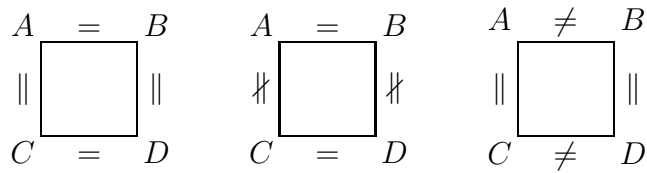
\cell{A}{=}{B}{\parallel}{}{\parallel}{C}{=}{D}
\qquad
\cell{A}{=}{B}{\nparallel}{}{\nparallel}{C}{=}{D}
\qquad
\cell{A}{\neq}{B}{\parallel}{}{\parallel}{C}{\neq}{D}
$$

The rest of the cases will be displayed visually.
The parenthesized symbol preceding each picture
 is the name of the case.
The symbols $\coverrowbelow{k}$ and $\coverrow{k}$ are synonymous.

\noindent\rom(The R-S group\rom)
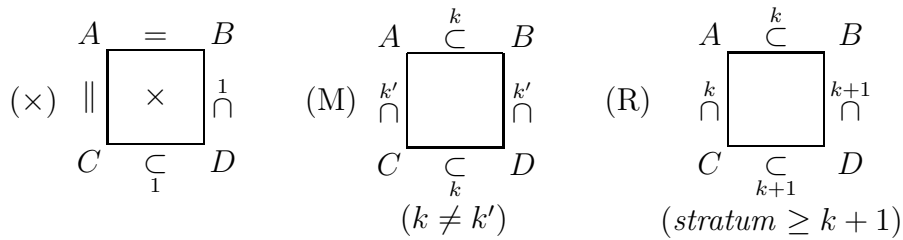
$$
\text{\rom{($\cross$)}}\;
 \cell{A}{\overset{\phantom{1}}{=}}{B}
  {\parallel}{\cross}{\vcoveredrow{1}}
  {C}{\coveredrowbelow{1}}{D}
\qquad
\text{\rom{(M)}}\;\underset{\textstyle\text{\rom{($k\neq k'$)}}}{
 \cell{A}{\coveredrow{k}}{B}
  {\vcoveredrow{k'}}{}{\vcoveredrow{k'}}
  {C}{\coveredrowbelow{k}}{D}
 }
\qquad
\text{\rom{(R)}}\;\underset{\textstyle\text{\rom(stratum $\ge k+1$\rom)}}{
 \cell{A}{\coveredrow{k}}{B}
  {\vcoveredrow{k}}{}{\vcoveredrow{k+1}}
  {C}{\coveredrowbelow{k+1}}{D}
 }
$$

\noindent\rom(Cancellation\rom)
$$
\text{\rom{($\bigcirc$)}}\;\underset{\textstyle\text{\rom(stratum $\le\bar k$\rom)}}
 {\cell{A}{\coveredrow{k}}{B}{\vcoveredrow{k}}{\bigcirc}{\vcoverrow{k}}
 {C}{\coverrowbelow{k}}{D}}
$$

\rom(The sign $\bigcirc$ inside the cell is just for easy recognition;
 it is not part of the initial data.\rom)

\noindent\rom(The jeu de taquin group\rom)

\begin{alignat*}{2}
\text{\rom{(\B J)}}\;&
\cell{A}{\covered}{B}{\vcover}{}{\vcover}{C}{\covered}{D}\;
\left(\begin{aligned}&\Lam(B)/\Lam(C)\\ &\;\;=\hdom\text{ or }\vdom\\
        &\Lam(A)=\Lam(D)\end{aligned}\right)
&\qquad
\text{\rom{(J)}}\;&
\cell{A}{\covered}{B}{\vcover}{}{\vcover}{C}{\covered}{D}\;
\left(\begin{aligned}&\Lam(B)/\Lam(C)\\ &\;\;\neq\hdom\text{ or }\vdom\\
        &\Lam(A)\neq\Lam(D)\end{aligned}\right)
\\
\text{\rom{(\B J$'$)}}\;&
\cell{A}{\cover}{B}{\vcovered}{}{\vcovered}{C}{\cover}{D}\;
\left(\begin{aligned}&\Lam(C)/\Lam(B)\\ &\;\;=\hdom\text{ or }\vdom\\
        &\Lam(A)=\Lam(D)\end{aligned}\right)
&\qquad
\text{\rom{(J$'$)}}\;&
\cell{A}{\cover}{B}{\vcovered}{}{\vcovered}{C}{\cover}{D}\;
\left(\begin{aligned}&\Lam(C)/\Lam(B)\\ &\;\;\neq\hdom\text{ or }\vdom\\
        &\Lam(A)\neq\Lam(D)\end{aligned}\right)
\end{alignat*}

\noindent\rom(The reverse R-S group\rom)
$$
\text{\rom{(W)}}\;
\underset{\textstyle\text{\rom{($k\neq k'$)}}}{
 \cell{A}{\coverrow{k}}{B}{\vcoverrow{k'}}{}{\vcoverrow{k'}}
 {C}{\coverrowbelow{k}}{D}
}
\qquad
\text{\rom{({\cyr Ya})}}\;
\underset{\textstyle\text{\rom{($k\ge 2$)}}}{
\cell{A}{\coverrow{k}}{B}{\vcoverrow{k}}{}{\vcoverrow{k-1}}
 {C}{\coverrowbelow{k-1}}{D}
}
$$

\rom{(2)}
The three shapes $\Lam(A)$, $\Lam(B)$, $\Lam(C)$,
 and the stratum containing the cell $ABCD$, together with
 the contents of the cell,
 determines which of the above cases the cell belongs to,
 and the shape $\Lam(D)$.
The list in \rom{(1)}, thus read as rules to determine $\Lam(D)$
 from the information stated immediately above,
 will be called the \em{local rules}.
We can recover the whole array of $\Lam(\cdot)$ from the picture of $w$
 by starting from the empty shapes on the top and the leftmost edges
 and applying these local rules in any possible order.

\rom{(3)}
The three shapes $\Lam(B)$, $\Lam(C)$, $\Lam(D)$
 and the stratum containing the cell $ABCD$
 determines which of the above cases the cell belongs to,
 and accordingly the contents of the cell and the shape $\Lam(A)$.
In other words, the local rules are invertible.
We can recover the whole array of $\Lam(\cdot)$
 and the positions of the $\cross$'s \rom(i.e.\ the word $w$\rom)
 if the shapes on the bottom and the rightmost edges
 are correctly given, together with the stratification.
In other words, the map which takes $w$ to the shapes on the bottom and the
 rightmost edges is injective.

\rom{(4)}
The sequence of shapes on the bottom edge
 equals the up-down tableau of degree $f$
 obtained from $w$ by Berele's algorithm,
 namely the Berele $Q$-symbol of $w$.

\rom{(5)}
Let $1\le k\le n$.
The sequence of shapes on the rightmost edge in the $k$-stratum
 represents a horizontal strip growing from left to right.

\rom{(6)}
The sequence of shapes on the rightmost edge in the $\bar k$-stratum
 represents a shrink by a horizontal strip from right to left,
 followed by a growth by another horizontal strip from left to right.
Moreover, if one puts $\lam\th k=\Lam(\ord(k_1)-1,f)$ and
 $\mu\th k=\Lam(v_k,f)$, where $v_k$ is the row coordinate of
 the turning point from shrink to growth,
 then $\mu\th k/\lam\th k$ is also a horizontal strip.

\rom{(7)}
The tableau of shape $\Lam(f,f)$
 in which $\mu\th k/\lam\th k$ is filled by the symbol $k$
 and $\lam\th{k+1}/\mu\th k$ is filled by the symbol $\bar k$
 \rom($k=1$, $2$, \dots, $n$,
 where $\lam\th{n+1}$ is understood to be $\Lam(f,f)$\rom)
 is the $Sp(2n)$-tableau obtained from $w$ by Berele's correspondence,
 namely the Berele $P$-symbol of $w$.

\end{thm}

\begin{rem}\label{rem:loc}

(1)
Theorem~\ref{ThmLocRules}\ says that the result of Berele's correspondence
 can be completely determined by the ``local rules'' listed in (1).

(2) The above set of local rules is an expansion of Fomin's local rules
 in \cite{Fomin} for the Robinson-Schensted correspondence.
The latter consist of the rules in the carry-over group and the R-S group only,
 in which the stratum condition in (R) does not appear.

(3) The rules in the jeu de taquin group were used by
S.~Fomin \cite{FominGen} and M.~van Leeuwen \cite{vanLeeuwen}, to give a
pictorial presentation of Sch\"utzenberger's  involution.

(4) The local rules thus expanded have certain restricted symmetries,
namely a $180^o$ rotation and reflection in one diagonal.
The restriction derives from dependence on stratification
in distinguishing between cases (R) and ($\bigcirc$).

(5) The local rules lead to only two possible types of rows
 in the picture of $w$, namely:

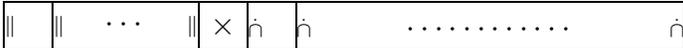
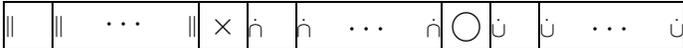
\def\vr{\vrule height18.8pt}
\def\cl#1#2#3#4{\def\hr{\hrule width#1}\def\hb{\hbox to#1}%
 \vbox to18.8pt{\hr\vfil\hb{\rlap{$\ssize#2$}\hfil$#3$\hfil\llap{$\ssize#4$}}%
 \vfil\hr}}
\def\cell#1#2#3{\cl{18pt}{#1}{#2}{#3}}
\def\longcell#1#2#3{\cl{54.8pt}{#1}{#2}{#3}}
\def\verylongcell#1#2#3{\cl{146.8pt}{#1}{#2}{#3}}
\def\veq{\parallel}
\begin{align*}
&\hbox{\vr\cell{\veq}{}{}\vr\longcell{\veq}{\cdots}{\veq}%
  \vr\cell{}{\cross}{}\vr\cell{\vcovered}{}{}%
  \vr\verylongcell{\vcovered}{\cdots\cdots\cdots\cdots}{\vcovered}\vr}
\\
\intertext{or}
&\hbox{\vr\cell{\veq}{}{}\vr\longcell{\veq}{\cdots}{\veq}%
  \vr\cell{}{\cross}{}\vr\cell{\vcovered}{}{}%
  \vr\longcell{\vcovered}{\cdots}{\vcovered}\vr\cell{}{\bigcirc}{}%
  \vr\cell{\vcover}{}{}\vr\longcell{\vcover}{\cdots}{\vcover}\vr}.
\end{align*}

The same statement applies to columns of the picture of $w$.

(6)
The local rules also guarantee that the shapes in the $k$- and $\bar k$-%
 strata (including the bottom lines thereof) cannot have more than $k$ parts.
This can be shown by induction on $k$, assuming its validity
 at the top of the $k$-stratum, and proceeding row by row in the $k$- and
 $\bar k$-strata as follows.
Suppose the property is satisfied for the vertices $(i-1,j')$, $0\le j'\le f$
 and $i$ is still in the $k$- or $\bar k$-stratum.
It is sufficient to show that there is no growth in the $(k+1)$st or lower row
 of the Young diagram along the vertical segment $( i - 1, j' )$--$( i, j' )$
 in the picture.
By Remark~\ref{rem:loc}~(5), we can concentrate on the interval
 starting at the right edge of the cell ($\cross$)
 (where the growth always starts in the 1st row) and ending at the left
 edge of the cell ($\bigcirc$) (or the rightmost edge of the picture
 if this row does not have ($\bigcirc$)).
Looking at the local rules, we know that the growth row number
 is either preserved (rules ($\phantom{\cross}$), (M), (J$'$),
 (\B J$'$) with $\Lam(C)/\Lam(B)=\hdom$), or changes by one
 (increases in (R), decreases in (\B J$'$) with $\Lam(C)/\Lam(B)=\vdom$),
 as we cross over a cell.
Therefore the growth row number must turn from $k$ to $k+1$ at some stage,
 if it ever exceeds $k$.
However, such change is not allowed in case (R) by the stratum condition.
So we cannot have growths in rows below the $k$th,
 and since we do not have more than $k$ rows on the vertices $(i-1,j')$,
 the same holds for the vertices $(i,j')$.
\end{rem}

\subsection{Structure of proof of Theorem~\ref{ThmLocRules}}\label{subsec:str}
The rest of this section is devoted to the proof of Theorem~\ref{ThmLocRules}.

The proof proceeds by induction based on a natural poset structure
 defined on the set of lattice points $[0,f]\times[0,f]$:
 if $(i',j')$ and $(i,j)$ are two lattice points,
 then $(i',j')\le(i,j)$ if and only if $i'\le i$ and $j'\le j$,
 in other words, $(i',j')$ lies in the closed rectangle with vertices
 $(0,0)$, $(0,j)$, $(i,0)$, and $(i,j)$.
We denote by $L$ the poset $[0,f]\times[0,f]$ defined in this manner.
An \em{order ideal} of $L$ is a subset $I$ of $L$
 for which $(i,j)\in I$ and $(i',j')\le (i,j)$ imply $(i',j')\in I$.
In the picture $I$ is a set of lattice points in $L$ which is saturated
 to the above and to the left.
We say that $(i,j)$ is a \em{cocorner vertex} of $I$
 if $(i,j)$, $(i,j+1)$, $(i+1,j)\in I$ but $(i+1,j+1)\not\in I$.
If $A$ is any vertex in $L$, we denote by $\tilde P(A)$ (resp.\ $P(A)$)
 the $P$-symbol %%%%%%%%%%%
 obtained by applying standardized Berele's correspondence %%%%%%%%%%
 (resp.\ original Berele's correspondence)
 to the word $\tilde w(A)$ (resp.\ $w(A)$).
Note that $\Lam(A)=\sh(\tilde P(A))=\sh(P(A))$.

We will show the following Lemma by induction on $I\in \mathcal{J}(L)$,
 the lattice of order ideals in $L$.
The lemma is concerned with all $\tilde P(A)$, $A\in I$
 as well as all $\Lam(A)$, $A\in I$.
This will readily imply Theorem~\ref{ThmLocRules}\ (1) by putting $I=L$.

\begin{lemma}
\label{LemLocRules}
Let $I$ be an order ideal in $L$.

\rom{(1)} If $A$ and $B$ are horizontally adjacent vertices in I,
 with $B$ to the right of $A$,
 then we have either $\Lam(A)=\Lam(B)$ \rom(called an \em{equal}\rom),
 $\Lam(A)\covered\Lam(B)$ \rom(a \em{growth}\rom),
 or $\Lam(A)\cover\Lam(B)$ \rom(a \em{shrink}\rom).

\rom{(2)} If $A$ and $C$ are vertically adjacent vertices in I,
 with $C$ below $A$,
 then we have either $\Lam(A)=\Lam(C)$ \rom(an \em{equal}\rom),
 $\Lam(A)\covered\Lam(C)$ \rom(a \em{growth}\rom),
 or $\Lam(A)\cover\Lam(C)$ \rom(a \em{shrink}\rom).
Moreover, the corresponding $P$-symbols satisfy
 one of the following relations.

\rom{(2a)} If $\Lam(A)=\Lam(C)$, then we have $\tilde P(A)=\tilde P(C)$.

\rom{(2b)} If $\Lam(A)\covered\Lam(C)$,
 and if $C$ has coordinates $(\ord(\gam_t),j)$, $\gam_t\in\tilde\Gam_w$,
 then one can obtain $\tilde P(C)$ from $\tilde P(A)$
 by filling the new cell $\Lam(C)\setminus\Lam(A)$ with $\gam_t$.

\rom{(2c)} If $\Lam(A)\cover\Lam(C)$,
 and if $C$ has coordinates $(\ord(\gam_t),j)$, $\gam_t\in\tilde\Gam_w$,
 then the following \rom{(2c1)}--\rom{(2c4)} hold:

\rom{(2c1)} We have $\gam=\bar k$ for some $k\in[1,n]$.

\rom{(2c2)} With $k$ defined as in \rom{(2c1)},
 the tableau $\tilde P(A)$ does not contain any $\bar k$,
 so that $k$ is the largest possible letter in $\tilde P(A)$.

\rom{(2c3)} If $\{ ( r, c ) \} = \Lam( A ) \setminus \Lam( C )$,
 then each of the bottom cells of the 1st through $c$th columns
 of $\tilde P(A)$ contains a $k$ \rom($k$'s can appear
 in other columns as well\rom).

\rom{(2c4)} The tableau $\tilde P(C)$ is obtained from $\tilde P(A)$
 by removing the $k$ in the 1st column \rom(which is the ``smallest'' $k$\rom),
 and then shifting each $k$ sitting at the bottoms
 of the 2nd through $c$th columns to the bottom of its left adjacent column.
If we discard the subscripts, $P(C)$ is simply obtained from $P(A)$
 by removing $k$ at $(r,c)$.

\rom{(3)} If $A$, $B$, $C$, and $D$ are the four vertices of a cell
 contained in $I$, with $D=(i,j)$,
 then the quadruple $(\Lam(A),\Lam(B),\Lam(C),\Lam(D))$ falls
 into exactly one of the cases listed in the local rules.

\end{lemma}

\begin{rem}\label{rem:AC}
(1)
The relation between $\tilde P(A)$ and $\tilde P(C)$ is not
trivial---nothing a priori assures any relation between $\tilde P(A)$
and $\tilde P(C)$, 
 as opposed to $\tilde P( A )$ and $\tilde P( B )$,
 which are directly connected by standardized Berele insertion.

(2)
The procedure to obtain $\tilde P(C)$ from $\tilde P(A)$ described in (2c4)
 can be understood to be ``column deletion,'' namely the tableau deletion
 procedure (as described in \cite[Section 5.1.5]{KnuthACP}),
 modified to serve as the inverse of the column insertion
 instead of the row insertion.  It is also a semistandard version of a
bijective tool used by Sundaram~\cite[Proof of Lemma~8.7]{Sundaram}.  
\end{rem}

\begin{proof}
We prove Lemma~\ref{LemLocRules} by induction using Lemma~\ref{LemLRStart} and Lemma~\ref{LemLRStep}.

\begin{lemma}
\label{LemLRStart}

Define an order ideal $I_0$ of $L$ by $I_0=\{\,(0,j)\mid j\in[0,f]\,\}\cup
 \{\,(i,0)\mid i\in[0,f]\,\}$.
Then Lemma~\ref{LemLocRules} holds for $I=I_0$.
\end{lemma}

This is clear since, by definition,
 the word $\tilde w(A)$ is the empty word for any $A\in I_0$,
 so that $\Lam(A)$ and $\tilde P(A)$ are all empty.
The essential point of the proof lies in the following inductive step.

\begin{lemma}
\label{LemLRStep}
{\bf (the key technical lemma)}
Let $I_1$ be an order ideal of $L$ which contains vertices $(i-1,j-1)$, 
$(i-1,j)$, and $(i,j-1)$, but not $(i,j)$. Note that
$I_2=I_1\cup\{(i,j)\}$ is also an order ideal of $L$. Then if
Lemma~\ref{LemLocRules} holds for $I_1$, then it also holds for $I_2$. 
\end{lemma}

With this admitted, Lemma~\ref{LemLocRules} is proved in the following manner.
If $I\neq L$, then $I$ necessarily has at least one cocorner.
Therefore, starting with $I_0$,
 we can continue to enlarge the region of validity of Lemma~\ref{LemLocRules} %
 by applying Lemma~\ref{LemLRStep} until we reach the case $I=L$.
\end{proof}

\subsection{Proof of Lemma~\ref{LemLRStep} (the key technical lemma)}\label{subsec:key}

Let $A=(i-1,j-1)$, $B=(i-1,j)$, $C=(i,j-1)$, and $D=(i,j)$.
Also let $i=\ord(\gam_t)$ where $\gam\in\Gam_n$ and $\gam_t\in\tilde\Gam_w$.
{\bf The symbols $\boldsymbol A$, $\boldsymbol B$, $\boldsymbol C$, $\boldsymbol D$, $\boldsymbol\gam$ and $\boldsymbol\gam_{\boldsymbol t}$ will be used throughout the proof in this sense.}

Since the claims already hold for segments and cells contained in $I_1$,
 the new claims to be proved are the claims (1) for the segment $CD$,
 (2) for the segment $BD$,
 and (3) for the cell $ABCD$.
Looking at the local rules in Theorem~\ref{ThmLocRules},
 one can see that (1) for $CD$ and the initial part of (2) for $BD$
 will follow from (3) for $ABCD$.
So we will concentrate on the validity of (3) for $ABCD$,
 and the latter part of (2) for $BD$, i.e., the relation between
 $\tilde P(B)$ and $\tilde P(D)$.

First we show that the three shapes $\Lam(A)$, $\Lam(B)$, and $\Lam(C)$,
 the stratum containing $ABCD$, and the contents of $ABCD$
 matches exactly one of the cases listed in Theorem~\ref{ThmLocRules}\ (1).
Since the cases are disjoint with respect to these data,
 what we have to show is that they cover all possibilities.
Since (1) holds for $AB$ and (2) holds for $AC$,
 we know that $\Lam(A)$ and $\Lam(B)$ (resp.\ $\Lam(A)$ and $\Lam(C)$)
 are either equal or one covers the other.  
Looking again at the list,
 we find that what we need to show is the validity of the following
 two statements:
\begin{enumerate}
\item[(a)] $ABCD$ can contain an $\cross$ {\it only if\/}
$\Lam(A)=\Lam(B)=\Lam(C)$.  
\item[(b)] If $\Lam(A)\coverrow{k}\Lam(B)=\Lam(C)$, then $k>1$.
\end{enumerate}

To see (a), first recall that an $\cross$ can appear only once
 in each row and each column of our picture of $w$.
Therefore if the cell $(i,j)$ of the picture of $w$ contains $\cross$,
 the cells $(i,1)$ through $(i,j-1)$ must be empty,
 so that $\tilde w(A)=\tilde w(C)$ and $\Lam(A)=\Lam(C)$.
Similarly the cells $(1,j)$ through $(i-1,j)$ must also be empty,
 so that $\tilde w(A)=\tilde w(B)$ and $\Lam(A)=\Lam(B)$.

To show (b) we first prove some easy facts, which will also prove useful
in later arguments.

\begin{lemma}
\label{LemCan}

In order that the standardized Berele insertion of the letter $\bet_s$
 into $\tilde P$ involves cancellation of a $k$-$\bar k$ pair,
 $\tilde P$ must contain at least one $\bar k$.
\end{lemma}

\begin{proof}[Proof of Lemma~\ref{LemCan}]
This is clear
 since the cancellation of $k$-$\bar k$ is first caused by a $\bar k$,
 bumped by a $k$ from the $k$th row into the $(k+1)$st row.
\end{proof}

\begin{lemma}
\label{LemAC}
Let $A$, $B$, $C$, $D$ be as in Figure~\ref{fig:abcd}.  

Assume that $\Lam(A)\neq\Lam(B)$ and $\Lam(A)\neq\Lam(C)$.
Then we have an $\cross$ representing $\gam_t$ straight to the left of $ABCD$.
Also there must be an $\cross$ straight above $ABCD$.
Let $\bet_s$ be the letter represented by the latter $\cross$.
Then $\bet$ must be strictly smaller than $\gam$ in $\Gam_n$.
\end{lemma}

\begin{proof}[Proof of Lemma~\ref{LemAC}]
Since $\bet_s < \gam_t$, $\bet \le \gam$ is clear.
The $\cross$'s in the same stratum must occur in increasing order,
 so the possibility of $\bet = \gam$ is rejected.

\end{proof}

\begin{lemma}
\label{LemYa}

In the notation and situation of Lemma~\ref{LemAC},
assume that $\Lam( A ) \cover \Lam( B )$ and $\Lam( A ) \cover \Lam( C )$.
Putting $X = \Lam( A ) \setminus \Lam( B ) = \{ ( r, c ) \}$ and
 $Y = \Lam( A ) \setminus \Lam( C )$, further assume that either $X = Y$
 or $X$ lies in a lower row than $Y$ in $\Lam( A )$.
Then the sliding path occurring in the standardized Berele insertion
 $\tilde P( A ) \StandBereleinsert \bet_s$
 must involve the square $(r-1,c)$ of the Young diagram.
\rom(This requires that the square $(r-1,c)$ exists,
 so that $r \ge 2$.\rom)
\end{lemma}

\begin{proof}[Proof of Lemma~\ref{LemYa}]
Since $\tilde P(A)$ contains no $\bar k$ (by (2c2) applied to $AC$),
 the cancellation involved in $\tilde P(A)\StandBereleinsert\bet_s$
 must be for a pair smaller than $k$-$\bar k$ due to Lemma~\ref{LemCan}.
Accordingly, the largest possible letter that gets moved (or removed) in
the bumping phase is $\overline{k-1}$
 (recall that the bumped letter moves a letter larger than itself).
Therefore, the $k$'s in $\tilde P(A)$ do not move in the bumping phase,
nor is any cancelled.
This assures that, when sliding starts, the bottoms of the 1st through
 $c$th columns are still occupied by $k$'s,
 as was the case before the bumping starts, by (2c3) applied to $AC$
 and our assumption on the positions of $X$ and $Y$.
Now in the sliding phase,
each of these $k$'s either stays put,
 moves left, or moves up.
When the sliding is finished, the tableau must be semistandard;
 in particular no two $k$'s can share a single column.
This is only possible if each of the $k$'s in the first $c$ columns
 sticks to its column.
Therefore the $k$ at $(r,c)$, which was the end of the sliding path,
 must have moved up into the cell $(r-1,c)$.
\end{proof}

Returning to the proof of Lemma~\ref{LemLRStep}, we find that (b) follows from the
above lemma, finishing the proof of the coverage of all possible cases.

It remains to show that the shape $\Lam(D)$
 agrees with that given by the local rules,
 and that the relation between $\tilde P(B)$ and $\tilde P(D)$ is
 as stated in (2).
We will do these together case by case.
First we deal with some easy cases, which are marked by (\phantom{$\cross$})
 in the list.

{\bf The carry-over group.  }\label{subsec:cog}
First assume $\Lam(A)=\Lam(C)$.
The cells straight to the left of $ABCD$ do not contain an $\cross$ and,
 by assumption, the cell $ABCD$ is also empty in this case.
Therefore we have $\tilde w(B)=\tilde w(D)$, so that $\tilde P(B)=\tilde P(D)$
 and $\Lam(B)=\Lam(D)$, validating the prescription by the local rule.
Also, $\tilde P(B)=\tilde P(D)$ validates Lemma~\ref{LemLocRules} (2a) for $BD$.

Next assume $\Lam(A)\neq\Lam(C)$ and $\Lam(A)=\Lam(B)$.
By a similar argument we have $\tilde w(C)=\tilde w(D)$,
 $\tilde P(C)=\tilde P(D)$, and $\Lam(C)=\Lam(D)$.
This validates the prescription for $\Lam(D)$.
Also, the relation between $\tilde P(B)$ and $\tilde P(D)$
 ((2b) or (2c) for $BD$)
 is validated because the same pair of tableaux occur on $AC$,
 which is contained in $I_1$.

{\bf The R-S group.  }\label{subsec:rsg}
{\sc The case ($\cross$)}.
We have $\tilde w(A)=\tilde w(B)=\tilde w(C)$,
 which implies $\tilde P(A)=\tilde P(B)=\tilde P(C)$.
Since $\tilde P( D ) = \tilde P( C ) \StandBereleinsert \gam_t$,
 it also equals $\tilde P( B ) \StandBereleinsert \gam_t$.
Note that $\gam_t$ is larger than any letter in $\tilde w(B)$.
This means that $\tilde P(D)$ is obtained from $\tilde P(B)$ by adding $\gam_t$
 at the end of the 1st row.
Therefore we have $\Lam( B ) = \Lam( C ) \coveredrow 1 \Lam( D )$,
 and (2b) holds for $BD$.

{\bf More Notation.}
In the remaining cases, we have $\Lam(A)\neq\Lam(B)$ and $\Lam(A)\neq\Lam(C)$.
Therefore we have an $\cross$ straight to the left of $ABCD$ representing
 the letter $\gam_t$.
Also we have an $\cross$ straight above $ABCD$.
Let $\bet_s$ be the letter represented by this $\cross$.
{\bf In the rest of the proof, $\boldsymbol\bet$ and $\boldsymbol\bet_{\boldsymbol s}$ will always be used in this sense.}

{\sc The case (M)}.
We have an $\cross$ straight to the left of $ABCD$.
By (2b) applied to $AC$ (which is in $I_1$),
 the difference between $\tilde P(C)$ and $\tilde P(A)$ is
 that $\tilde P(C)$ has an extra $\gam_t$
 at position $\Lam(C)\setminus\Lam(A)$,
 and $\gam_t$ is the largest letter appearing in $\tilde P(C)$.
By definition we have $\tilde P(B)=\tilde P(A)\StandBereleinsert\bet_s$ and
 $\tilde P(D)=\tilde P(C)\StandBereleinsert\bet_s$.
Moreover, the assumption of the case $\Lam(A)\covered\Lam(B)$ assures
 that the operation $\tilde P(A)\StandBereleinsert\bet_s$
 is an ordinary insertion.

We want to know how $\tilde P(C)\StandBereleinsert\bet_s$ differs
 from $\tilde P(A)\StandBereleinsert\bet_s$.
Since the position $\Lam(C)\setminus\Lam(A)$ is a cocorner of $\Lam(A)$,
 it can be a part of the bumping path for the insertion $\tilde P(A)
 \StandBereleinsert\bet_s$ only as its final point where a new cell is added.
However, since $\Lam(B)\neq\Lam(C)$,
 this bumping path ends at some other cocorner of $\Lam(A)$.
If this cocorner is in a row above $\Lam(C)\setminus\Lam(A)$,
 then the bumping path remains the same for the insertion
 $\tilde P(C)\StandBereleinsert\bet_s$.
If this cocorner is in a row below, the bumping path also remains the same
 for $\tilde P(C)\StandBereleinsert\bet_s$ since,
 when it passes the row of $\Lam(C)\setminus\Lam(A)$,
 it can already bump an element of $\tilde P(A)$ and the extra $\gam_t$
 remains intact since it is definitely greater than that element.
In either case the insertion of $\bet_s$ into $\tilde P(B)$ causes exactly
 the same change as it would cause to $\tilde P(A)$,
 leaving $\gam_t$ where it is.
This validates the local rule for $\Lam(D)$ and (2b).

{\sc The case (R)}.
The relation between $\tilde P(A)$ and $\tilde P(C)$,
 $\tilde P(A)$ and $\tilde P(B)$, $\tilde P(C)$ and $\tilde P(D)$
 are the same as in the previous case.
This time we are under the assumption that $\Lam(B)=\Lam(C)$,
 so the bumping path of inserting
 $\bet_s$ into $\tilde P(C)$ hits $\Lam(C)\setminus\Lam(A)$,
 and tries to bump the extra $\gam_t$, say in the $k$th row, to row $k+1$.
Under the stratum condition attached to (R), $\gam_t$ can sit in row $k+1$
 without violating the symplectic condition,
 so $\tilde P(C)\StandBereleinsert\bet_s$ ends in an ordinary insertion,
 the result being $\tilde P(B)$ plus an extra $\gam_t$
 at the end of row $k+1$.
This validates the local rule for $\Lam(D)$ and (2b).

\begin{figure}
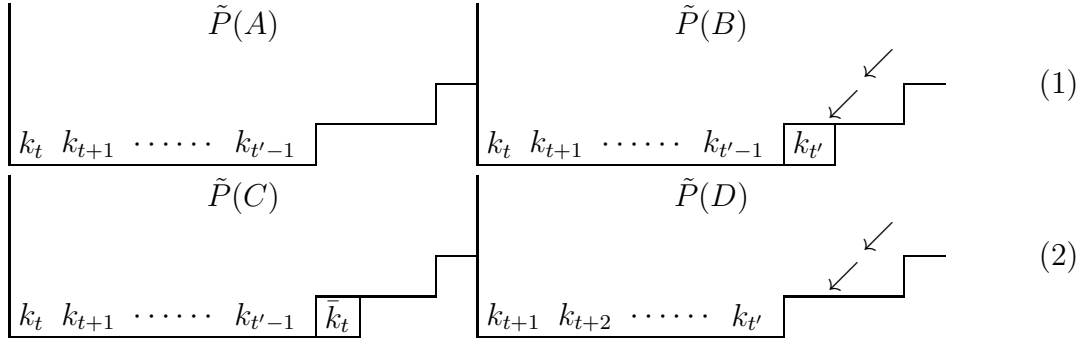

\caption{The cancellation case}
\label{fig:pfo}
\def\vcell#1{\vbox to 15pt{\vss\hbox{\hss\;$#1$\;\hss}\vss}}
\allowdisplaybreaks
\begin{alignat}{2}
&
\vcenter{\packed
 \halign{&#\cr
%  \vrule height30pt&\multispan9\hfil\vcell{\tilde P(A)}\hfil\cr
  \vrule height15pt&\multispan9\hfil\vcell{\tilde P(A)}\hfil\cr
  \vrule height15pt\cr
  \vrule&\multispan6\hfill&\vrule&\vbox{\hrule width15pt}&\vrule\cr
  \vrule&\multispan6\hfill&\vrule height15pt\cr
  \vrule&\multispan4\hfill&\vrule&\vbox{\hrule width45pt}&\vrule\cr
  \vrule&\vcell{k_t}&\vcell{k_{t+1}}&\vcell{\cdots\cdots}&\vcell{k_{t'-1}}
   &\vrule\cr
  \multispan6\hrulefill\cr
 }
}
&
%\qquad
&
\vcenter{\packed
 \halign{&#\cr
%  \vrule height30pt&\multispan{11}\hfil\vcell{\tilde P(B)}\hfil\cr
  \vrule height15pt&\multispan{11}\hfil\vcell{\tilde P(B)}\hfil\cr
  \vrule&\multispan5\hfill&\multispan3\hfil\vcell{\swarrow}\cr
  \vrule&\multispan8\hfill&\vrule&\vbox{\hrule width15pt}&\vrule\cr
  \vrule&\multispan5\hfill&\multispan3\hfil\vcell{\swarrow}\hfil&\vrule\cr
  \vrule&\multispan4\hfill&\vrule&\multispan3\vbox{\hrule width45pt}&\vrule\cr
  \vrule&\vcell{k_t}&\vcell{k_{t+1}}&\vcell{\cdots\cdots}&\vcell{k_{t'-1}}
   &\vrule&\vcell{k_{t'}}&\vrule\cr
  \multispan8\hrulefill\cr
 }
}
\\
\vspace{\medskipamount}
&
\vcenter{\packed
 \halign{&#\cr
%  \vrule height30pt&\multispan{11}\hfil\vcell{\tilde P(C)}\hfil\cr
  \vrule height15pt&\multispan{11}\hfil\vcell{\tilde P(C)}\hfil\cr
  \vrule height15pt\cr
  \vrule&\multispan8\hfill&\vrule&\vbox{\hrule width15pt}&\vrule\cr
  \vrule&\multispan8\hfill&\vrule height15pt\cr
  \vrule&\multispan4\hfill&\vrule&\multispan3\vbox{\hrule width45pt}&\vrule\cr
  \vrule&\vcell{k_t}&\vcell{k_{t+1}}&\vcell{\cdots\cdots}&\vcell{k_{t'-1}}
   &\vrule&\vcell{\bar k_t}&\vrule\cr
  \multispan8\hrulefill\cr
 }
}
&
%\qquad
&
\vcenter{\packed
 \halign{&#\cr
%  \vrule height30pt&\multispan9\hfil\vcell{\tilde P(D)}\hfil\cr
  \vrule height15pt&\multispan9\hfil\vcell{\tilde P(D)}\hfil\cr
  \vrule&\multispan5\hfill&\hfil\vcell{\swarrow}\cr
  \vrule&\multispan6\hfill&\vrule&\vbox{\hrule width15pt}&\vrule\cr
  \vrule&\multispan5\hfill&\hfil\vcell{\swarrow}\hfil&\vrule\cr
  \vrule&\multispan4\hfill&\vrule&\vbox{\hrule width45pt}&\vrule\cr
  \vrule&\vcell{k_{t+1}}&\vcell{k_{t+2}}&\vcell{\cdots\cdots}&\vcell{k_{t'}}
   &\vrule\cr
  \multispan6\hrulefill\cr
 }
}
\end{alignat}
\end{figure}

{\bf The cancellation group.  }\label{subsec:can}
{\sc The case ($\bigcirc$)}.
We proceed as in the previous case.
This time placing $\gam_t$ in row $k+1$ would violate
 the symplectic condition.
$\gam$ cannot be less than $k$
 since $\gam_t$ is allowed to sit in the $k$th row in $\tilde P(C)$,
 so $\gam$ is either $k$ or $\bar k$.
On the other hand, let $\gam'_{t'}$ denote the letter which bumped $\gam_t$.
This letter is also allowed to sit in the $k$th row in $\tilde P(B)$,
 and since it bumps $\gam_t$ we must have $\gam'<\gam$.
The only possible combination is $\gam'=k$ and $\gam=\bar k$.
Then, as explained in the third paragraph of
Subsection~\ref{subsec:stber} (p.~\pageref{1shiftpath}), the cells to
the left of $\Lam(C)\setminus 
 \Lam(A)$ in this row all contain $k$'s.
Moreover, by the proof of Lemma~\ref{LemStand}, $\bar k_t$ is the smallest $\bar k$,
 and the smallest $k$ in $\tilde P(C)$ has the same subscript $t$.
The instruction is to remove $k_t$ and shift the remaining $k$'s in this row
 to the left as part of the sliding.
Then sliding stops here since $\Lam(C)\setminus\Lam(A)$ is a corner
 of $\Lam(C)$.
The tableau $\tilde P(D)$ thus obtained matches the description
 of (2c4) relative to $\tilde P( B )$.
We have $\Lam(D)=\Lam(A)$, validating the local rule for the shape.
See Figure~\ref{fig:pfo}
Here $\swarrow$ signifies the bumping path that occurs in the first
$k-1$ rows of $\tilde P(B)$ and $\tilde P(D)$, that distinguish them
from $\tilde P(A)$ and $\tilde P(C)$ (respectively).  The two bump paths
are identical in this case, as shown above.

Now (2c1) is already validated since $\gam=\bar k$.
$\bar k_t$ is not only the smallest $\bar k$, but the only $\bar k$
 in $\tilde P(C)$ since it is the largest letter in $\tilde w(C)$.
Then (2b) applied to $AC$ implies
 that $\tilde P(A)$ does not contain $\bar k$.
In addition, the letter $\bet_s$ inserted into $\tilde P(A)$
 (to produce $\tilde P(B)$) cannot be a $\bar k$ due to Lemma~\ref{LemAC}.
Therefore $\tilde P(B)$ does not contain $\bar k$, which validates (2c2)
 for $BD$.
The argument in the previous paragraph validates (2c3) and (2c4) for $BD$.

{\bf The jeu de taquin group.}\label{subsec:jeug}
{\sc The case (\B J$'$)}.
Since $AC$ is a growth, by (2b) $\tilde P(C)$ is $\tilde P(A)$ plus
 an extra $\gam_t$ at some cocorner of $\Lam(A)$.
Put $Y = \Lam( C ) \setminus \Lam( A )$.
Since $AB$ is a shrink, the Berele insertion of $\bet_s$ into $\tilde P(A)$
 involves sliding, ending at some corner of $\Lam(A)$.
Put $X = \Lam( A ) \setminus \Lam( B )= \{ ( r, c ) \}$.
Since we are in case (\B J$'$),
 where $\Lam( C ) \setminus \Lam( B )$ is a domino,
 $Y$ is either immediately to the right of or below $X$.

If one starts with $\tilde P(C)$ instead of $\tilde P(A)$
 and Berele inserts the letter $\bet_s$,
 then bumping, cancellation, and sliding proceeds without any alteration
 until the hole arrives at $X$.
Note that since $X$ is a corner of $\Lam(A)$,
 only one of $(r+1,c)$ or $(r,c+1)$ can be a part of $\Lam(C)$.
One of them is $Y$, and the other is not in $\Lam(C)$;
hence, the entry at $Y$, namely $\gam_t$, slides into the hole at $X$.
Therefore, $\tilde P(D)$ has the same shape as $\tilde P(A)$,
 and is obtained from $\tilde P(B)$ by adding $\gam_t$ at $X$.
This validates the local rule for $\Lam(D)$ as well as (2b) for $BD$.

{\sc The case (J$'$)}.
Let $X$, $Y$ be as in the previous case.
This time $\Lam(C)\setminus\Lam(B)$ is not a domino,
 so $X$ is not adjacent to $Y$.

We consider how the operation $\tilde P(C)
 \StandBereleinsert\bet_s$ differs from $\tilde
P(A)\StandBereleinsert\bet_s$. 
In the bumping phase, the existence or nonexistence of $\gam_t$ at $Y$
 does not affect anything for the same reason as in case (M)
 (even if $X$ is in a row below than $Y$).
In the sliding phase, if the sliding path for $\tilde P(A)$
 does not pass through the cells adjacent to $Y$,
 then clearly the existence or nonexistence of $\gam_t$ at $X$ has no affect.
If the sliding path passes through some cell $Z$ adjacent to $Y$,
 then that cell cannot coincide with the end point $X$ of the sliding,
 under the current assumption.
Then in the next step, the letter that slides into $Z$ in $\tilde P(A)
 \StandBereleinsert\bet_s$ also slides there in $\tilde P(C)\StandBereleinsert
 \bet_s$ regardless of the existence of $\gam_t$
 (and $\gam_t$ does not move),
 since that letter is smaller than $\gam_t$.
The rest of the sliding is not affected either.
Therefore $\Lam( D ) = \Lam( A ) \cup Y \setminus X$,
 and $\tilde P(D)$ is obtained from $\tilde P(B)$ by adding $\gam_t$ at $Y$.
This validates the local rule for $\Lam(D)$ and (2b) for $BD$.

{\bf Still More Notation.}
In the remaining cases, $AC$ is always a shrink.
Put $Y = \Lam( A ) \setminus \Lam( C )$, and $Y = \{ ( r, c ) \}$.
By (2c) applied to $AC$, we have $\gam=\bar k$ for some $k$,
 $\tilde P(A)$ contains some $k_u$ (with some index $u$) at $Y$,
 and $\tilde P(C)$ is obtained from $\tilde P(A)$ by removing
 the smallest $k$ and shifting the $k$'s to the left
 inside the horizontal strip for $k$ up to $Y$, vacating the position $Y$.
{\bf In the rest of the proof, $\boldsymbol Y$, $\boldsymbol r$, $\boldsymbol c$,
 $\boldsymbol k$ and $\boldsymbol k_{\boldsymbol u}$ will always be used
 in this sense.}

\begin{rem}\label{rem:ACshrink}
In the remaining cases, the following points easily follow regardless of further case distinction.

(a) $\tilde P(B)$ does not contain any $\bar k$.
This is because $\tilde P(A)$ does not contain any $\bar k$
 due to (2c2) applied to $AC$,
 and the inserted letter $\bet_s$ cannot be a $\bar k$
 due to Lemma~\ref{LemAC}.
In particular, if $BD$ is shown to be a shrink, then (2c2) for $BD$ follows.

(b) If $BD$ is shown to be a shrink, then (2c1) for $BD$ follows
 from that of $AC$ since the segments $AC$ and $BD$ are in the same row.
\end{rem}

With these at hand, we continue our case-by-case analysis.

{\sc The case (\B J)}.
Put $X = \Lam( B ) \setminus \Lam( A )$.
Since $AB$ is a growth, the Berele insertion of $\bet_s$ into $\tilde P(A)$
 ends in an ordinary insertion at $X$.
Since we are in case (\B J), where $\Lam( B ) \setminus \Lam( C )$ is assumed
 to be a domino, $X$ is either immediately to the right of or below $Y$.

First assume that $X$ is to the right of $Y$.  See Figure~\ref{fig:pfbjh} 
Note that the occupant of $X$ in $\tilde P(B)$ is also a $k$, say $k_v$,
 since it has to be greater than the occupant of $Y$, which is out of
 the bumping path and hence still contains $k_u$,
 and it cannot be a $\bar k$ because of Remark~\ref{rem:ACshrink}~(a).
Also note that $k_v$ is the letter bumped from the previous row
 (or the inserted letter if $r=1$).
If one inserts $\bet_s$ into $\tilde P(C)$ instead of $\tilde P(A)$,
 the bumping proceeds up to the row of $Y$,
 where $k_v$ now lands at $Y$ instead of $X$.
Therefore the shape of $\tilde P(D)$ equals that of $\tilde P(A)$,
 validating the local rule for $\Lam(D)$.
Since $BD$ is a shrink, (2c1) and (2c2) follow by Remark~\ref{rem:ACshrink}.
Then (2c3) follows since the $k$'s in the 1st through $c$th columns
 do not move in the operation $\tilde P(A)\StandBereleinsert\bet_s$,
 and $\tilde P(B)$ has a $k$ at $X$ as well.
Comparing $\tilde P(D)$ with $\tilde P(B)$,
 the $k$'s in the 1st through $c$th columns shifted to the left,
 since they are not affected by insertion of $\bet_s$,
 and $k_v$ at $X$ in $\tilde P(B)$ also shifts to $Y$ in $\tilde P(D)$.
This shows (2c4) for $BD$.

\begin{figure}
\caption{The case (\B J), horizontal domino}
\label{fig:pfbjh}
\def\cell#1{\vbox to 1.5em{\vss\hbox to 1.5em{\hss$#1$\hss}\vss}}
\def\hcell#1{\vbox{\hbox to 1.5em{\hss$#1$\hss}}}
\begin{alignat*}{2}
&
\vbox{\packed\halign{&#\cr
 &&\hcell{Y}&&\hcell{X}&&\cell{\revddots}\cr
        &\multispan5\hrulefill\cr
 \cell{}&\vrule&\cell{k_u}&\vrule&\cell{\emp}&\vrule&\cell{}\cr
        &\multispan5\hrulefill\cr
 \cell{\revddots}\cr
 \multispan7\hfil$\tilde P(A)$\hfil\cr
}}
&&
\vbox{\packed\halign{&#\cr
 &&&&&&\cell{\revddots}\cr
        &\multispan5\hrulefill\cr
 \cell{}&\vrule&\cell{k_u}&\vrule&\cell{k_v}&\vrule&\cell{}\cr
        &\multispan5\hrulefill\cr
 \cell{\revddots}\cr
 \multispan7\hfil$\tilde P(B)$\hfil\cr
}}
\\
\vspace{\smallskipamount}
&
\vbox{\packed\halign{&#\cr
 &&&&&&\cell{\revddots}\cr
        &\multispan5\hrulefill\cr
 \cell{}&\vrule&\cell{\emp}&\vrule&\cell{\emp}&\vrule&\cell{}\cr
        &\multispan5\hrulefill\cr
 \cell{\revddots}\cr
 \multispan7\hfil$\tilde P(C)$\hfil\cr
}}
&&
\vbox{\packed\halign{&#\cr
 &&&&&&\cell{\revddots}\cr
        &\multispan5\hrulefill\cr
 \cell{}&\vrule&\cell{k_v}&\vrule&\cell{\emp}&\vrule&\cell{}\cr
        &\multispan5\hrulefill\cr
 \cell{\revddots}\cr
 \multispan7\hfil$\tilde P(D)$\hfil\cr
}}
\end{alignat*}
\end{figure}

Next assume that $X$ is below $Y$.  
This can only happen if the bumping procedure in inserting $\bet_s$ to
 $\tilde P(A)$ bumps $k_u$ at $Y$ (row $r$) into row $r+1$,
 which happens to have been shorter than row $r$ by one cell.
If one inserts $\bet_s$ into $\tilde P(C)$ instead of $\tilde P(A)$,
 the letter which would bump $k_u$ at $Y$ for $\tilde P(A)$ just stays at $Y$,
 forming the end point of the bumping path.
Therefore the shape of $\tilde P(D)$ is the same as that of $\tilde P(A)$,
 validating the local rule for $\Lam(D)$.
Since $BD$ is a shrink, (2c1) and (2c2) follow by Remark~\ref{rem:ACshrink}.
Then (2c3) for $BD$ follows from that for $AC$ since the shrinking cells
 of $AC$ ($Y$) and $BD$ ($X$) are both in the same column (column $c$),
 and the columns to the left of this are intact.
Finally (2c4) follows from that of $AC$ for the same reason.

\begin{figure}
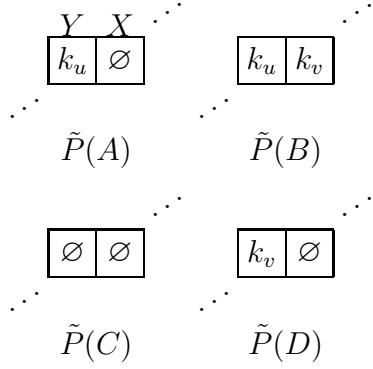

\caption{The case (\B J), vertical domino}
\label{fig:pfbjv}
\def\cell#1{\vbox to 1.5em{\vss\hbox to 1.5em{\hss$#1$\hss}\vss}}
\def\vdompic#1#2#3{%
\vbox{\packed\halign{&##\cr
 &&&&\cell{\revddots}\cr
        &\multispan3\hrulefill\cr
        &\vrule&\cell{#1}&\vrule\cr
        &\multispan3\hrulefill\cr
        &\vrule&\cell{#2}&\vrule\cr
        &\multispan3\hrulefill\cr
 \cell{\revddots}\cr
 \multispan5\hfil$#3$\hfil\cr
}}}
\begin{alignat*}{2}
\vbox{\packed\halign{&#\cr
 &&&&\cell{\revddots}\cr
        &\multispan3\hrulefill\cr
\cell{Y}&\vrule&\cell{k_u}&\vrule\cr
        &\multispan3\hrulefill\cr
\cell{X}&\vrule&\cell{\emp}&\vrule\cr
        &\multispan3\hrulefill\cr
 \cell{\revddots}\cr
 \multispan5\hfil$\tilde P( A )$\hfil\cr
}}
&&
\vdompic{\bullet}{k_u}{\tilde P(B)}
\\
\vspace{\smallskipamount}
&
\vdompic{\emp}{\emp}{\tilde P(C)}
&&
\vdompic{\bullet}{\emp}{\tilde P(D)}
\end{alignat*}
\end{figure}

{\sc The case (J)}.  Let $X$ be as in the previous case.  $X$ is the end
point of the bumping path involved in the operation $\tilde
P(A)\StandBereleinsert\bet_s$, and $Y$ is the end point of the shift of
$k$'s observed in comparison of $\tilde P(A)$ and $\tilde P(C)$.  If $X$
is in a row above $Y$, then it is clear that the shift of $k$'s does not
affect the bumping procedure for inserting $\bet_s$.  If $X$ is in a row
below $Y$, then the bumping path intersects with the shifting path at
only one cell.  This cell must contain the leftmost occurrence of $k$ in
that row by definition of the bumping procedure.  All columns to the
left are longer (since each one by assumption has a $k$, whose position
is lower because the $k$'s form a horizontal strip), and all their
contents are $\leq k$ (involving $k$'s with smaller indices only).  This
forces the bumped $k$ to sit in the empty cell immediately below the
intersection, marking the end of the bumping path.   
 
Now the intersection cannot occur at $Y$, since 
if it did, then $X$ must have been its lower adjacent,
which would put us in case (\B J) rather than (J).  
Now if $\bet_s$ is inserted to $\tilde P(C)$, the bumping procedure
 is the same until it comes to the intersection.
At the intersection a different $k$ is bumped, and it sits straight below
at $X$.
Therefore $\Lam( D ) = \Lam( C ) \cup X$, as prescribed in the local rule,
and $BD$ is a shrink, so the latter half of (2) for $BD$ consists of
(2c1)--(2c4).
Now (2c1) and (2c2) follow from Remark~\ref{rem:stratum}; further
(2c3) follows since $\tilde P(A)$ has the same property and the difference
 in $\tilde P(C)$ concerning these $k$'s is just one vertical movement
 of one of the $k$'s.
Finally, (2c4) also follows from the above description.  

\begin{figure}
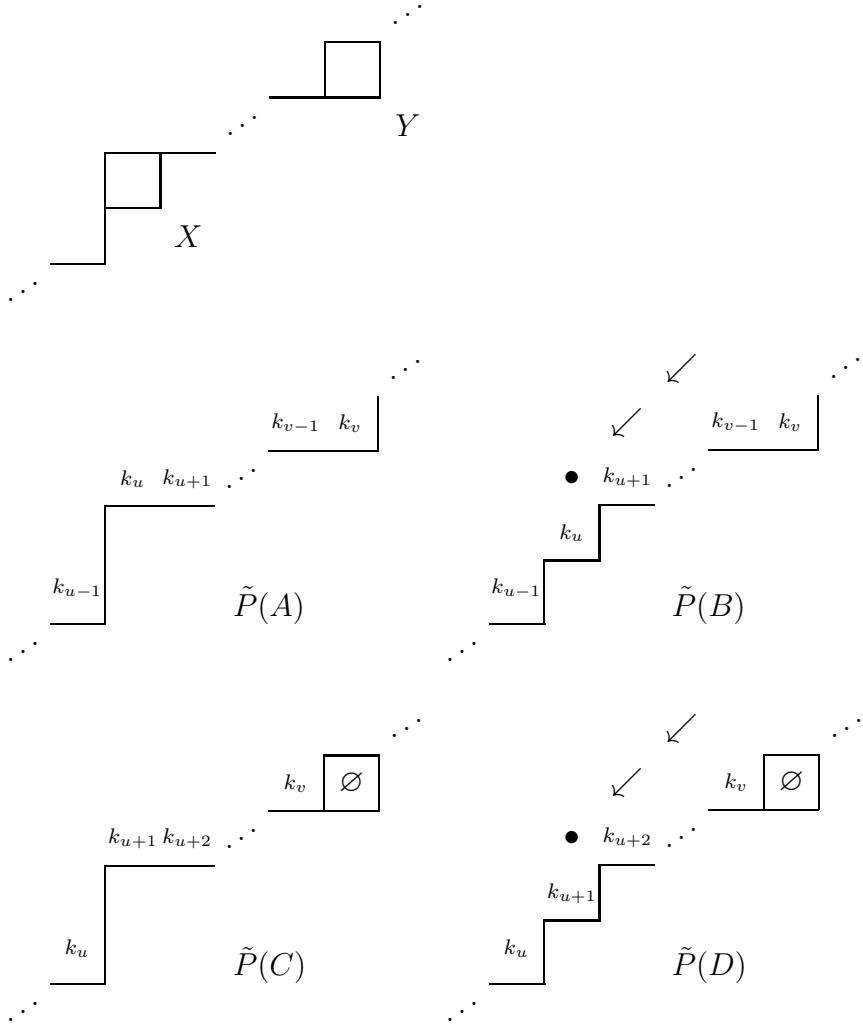

\caption{The case (J)}
\label{fig:pfj}
\def\cell#1{\vbox to 1.75em{\vss\hbox to 1.75em{\hss$#1$\hss}\vss}}
\allowdisplaybreaks
\begin{alignat*}{2}
&
\vbox{\packed\halign{&#\cr
 &&&&&&&&&&&\cell{\revddots}\cr
 &&&&&&&&\multispan3\hrulefill\cr
 &&&&&&&\cell{}&\vrule&\cell{}&\vrule\cr
 &&&&&&&\multispan4\hrulefill\cr
 &&&\cell{}&&\cell{}&\cell{\revddots}&&&&&\cell{Y}\cr
 &&\multispan4\hrulefill\cr
 &\cell{}&\vrule&\cell{}\vrule\cr
 &&\multispan3\hrulefill\cr
 &\cell{}&\vrule&&&\cell{X}\cr
 &\multispan2\hrulefill\cr
 \cell{\revddots}\cr
}}
&&
\\
\vspace{\medskipamount}
&
\vbox{\packed\halign{&#\cr
 &&&&&&&&&&&\cell{\revddots}\cr
 &&&&&&&\cell{\ssize k_{v-1}}&&\cell{\ssize k_v}&\vrule\cr
 &&&&&&&\multispan4\hrulefill\cr
 &&&\cell{\ssize k_u}&&\cell{\ssize k_{u+1}}&\cell{\revddots}\cr
 &&\multispan4\hrulefill\cr
 &\cell{}&\vrule\cr
 &\cell{\ssize k_{u-1}}&\vrule&\multispan9\hfil$\tilde P(A)$\hfil\cr
 &\multispan2\hrulefill\cr
 \cell{\revddots}\cr
}}
&&
\vbox{\packed\halign{&#\cr
 &&&&&&\cell{\swarrow}&&&&&\cell{\revddots}\cr
 &&&&&\cell{\swarrow}&&\cell{\ssize k_{v-1}}&&\cell{\ssize k_v}&\vrule\cr
 &&&&&&&\multispan4\hrulefill\cr
 &&&\cell{\bullet}&&\cell{\ssize k_{u+1}}&\cell{\revddots}\cr
 &&&&\multispan2\hrulefill\cr
 &\cell{}&&\cell{\ssize k_u}&\vrule\cr
 &&\multispan3\hrulefill\cr
 &\cell{\ssize k_{u-1}}&\vrule&\multispan9\hfil$\tilde P(B)$\hfil\cr
 &\multispan2\hrulefill\cr
 \cell{\revddots}\cr
}}
\\
\vspace{\smallskipamount}
&
\vbox{\packed\halign{&#\cr
 &&&&&&&&&&&\cell{\revddots}\cr
 &&&&&&&&\multispan3\hrulefill\cr
 &&&&&&&\cell{\ssize k_v}&\vrule&\cell{\emp}&\vrule\cr
 &&&&&&&\multispan4\hrulefill\cr
 &&&\cell{\ssize k_{u+1}}&&\cell{\ssize k_{u+2}}&\cell{\revddots}\cr
 &&\multispan4\hrulefill\cr
 &\cell{}&\vrule\cr
 &\cell{\ssize k_u}&\vrule&\multispan9\hfil$\tilde P(C)$\hfil\cr
 &\multispan2\hrulefill\cr
 \cell{\revddots}\cr
}}
&&
\vbox{\packed\halign{&#\cr
 &&&&&&\cell{\swarrow}&&\multispan3\hrulefill&\cell{\revddots}\cr
 &&&&&\cell{\swarrow}&&\cell{\ssize k_v}&\vrule&\cell{\emp}&\vrule\cr
 &&&&&&&\multispan3\hrulefill\cr
 &&&\cell{\bullet}&&\cell{\ssize k_{u+2}}&\cell{\revddots}\cr
 &&&&\multispan2\hrulefill\cr
 &\cell{}&&\cell{\ssize k_{u+1}}&\vrule\cr
 &&\multispan3\hrulefill\cr
 &\cell{\ssize k_u}&\vrule&\multispan9\hfil$\tilde P(D)$\hfil\cr
 &\multispan2\hrulefill\cr
 \cell{\revddots}\cr
}}
\end{alignat*}
\end{figure}

{\bf The reverse R-S group.  }\label{subsec:revrsg}
{\sc The case (W)}.
The bumping path involved in the insertion of $\bet_s$ into $\tilde P(A)$
 does not intersect the shifting path, since if it did,
 the insertion would end in an ordinary insertion after bumping one of
the $k$'s in the shifting path (as explained in case (J)),
 contradicting the assumption that $\Lam(A)\cover\Lam(B)$.
So if one inserts $\bet_s$ into $\tilde P(C)$ instead of $\tilde P(A)$,
 the bumping phase is exactly the same.
The assumption of the case (W) is that
 the end point of the sliding path involved in $\tilde P(A)\StandBereleinsert
 \bet_s$, which we again call $X$, does not coincide with $Y$.

First assume that $X$ is higher than $Y$ in $\Lam( A )$.
In this case the whole sliding path is higher than the shifting path,
 and there is no intersection.
We claim that the insertion of $\bet_s$ into $\tilde P( C )$ causes
 the same sliding as that into $\tilde P( A )$.
Look at each step of the sliding.
If the hole is at least two rows above the shifting path, clearly the sliding
 occurs in the same direction regardless of the shift.
Assume that the hole is immediately above the shifting path.
Because of their positions, the sliding path must be purely horizontal
 after this point to the end of the row, and if the shifting path has
more squares to the right of the column of the hole,
 then that portion of the shifing path also must be purely horizontal.
We distinguish two cases:
(1)
If the shifting path continues more to the right, then the right adjacent of
 the hole contains a letter strictly smaller than $k$.
If one shifts the $k$'s in the shifting path (to make it into $\tilde P( C )$),
 the cell below the hole still contains some $k$,
 so the sliding in $\tilde P(C)$ proceeds to the right, as desired.
(2)
If the shifting path ends exactly at the cell below the hole,
 then in $\tilde P( C )$ that place is vacant.
There is no comparison involved in this case, and the sliding again proceeds to the
 right as desired.
In these cases, the remaining claims are easy, so we omit the details.

On the other hand, if $X$ is lower, then by Lemma~\ref{LemYa} the sliding path reaches
 $X$ from above, and $X$ is the only intersection of these two paths.
Let $k_v$ denote the contents of $X$ in $\tilde P( A )$.
Since $Y$ is higher than $X$, we must have at least $k_{ v + 1 }$
 in the shifting path, which must be in the column immediately to the right
 of $X$.
Now we again claim that the insertion of $\bet_s$ into $\tilde P( C )$
 causes the same sliding as it causes to  $\tilde P( A )$,
 except that the contents of the end point is a different $k$.
Look at each step of the sliding.
So long as none of the candidates to slide into the hole belongs to
 the shifting path, the sliding is not affected by the shift.
Suppose the hole is the left adjacent of $k_{ v + 1 }$.
The hole is straight above $X$, so in $\tilde P( A ) \StandBereleinsert
 \bet_s$ the sliding proceeded below.
After the shift, the lower adjacent of the hole is either $k_{ v + 1 }$
 (when the lower adjacent is $X$) or strictly smaller than $k$ (when the hole
 is at least two rows above $X$).
The right adjacent of the hole is either $k_{ v + 2 }$ (when $Y$ is at least
 two cells apart from the hole) or vacant (when the right adjacent is $Y$).
In either case, the sliding proceeds below in $\tilde P( C ) \StandBereleinsert
 \bet_s$.
Finally assume that the hole is immediately above $X$ (just before
sliding ends).
The case where the right adjacent of the hole is $k$ was discussed just above,
 and otherwise the right adjacent is vacant so that the sliding must land
 into $X$ regardless of whether it is for $\tilde P( A )$ or $\tilde P( C )$.
In summary, the location of the sliding path for $\tilde P( C )
 \StandBereleinsert \bet_s$ is exactly the same as that for $\tilde P( A )
 \StandBereleinsert \bet_s$, whose end point contains $k_{ v + 1 }$ instead
 of $k_v$.
The tableau thus obtained is the same tableau as one obtains
 from $\tilde P( B ) = \tilde P( A ) \StandBereleinsert \bet_s$
 by the procedure described in~(2c4).
Now the remaining claims follow easily.

{\sc The case ({\cyr Ya})}.
Again by Lemma~\ref{LemYa},
 the bumping path and the sliding path do not intersect
 the shifting path before reaching the cell $X=Y$, so up to that point
insertion of $\bet_s$ 
 to $\tilde P(C)$ proceeds in the same manner as that to $\tilde P(A)$.
Now the sliding for $\tilde P(A)\StandBereleinsert\bet_s$
 involves the cell $(r-1,c)$.
If $(r-1,c+1)$ is empty, then the sliding in $\tilde P(C)\StandBereleinsert
 \bet_s$ ends at $(r-1,c)$, and the difference between $\tilde P(D)$ and
 $\tilde P(B)$ is, as expected, the shifting of $k$'s at the bottoms of
 1st through $c$th columns.
If $(r-1,c+1)$ is not empty
 (see Figure~\ref{fig:pfya} below, where $\nwarrow$ signifies that the
difference with 
 the left tableau is a sequence of sliding, and a bumping preceding it,
 which is not mentioned in the picture),
 then this cell as well as all cells in this row to the right must contain
 something $\ge k$, because the sliding in $\tilde P( A ) \StandBereleinsert
 \bet_s$ proceeded to $X=Y$ which contained $k_u$.
Since no letters $\ge \bar k$ can appear in $\tilde P(A)$, they are all
 $k$'s, and remain the same in $\tilde P(C)$,
 and in $\tilde P(C)\StandBereleinsert\bet_s$ the sliding proceeds
 towards the end of row $r-1$.
In this case the difference of $\tilde P(D)$ from $\tilde P(B)$ spreads
 to the end of row $r - 1$, where the difference is also the shifting
 of $k$'s at the bottom cells.
Thus all claims now follow for this case.

\begin{figure}
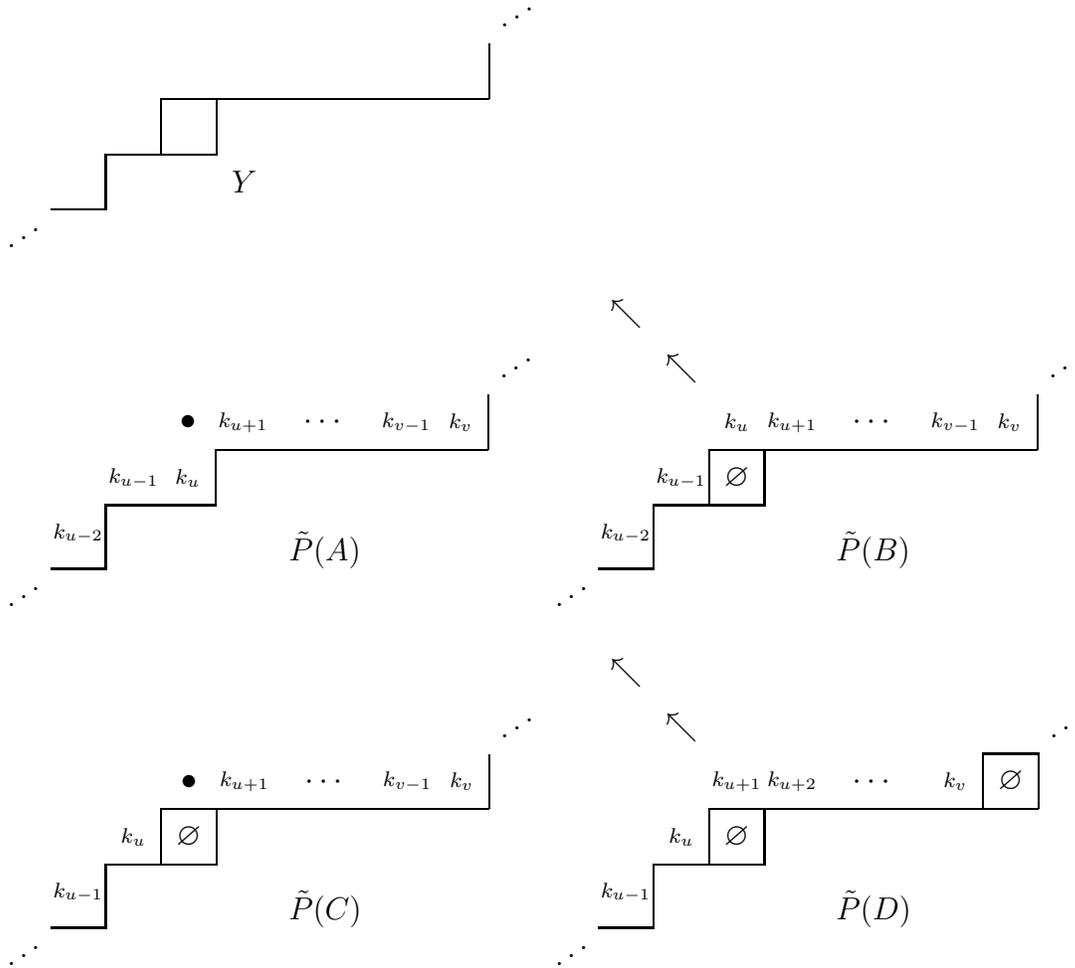

\caption{The case ({\cyr Ya})}
\label{fig:pfya}
\def\cell#1{\vbox to 1.75em{\vss\hbox to 1.75em{\hss$#1$\hss}\vss}}
\allowdisplaybreaks
\begin{alignat*}{2}
&\vbox{\packed\halign{&#\cr
 &\cell{}\cr
 &&&\cell{}&&&&&\cell{}&\cell{}&&&&&\cell{\revddots}\cr
 &&&&&\cell{}&&\cell{}&\multispan2\hfil\cell{}\hfil&
  \cell{}&&\cell{}&\vrule\cr
 &&&&\multispan9\hrulefill\cr
 &&&\cell{}&\vrule&\cell{}&\vrule\cr
 &&\multispan5\hrulefill\cr
 &\cell{}&\vrule&\cell{}&&\cell{}&&\cell{Y}\cr
 &\multispan2\hrulefill\cr
 \cell{\revddots}\cr
}}
&&
\\
&
\vbox{\packed\halign{&#\cr
 &\cell{}\cr
 &&&\cell{}&&&&&\cell{}&\cell{}&&&&&\cell{\revddots}\cr
 &&&&&\cell{\bullet}&&\cell{\ssize k_{u+1}}&\multispan2\hfil\cell{\cdots}\hfil&
  \cell{\ssize k_{v-1}}&&\cell{\ssize k_v}&\vrule\cr
 &&&&&&\multispan7\hrulefill\cr
 &&&\cell{\ssize k_{u-1}}&&\cell{\ssize k_u}&\vrule\cr
 &&\multispan5\hrulefill\cr
 &\cell{\ssize k_{u-2}}&\vrule&\multispan{12}\hfil$\tilde P(A)$\hfil\cr
 &\multispan2\hrulefill\cr
 \cell{\revddots}\cr
}}
&&
\vbox{\packed\halign{&#\cr
 &\cell{\nwarrow}\cr
 &&&\cell{\nwarrow}&&&&&\cell{}&\cell{}&&&&&\cell{\revddots}\cr
 &&&&&\cell{\ssize k_u}&&\cell{\ssize k_{u+1}}&\multispan2\hfil\cell{\cdots}\hfil&
  \cell{\ssize k_{v-1}}&&\cell{\ssize k_v}&\vrule\cr
 &&&&\multispan9\hrulefill\cr
 &&&\cell{\ssize k_{u-1}}&\vrule&\cell{\emp}&\vrule\cr
 &&\multispan5\hrulefill\cr
 &\cell{\ssize k_{u-2}}&\vrule&\multispan{12}\hfil$\tilde P(B)$\hfil\cr
 &\multispan2\hrulefill\cr
 \cell{\revddots}\cr
}}
\\
\vspace{\smallskipamount}
&
\vbox{\packed\halign{&#\cr
 &\cell{}\cr
 &&&\cell{}&&&&&\cell{}&\cell{}&&&&&\cell{\revddots}\cr
 &&&&&\cell{\bullet}&&\cell{\ssize k_{u+1}}&\multispan2\hfil\cell{\cdots}\hfil&
  \cell{\ssize k_{v-1}}&&\cell{\ssize k_v}&\vrule\cr
 &&&&\multispan9\hrulefill\cr
 &&&\cell{\ssize k_{u}}&\vrule&\cell{\emp}&\vrule\cr
 &&\multispan5\hrulefill\cr
 &\cell{\ssize k_{u-1}}&\vrule&\multispan{12}\hfil$\tilde P(C)$\hfil\cr
 &\multispan2\hrulefill\cr
 \cell{\revddots}\cr
}}
&&
\vbox{\packed\halign{&#\cr
 &\cell{\nwarrow}\cr
 &&&\cell{\nwarrow}&&&&&\cell{}&\cell{}&&\multispan3\hrulefill&\cell{\revddots}\cr
 &&&&&\cell{\ssize k_{u+1}}&&\cell{\ssize k_{u+2}}&\multispan2\hfil\cell{\cdots}\hfil&
  \cell{\ssize k_v}&\vrule&\cell{\emp}&\vrule\cr
 &&&&\multispan9\hrulefill\cr
 &&&\cell{\ssize k_u}&\vrule&\cell{\emp}&\vrule\cr
 &&\multispan5\hrulefill\cr
 &\cell{\ssize k_{u-1}}&\vrule&\multispan{12}\hfil$\tilde P(D)$\hfil\cr
 &\multispan2\hrulefill\cr
 \cell{\revddots}\cr
}}
\end{alignat*}
\end{figure}

\noindent(End of Proof of Lemma~\ref{LemLRStep})

\subsection{Concluding the proof of Theorem~\ref{ThmLocRules}}\label{subsec:con}
Putting $I=L$ in Lemma~\ref{LemLocRules}, we have actually already shown (1).
(2) is also contained in the first part of the proof of Lemma~\ref{LemLocRules}.
Note that the irrelevance of the order of application is automatic
 since we have proved that $\Lam(D)$ is always equal to the shape of
 the $P$-symbol of $\tilde w(D)$, which is determined independently of
 the local rules.
To see (3), we just have to look again at the local rules
 to check that the cases are disjoint with respect to the shapes $\Lam(B)$,
 $\Lam(C)$, $\Lam(D)$ and the stratum, which is easy.
Note that we did not claim that we can start from any sequences of shapes
 along the rightmost edge and the bottom edge; what we show here is that
 if an array of shapes is obtained from a $w$, then it can be recovered
 from the rightmost and bottom edges.
No statement as to ``coverage of all cases'' in (3) is necessary
 for showing this type of ``injectivity'' result.
(4) is actually just the definition of the shapes $\Lam(f,j')$, $0\le j'\le f$.
To show (5), suppose the bottom row of the $k$-stratum occurs in row
$u_k$. 
By Lemma~\ref{LemStand}, in $P(u_k,f)$, the letters $k_1$, $k_2$, \dots,
$k_{m_w(k)}$ 
 appear from the left to the right forming a horizontal strip.
By Lemma~\ref{LemLocRules} (2c1), shrinks do not occur in the $k$-stratum,
 and by Lemma~\ref{LemLocRules} (2b)
 the difference of the growth from $(i'-1,f)$ to $(i',f)$
 (where $i'$ is in the $k$-stratum) is the position of the corresponding
 subscripted $k$ in $P(i,f)$.
Hence (5) follows.
Next, if $A$, $B$, $C$ are three vertically contiguous vertices
 on the rightmost side of a $\bar k$-stratum,
 it cannot happen that $AB$ is a growth and $BC$ is a shrink since,
 by Lemma~\ref{LemLocRules} (2b), $\tilde P(B)$ contains a $\bar k$
 and this would violate Lemma~\ref{LemLocRules} (2c2).
Therefore the rightmost side of a $\bar k$-stratum
 consists of a series of shrinks first, followed by a series of growths.
By Lemma~\ref{LemLocRules} (2c4), these shrinks 
are a part of the horizontal strip gained in the $k$-stratum.
Also by Lemma~\ref{LemLocRules} (2c3),
 these shrinks must happen from right to left.
Hence the initial part of (6) follows.
The last part of (6) also follows.
The part concerning the growth part in the $\bar k$-stratum
 follows by the same argument as that for (5).
Finally (7) follows by combining (5) and (6).
\qed

%%My additions
%%BEGIN FILE newyoungmacros.tex
%%%%%%%%%%%% NEWYOUNGMACROS.TEX
%%This one modified to entertain the other orientation as well.                                  
\newdimen\Squaresize \Squaresize=20pt
\newdimen\thickness \thickness=1pt         
                                                    
\def\Square#1{\hbox{\vrule width \thickness
   \vbox to \Squaresize{\hrule height \thickness\vss                                  \hbox to \Squaresize{\hss#1\hss}
   \vss\hrule height\thickness} 
\unskip\vrule width \thickness} 
\kern-\thickness}                                                            
                               
\def\vsquare#1{\vbox{\Square{$#1$}}\kern-\thickness}
\def\blank{\omit\hskip\Squaresize}

\def\young#1{\let\\=\cr	%added  let\\ =\cr 
\vbox{\smallskip\offinterlineskip
\halign{&\vsquare{##}\cr #1}}}

%%For making smaller young diagrams:
\def\smyoung#1{\let\\=\cr\Squaresize=15pt
\vbox{\smallskip\offinterlineskip
\halign{&\vsquare{##}\cr #1}}}

%For making diagrams for the Young-Fibonacci lattice
\def\fibyoung#1{\let\\=\cr		%added  let\\ =\cr 
\vbox{\smallskip\offinterlineskip
\halign{&\vsquare{##}\cr #1}}\,}

%%%%  example  %%%%%
%%\blank means ``don't cellify''
%$$\young{\blank &\blank &\blank &\blank & & & \cr
%\blank &\blank &\blank & & & & & \cr
%\blank &\blank & &e &f &g & & \cr
%\blank &\blank & &\blank &\blank &5 & \cr
%\blank &\blank & &\blank &7 & &\blank & \cr
%\blank &\blank & &3 & &33 & & \cr
%\blank &\blank & & &\blank & &\blank & \cr
%\blank &\blank & &\blank & & & & \cr
%}$$
%
%Here's something like a \strut turned on it's side (see p.353 of
 %TeXbook).  We need it so the border's look right
\newbox\vstrutbox
\setbox\vstrutbox=\hbox{\vrule height 0pt depth 0pt width 4pt}
\def\vstrut{\relax \ifmmode\copy\vstrutbox\else\unhcopy\vstrutbox\fi }
%And another one to keep the growths from running into the right hand
 %edge of the cell.  
\newbox\smvstrutbox
\setbox\smvstrutbox=\hbox{\vrule height 0pt depth 0pt width 1pt}
\def\smvstrut{\relax \ifmmode\copy\smvstrutbox\else\unhcopy\smvstrutbox\fi }

%%It turns out that doublespace.sty redefines strut, which screws up
 %%these macros by making elements which represent the growth in the
 %%upper right hand corners move too far down.  Hence, we recopy the
 %%original definition of \strut from the TeXbook, but we call it
 %%something else, just in case we actually want to use the dblespace
 %%modified strut later (unlikely, but . . .).
\newbox\reghstrutbox
\setbox\reghstrutbox=\hbox{\vrule height 9.5pt depth 3.5pt width 0pt}
\def\reghstrut{\relax \ifmmode\copy\reghstrutbox\else\unhcopy\reghstrutbox\fi }
  %And a smaller one
\newbox\smhstrutbox
\setbox\smhstrutbox=\hbox{\vrule height 4pt depth 2pt width 0pt}
\def\smhstrut{\relax \ifmmode\copy\smhstrutbox\else\unhcopy\smhstrutbox\fi }
  %End horizontal strut redefined.  
\newdimen\fsquaresize \fsquaresize=40pt
\newdimen\fthickness 

\def\fsquare#1{\hbox{\vrule width \fthickness
   \vbox to \fsquaresize{\hrule height \fthickness
	\hbox to \fsquaresize{\hss\hstrut #1\smvstrut}
   \vss\hrule height\fthickness} 
\unskip\vrule width \fthickness} 
\kern-\fthickness}                                                            
                               
\def\fvsquare#1{\vbox{\fsquare{$#1$}}\kern-\fthickness}
%
%%TOP-DOWN VERSION

\def\dfsquare#1{\hbox{\vrule width \fthickness
   \vbox to \fsquaresize{\hrule height \fthickness
  \vfill \hbox to \fsquaresize{\hfill\hstrut \hbox{\cornerfont #1}\smvstrut}
   \hrule height\fthickness} 
\unskip\vrule width \fthickness} 
\kern-\fthickness}                                                            
                               
\def\dfvsquare#1{\vbox{\dfsquare{#1}}\kern-\fthickness}
%

%macros to make the border, just omit rules?
\def\fborder#1{\hbox{\vrule width 0pt
   \vbox to \fsquaresize{\hrule height 0pt
	\hbox to \fsquaresize{\hss \hstrut #1\vstrut }	
   \vss\hrule height 0pt} 		%had to mess with kerns;
\unskip\vrule width 0pt} 
\kern-4\fthickness} 
                                                         
\def\ftopborder#1{\hbox{\vrule width 0pt
   \vbox to \fsquaresize{\hrule  height 0pt\vss
	\hbox to \fsquaresize{ \hstrut #1\smvstrut \hss }	
   \hrule height 0pt} 		%had to mess with kerns;
\unskip\vrule width 0pt} 
\kern-\fthickness}

\def\fvborder#1{\vbox{\fborder{$#1$}}\kern-4\fthickness}

\def\fvtopborder#1{\vbox{\ftopborder{$#1$}}\kern-4\fthickness}
%for vertical centering along the right side
\def\fvcborder#1{\hbox{\vrule width 0pt
   \vbox to \fsquaresize{\hrule  height 0pt\vss
	\hbox to \fsquaresize{ \hstrut #1\vstrut \hss }	
   \vss\hrule height 0pt} 		%had to mess with kerns;
\unskip\vrule width 0pt} 
\kern-\fthickness}

\def\fvborder#1{\vbox{\fborder{$#1$}}\kern-4\fthickness}

\def\fvvcborder#1{\vbox{\fvcborder{$#1$}}\kern-4\fthickness}
%End border macros
%
                             
%Macros to put an ``X'' in the square as well (or even better, a
 %number).  Sanserif font looks better than the alternatives.
\font\cellfontfive=cmssbx10 scaled \magstep5
\font\cellfontfour=cmssbx10 scaled \magstep4
\font\cellfontthree=cmssbx10 scaled \magstep3
\font\cellfonttwo=cmssbx10  scaled \magstep2
\font\cellfontone=cmssbx10  scaled \magstep1
\font\cellfontzero=cmssbx10  
%\font\cellfontmone=cmssbx8
\font\sevenrm=cmr7
\font\eightrm=cmr8
\font\sixrm=cmr6
\font\fiverm=cmr5
\font\sevenit=cmti7
\font\eightit=cmti8
\def\cellfont{\cellfontfive}

\def\xsquare#1#2{\hbox{\vrule width \fthickness
   \vbox to \fsquaresize{\hrule height \fthickness
	\vbox {
	   \hbox to \fsquaresize {\hss\hstrut {\cornerfont #2}\smvstrut }
	   \vss}
	\vss
	   \hbox to \fsquaresize{\hss\hbox{\cellfont #1}\hstrut \hss}
   	\vss  \hrule height\fthickness
	} 
\unskip\vrule width \fthickness} 
\kern-\fthickness}                                                            
                               
\def\xvsquare#1.#2.{\vbox{\xsquare{#1}{#2}}\kern-\fthickness} %?Why $$ ?
%%Now that I am using the other orientation, I need a new macro.

\def\Xsquare#1#2{\hbox{\vrule width \fthickness
   \vbox to \fsquaresize{\hrule height \fthickness \kern -5.5\fthickness
	\hbox to \fsquaresize{\hss 
	   \vbox to \fsquaresize{\vss\hbox to 0pt{\cellfont #1\hss }\vss} 
	  \hss
   	   \vbox to \fsquaresize{\vss \hbox{\hstrut \cornerfont #2}}
	}
    \hrule height\fthickness}
\unskip\vrule width \fthickness} 
\kern-\fthickness}

\def\Xvsquare#1.#2.{\vbox{\Xsquare{#1}{#2}}\kern-\fthickness} %?Why $$ ?

%End X macros 

\def\fyoung#1{  
	\def\hstrut{\reghstrut}
	\fthickness=1pt  		%.025\fsquaresize  %%Don't need.
\def\>{\omit\fvborder}
\def\vc{\omit\fvvcborder}
\def\<{\omit\fvtopborder}
\def\x{\omit\xvsquare}	
\def\X{\omit\Xvsquare}	
\let\\=\cr
\def\blank{\omit\hskip\fsquaresize}
\vbox{\smallskip\offinterlineskip
\halign{&\fvsquare{##}\cr #1}}}	

%%THE TOP-DOWN VERSION

\def\dfyoung#1{  
	\def\hstrut{\reghstrut}
	\fthickness=1pt  		%.025\fsquaresize  %%Don't need.
\def\>{\omit\fvborder}
\def\vc{\omit\fvvcborder}
\def\<{\omit\fvtopborder}
\def\x{\omit\xvsquare}	
\def\X{\omit\Xvsquare}	
\let\\=\cr
\def\blank{\omit\hskip\fsquaresize}
\vbox{\smallskip\offinterlineskip
\halign{&\dfvsquare{##}\cr #1}}}	

%SMALLER DIAGRAMS

\def\smfyoung{\fsquaresize=30pt 
		\def\cellfont{\cellfontthree}  \def\cornerfont{\sevenrm}
				\fyoung}
  %End \smfyoung

%%SMALLER, TOP-DOWN VERSION
\def\smdfyoung{\fsquaresize=19pt 
		\def\cellfont{\cellfontzero}
		\def\borderfont{\sevenit}
		\def\cornerfont{\sevenrm}\dfyoung} 
  %End \smdfyoung

\def\textfyoung{\fsquaresize=12pt 
		\def\cellfont{\fiverm}  \fyoung}

\def\smxsquare#1#2{\hbox{\vrule width \fthickness
   \vbox to \fsquaresize{\hrule height \fthickness
	\vss	\hbox to \fsquaresize{\hss\hbox{\cellfont #1}\hstrut \hss}
   \vss  \hrule height\fthickness} 
\unskip\vrule width \fthickness} 
\kern-\fthickness}

\def\smXsquare#1#2{\hbox{\vrule width \fthickness
   \vbox to \fsquaresize{\hrule height \fthickness\kern %2\fthickness
	\vbox {\vss
	 	   \hbox to \fsquaresize{\hstrut\hss\hbox{\cellfont #1} \hss}
   	 \hbox to \fsquaresize {\hss\hstrut {\mit #2}\smvstrut }
	   }
	 %\hrule height\fthickness
	} 
\unskip\vrule width \fthickness} 
\kern-\fthickness}                                                            

\def\ssmfyoung{\fsquaresize=22pt \def\hstrut{\reghstrut}
		\def\cellfont{\cellfontone}
	\def\xsquare{\smxsquare}\def\Xsqaure{\smXsquare }
\fyoung}
%\tracingcommands=1
%Just so that I can use this name:
\def\fomin{\fyoung}

%%END FILE newyoungmacros.tex

\def\zz{\mathbb}
\def\mit{\rm}
\def\rs{\quad \stackrel{\rm R-S}{\longleftrightarrow}\quad }
\def\berele{\underset{\Cal B}\longleftrightarrow}%%  \quad seemed too big (IT)
\def\em{\sl}
\def\Bereleinsert{\underset{\Cal B}{\leftarrow}}
\def\cover{\overset{\cdot}{\supset}}
\def\covered{\overset{\cdot}{\subset}}
\def\rshrink#1{\overset{#1}{\supset}}
\def\rgrowth#1{\overset{#1}{\subset}}
\def\Cpx{\mathbb C}
\def\cross{{\times}}
\def\emp{\varnothing}
\def\gam{\gamma}
\def\Gam{\varGamma}
\def\inv{^{-1}}
\def\kap{\kappa}
\def\lam{\lambda}
\def\Lam{\varLambda}
\def\Nat{\mathbb N}
\def\Par{\Cal P}
\def\th#1{^{(#1)}}

\section{The Reverse Correspondence}\label{sec:rev}

As stated in (3) of Theorem~\ref{ThmLocRules},
if we know the up-down tableau on the bottom and the
sequence of shapes on the rightmost edge,
both properly produced from a word $w$,
we can recover the whole array of
shapes as well as the word $w$.  On the other hand, suppose we are given
a pair $(P,Q)$ of an $Sp(2n)$-tableau $P$ and an up-down tableau $Q$
produced from a word $w$ by Berele's correspondence, and we wish to recover $w$
from $P$ and $Q$ in a pictorial manner.  We can put $Q$ along
the bottom edge; however, the tableau $P$ only tells us the
shapes on the rightmost edge at certain points:
at the border between the ($\overline{k-1}$)- and $k$-strata for each $k$,
and at some point in the $\bar k$-stratum where the sequence of
shrinks turns into a sequence of growths for each $k$.
This of itself is insufficient
to determine the rightmost edge, since we do not know directly from $P$ and $Q$
 how many $k$-$\bar k$ cancellations occur in the $\bar k$-stratum,
 nor what shape should be put at the border between $k$ and $\bar
k$.  Nonetheless, Berele's correspondence
in its original form is reversible, so it should be possible to see how
to do this from our pictorial point of view.

Suppose we are given a pair $(P,Q)$ as above, which by the usual Berele
correspondence corresponds to a word $w$.  If we apply our algorithm to
this word, we obtain $Q$ along the bottom edge, and along the right edge
an up-down tableau
$T = ( \tau \th 0, \tau \th 1, \ldots, \tau \th f )$ which, together with the
stratification datum, determines $P$ in the manner described in (6) and (7)
of Theorem~\ref{ThmLocRules}.  If we use the notation in the previous section for this $w$,
then $\tau \th i = \Lam( i, f )$.  In fact,
not merely a sequence of shapes, but a sequence of
symplectic tableaux $( P( 0, f ), P( 1, f ), \ldots, P( f, f ) )$ is attached
to this edge; we put $P \th i = P( i, f )$ for simplicity.  $P \th i$ is
the symplectic tableau determined by the up-down tableau $( \tau
\th 0, \tau \th 1, \ldots, \tau \th i ) $ in the same manner.
%\begin{comment}
%(Note that this sequence of symplectic tableaux is certainly not the
%same as the sequence obtained in the usual correspondence, whose shapes
%are given by the up-down tableau $Q$.)
%\end{comment}

We claim that we can recover $P \th{ f - 1 }$ from $P$ and $Q$ as in the
 following Theorem~\ref{thm:revmain}, so that inductively we can recover the whole array
 of shapes without going back to $w$ in the first place.

\begin{rem}\label{rem:stratum}
(5) and (6) of Theorem~\ref{ThmLocRules}, together with (2c2) of
Lemma~\ref{LemLocRules} imply 
that, if one knows which stratum
 the bottom row cells of the picture belong to, then $P \th{ f - 1 }$ is
 completely determined from $P$ as follows:
\begin{enumerate}
\item If the bottom row belongs to $k$-stratum, with $k \in \Nat$, then
 $P$ contains a letter $k$, and
 $P \th{ f - 1 }$ is obtained from $P$ by removing the rightmost occurrence
 of the letter $k$.
 (The rightmost segment is a growth.)
\item If the bottom row belongs to $\bar k$-stratum, with $k \in \Nat$,
 and if $P$ contains the letter $\bar k$, then $P \th{ f - 1 }$ is obtained
 from $P$ by removing the rightmost occurrence of the letter $\bar k$.
 (The rightmost segment is a growth.)
\item If the bottom row belongs to $\bar k$-stratum, with $k \in \Nat$,
 and if $P$ does not contain the letter $\bar k$, then $P \th{ f - 1 }$ is
 obtained from $P$ by adding a letter $k$ to the bottom of the first column
 of $P$ which does not already contain $k$.
 (The rightmost segment is a shrink.)
\end{enumerate}
\end{rem}

\begin{thm}\label{thm:revmain}
Let $P$, $Q$ be as above, and put $Q = ( \kap \th 0, \kap \th 1, \dots,
 \kap \th f )$.
Put $l = \max \{\, l( \kap \th j ) \mid 0 \le j \le f \,\}$, and let $\gam =
 m$ or $\bar m$ \rom($m \in \Nat$\rom) be the largest letter in $P$.
\begin{enumerate}
\item If $m > l$, then $\gam$ is the largest letter in $w$,
 so that the bottom row of cells in the picture belongs to the $\gam$-stratum,
 and $P \th{ f - 1 }$ is obtained from $P$ by deleting the rightmost occurrence
 of $\gam$.
\end{enumerate}
\indent In the remaining cases, assume that $m \le l$.
\begin{enumerate}
\item[(2)] If $P$ contains a letter $\bar l$,
 then it is also the largest letter in $w$,
 so that the bottom row of the picture belongs to the $\bar l$-stratum,
 and $P \th{ f - 1 }$ is obtained from $P$ by deleting the rightmost occurrence
 of $\bar l$.
\item[(3)] If $P$ does not contain the letter $\bar l$,
 and if there is no $l$-shrink in $Q$
 \rom(i.e.\ $\kap \th{ j - 1 } \rshrink{ l } \kap \th j$ never occurs\rom),
 then $l$ is the largest letter in $w$ \rom(also in $P$\rom),
 so that the bottom row of the picture belongs to the $l$-stratum,
 and $P \th{ f - 1 }$ is obtained from $P$ by deleting the rightmost occurrence
 of $l$.
\end{enumerate}
\indent Now suppose that $P$ does not contain the letter $\bar l$, and that
 there is an $l$-shrink in $Q$.
Let $\bar P$ be the tableau obtained from $P$ by adding an $l$ to the bottom
 of the first column of $P$ which does not already contain $l$.
Put $Q$ on the bottom edge of a $1 \times f$ grid, and put the shape of
 $\bar P$ at the rightmost vertex of the upper level.
Work backwards by local rules from right to left so long as
 the vertical segment remains a shrink, assuming that we are in the
 $\bar l$-stratum.
\begin{enumerate}
\item[(4)]
If we reach the leftmost edge in this test, then $P \th{ f - 1 }$
 cannot be equal to $\bar P$.
In this case $l$ is the largest letter in $w$ \rom(also in $P$\rom),
 so that the bottom row of the picture must belong to the $l$-stratum,
 and $P \th{ f - 1 }$ is obtained from $P$ by deleting the rightmost
 occurrence of $l$.
\item[(5)]
If we encounter the case \rom($\bigcirc$\rom) before reaching
 the leftmost edge, then we have $P \th{ f - 1 } = \bar P$.
In this case the largest letter in $w$ is $\bar l$,
 so that the bottom row of the picture belongs to the $\bar l$-stratum.
\end{enumerate}
\end{thm}

\begin{rem} \label{rem:hyp}
In distingushing the cases (4) and (5),
what we propose to do in the above statement is to hypothetically assume
that the segment $( f - 1 , f )$--$( f, f )$ is a shrink,
and try to go backwards on the picture by one row to see if this assumption
yields an unacceptable row.

The situation in (4) is clearly a case of making the wrong assumption,
 since the leftmost edge can never be a shrink in a correct picture.

The case (5) is more delicate.
If we encounter $\bigcirc$ as in (5) and consequently
 reach a growth on the vertical segment, actually we can continue to
 apply local rules backwards to produce exactly one $\cross$ to the left
 of $\bigcirc$ and reach the leftmost edge with empty shapes on both ends
 of the vertical segment.
For let $\bar \kap \th 0$, $\bar \kap \th 1$, \dots, $\bar \kap \th f$
be the shapes obtained on the upper level.
All the
local rules we use with a growth on the right edge will yield
a growth again on the left edge, with the exception of ($\cross$), which sets
an equal left edge.
If we reach the latter case, the reverse computation of the
row is complete, since the local rules (\phantom{$\cross$}) insist that
all vertical segments further to the left are equals;
since $\kap \th 0$ is empty, $\bar \kap \th 0$ is also empty.
Now, if the right edge of the leftmost cell is still a growth,
then since $\kap \th 1 = ( 1 )$, we must have $\bar \kap \th 1 = \varnothing$.
Since $\kap \th 0$ is also $\varnothing$, the leftmost cell
falls into the case ($\cross$).
This means that we always will obtain exactly one cell in
case ($\cross$), and the shapes $\bar \kap \th 0$, $\bar \kap \th 1$,
\dots, $\bar \kap \th f$ form an acceptable row.

The above reasoning
also shows that if we hypothesize a growth along the segment $( f - 1, f
)$--$( f, f )$ that we will always get an apparently acceptable row of
the diagram---whether or not it is actually correct.  Thus, it is crucial
when in doubt to hypothesize a shrink along the right edge unless it
yields an unacceptable row.

Consider the following example, where it only becomes apparent we made
the wrong choice after two rows have been created working backwards.
Let $P = \smallmatrix 1 & 1 \\ 2 & 2 \endsmallmatrix$, and $Q$ be
 as in the bottom of the following partial picture (in which $f = 8$).
These numbers indicate the parts of partitions, not tableaux.
If one assumes that $\Lam( 7, 8 ) = ( 2, 1 )$ and $P \th 7 = \smallmatrix 1 & 1
 \\ 2 \endsmallmatrix$ so that the segment $( 7, 8 )$--$( 8, 8 )$ is a growth,
 one obtains the shapes on the $7$th row as below.
To go up one more level, one is forced to take $\Lam( 6, 8 ) = ( 2 )$ and
 $P \th 6 = \smallmatrix 1 & 1 \endsmallmatrix$,
 since by (6) of Theorem~\ref{ThmLocRules} a growth cannot follow a shrink
 within the $2$-stratum.
This gives the shapes on the $6$th row as indicated.
However, this is impossible, because the next row of cells above must belong
to $\bar 1$- or $1$-stratum, whereas one of the shapes has more than one
part:
$$
\setcounter{MaxMatrixCols}{30}
\def\nowidth#1{\hbox to0pt{\hss$#1$\hss}}
\newdimen\digitwidth
\setbox0=\hbox{1}
\digitwidth=\wd0
\def\thin#1{{\hbox to\digitwidth{\hss$#1$\hss}}}
\begin{matrix}
\text{row 6} \rightarrow &
 \varnothing && 1 && 1 && 1 && 2 && \thin{21} && \thin{31} && 3 && 2 && \\
&& && \nowidth \cross && && && && && && & \\
\text{row 7} \rightarrow &
 \varnothing && 1 && 2 && 2 && \thin{21} && \thin{22} && \thin{32} &&
  \thin{31} && \thin{21} && \\
&& && && \nowidth \cross && && && && && & \\
\text{row 8} \rightarrow &
 \varnothing && 1 && 2 && 3 && \thin{31} && \thin{32} && \thin{33} &&
  \thin{32} && \thin{22} &&
\end{matrix}
$$

The correct picture corresponding to this pair of tableaux
 is shown after the proof of Theorem~\ref{thm:revmain} is complete.
\end{rem}

\begin{proof}
First note that saying that $\gam$ is the largest letter in $w$ is the same
 thing as saying that the bottom row of the picture of $w$ belongs to the
 $\gam$-stratum.

We make a couple of easy points.

\begin{lemma}\label{lem:revone}
If the largest letter in $w$ survives in $P$,
 then $P \th{ f - 1 }$ is obtained from $P$ by deleting
 the rightmost occurrence of $\gam$.
\end{lemma}

\begin{proof}
This follows immediately from Remark~\ref{rem:stratum}.  
\end{proof}

\begin{lemma}\label{lem:revtwo}
$w$ contains at least one letter $\geq l$.
\end{lemma}

\begin{proof}
By the definition of symplectic tableau, we see that $P( f, j )$
 (where $j$ is such that $l( \kap \th j ) = l$) contains some letter $\geq l$.
Since $P( f, j )$ contains a subset of the letters in $w$,
 the claim follows.
\end{proof}

\begin{lemma}\label{lem:revthree}
Any letter $m'$ or $\overline{ m' }$ with $m' > l$ which
 appears in $w$ also appears in $P$.
\end{lemma}

\begin{proof}
$m'$ or $\overline{ m' }$ could only have been canceled
 {\it after\/} being pushed into the $m'$th row.
\end{proof}

Now we come back to the case-by-case analysis of Theorem~\ref{thm:revmain}.

(1) Since $m>l$, combining Lemmas~\ref{lem:revtwo} and
\ref{lem:revthree}, we conclude that the largest letter in $w$ is $\gam$.
The rest follows from Lemma~\ref{lem:revone}.    

Note that in all remaining cases, no letters $> \bar l$ can appear in $w$
by Lemma~\ref{lem:revthree}, so that the largest letter in $w$ is either
$l$ or $\bar l$. 

(2)
In this case $P$ contains $\bar l$, so it is the largest letter in $w$.
The rest follows by Lemma~\ref{lem:revone}.

(3)
If $P$ did not contain $l$, then the largest letter in $P$ would be
 less than $l$, so that $l( \kap \th f ) < l$.
However, this would necessitate an $l$-shrink on the part of $Q$,
 contradicting our assumption.
Therefore $P$ contains $l$.
If $w$ contained $\bar l$, then there must have been an $l$-$\bar l$
 cancellation, which would again cause an $l$-shrink,
 contradicting our assumption.
Therefore $l$ is also the largest letter in $w$.
The rest follows by Lemma~\ref{lem:revone}.

(4) 
This hypothetical row in the picture is clearly impossible since the
leftmost edge must be identically zero (empty).  By
Remark~\ref{rem:stratum} this eliminates the possibility of the largest
letter in $w$ being $\bar l$, so that the largest letter must be $l$.
The rest follows by Lemma~\ref{lem:revone}.

(5)
By an argument presented in Remark~\ref{rem:hyp}, 
 continuing backwards by local rules in this hypothetical row of cells
 produces exactly one $\cross$ to the left of the circle,
 and reaches the leftmost edge with empty shapes on both sides of the
 vertical edge.
Let $\bar \kap \th 0$, $\bar \kap \th 1$, \dots, $\bar \kap \th f$ be the
 shapes obtained on the upper level.
They determine an up-down tableau
$\bar Q$ if we identify the unique pair of consecutive identical
shapes $\bar \kap \th{ j - 1 } = \bar \kap \th j$
just above the cell labeled with $\times $.
Let $v$ be the word corresponding to the pair $( \bar P, \bar Q )$
by Berele's correspondence.
Our pictorial procedure applied to $v$ yields a valid $( f - 1 ) \times
( f - 1 )$ diagram.
We can then paste these shapes onto the $( f - 1 ) \times f$ diagram
by simply duplicating column $j - 1$.
All we need to show is
that this $( f - 1 ) \times f$ diagram and our hypothetical row fit
together properly to form a valid $f \times f$ diagram.  Since all cells
follow local rules, what this actually
means is that there is no inversion of strata between the
picture of $v$ and the hypothetical row, and that if the bottom row of $v$
 and the hypothetical row belongs to the same stratum, then
 the positions of $\cross$ in those two rows are in increasing order.
If these conditions are satisfied, we can define $w$ by putting $\bar l$
 between the ($j - 1$)st and the $j$th letters of $v$.
Then the synthesized $f \times f$ picture coincides with the picture of $w$,
 and the claim will be proved.

First we show the consistency of stratification.
Since the hypothetical row is produced by the local rules for the $\bar l$-%
 stratum, it is sufficient to see that the word $v$ contains
 no letters $\ge \bar l$.
Now inductively assume that the whole theorem is true for any pair $( P^\flat,
 Q^\flat )$ whose $Q$-part has order $f - 1$.
(The starting point of induction is the trivial case where $f = 0$.)
Note that the local rules for the $\bar l$-stratum, with shapes of length
 $\le l$ on the lower edge, assure that the shapes on the upper edge also have
 length $\le l$.
Therefore we have $\max \{\, l( \bar \kap \th j ) \mid 0 \le j \le f \,\} \le
 l$.
Also the largest letter in $\bar P$ is $l$ by construction.
Then the theorem applied to the pair $( \bar P, \bar Q )$ assures that
 the largest letter in $v$ $\le \bar l$.
Hence the consistency of stratification is proved.

Next suppose that row $f - 1$ of the picture of $v$ belongs to the $\bar l$-%
 stratum.
Again by the theorem applied to $( \bar P, \bar Q )$, this implies that
 $( f - 2, f - 1 )$--$( f - 1, f - 1 )$ of the picture of $v$ is a shrink
 (which is shifted to $( f - 2, f )$--$( f - 1, f )$ in the synthesized
 picture), so that row $f - 1$ of the picture of $v$ (as well as
 the synthesized picture) contains exactly
 one $\cross$ and $\bigcirc$.
We want to show that the position of $\cross$ in row $f$
 of the synthesized picture is to the right of that in row $f - 1$.

To analyze the relation between two consecutive rows in the picture,
 we will use the following two lemmas.
Consider two vertically contiguous cells, and suppose that the shapes $\lam$,
 $\mu$, $\nu$ and $\lam'$ are given. (See Figure~\ref{fig:2cells}.)

\begin{figure}
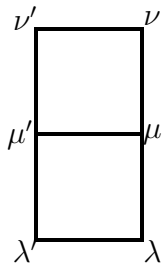

\caption{Two vertically contiguous cells}
\label{fig:2cells}
%\newdimen\squaresize\newdimen\Thickness
%\let\Young=\young%%%%%%%%%%%%%%%%%%%%%%%%%%%%%%%%%% (IT, but failed 960701)
%%%%%%%%%%%%%%%%%%%%%%%%%%%%%%%%%%%%%%%%%%%%%%%%%%% previously it was \def

%%BEGIN FILE YYoung.tex  %%
%%This is YYoung.tex, fancy macros for making young diagrams.  
%Make ! a letter so we can use it in the names of control sequences
\catcode `!=11

\newdimen\squaresize 
\newdimen\thickness 
\newdimen\Thickness
\newdimen\ll! \newdimen \uu! 
\newdimen\dd! \newdimen \rr! \newdimen \temp!

%parameters are left, up, down, right, and contents
\def\sq!#1#2#3#4#5{%
\ll!=#1 \uu!=#2 \dd!=#3 \rr!=#4
\setbox0=\hbox{%
%left edge
 \temp!=\squaresize\advance\temp! by .5\uu!
 \rlap{\kern -.5\ll! 
 \vbox{\hrule height \temp! width#1 depth .5\dd!}}%
%
%up edge
 \temp!=\squaresize\advance\temp! by -.5\uu!  
 \rlap{\raise\temp! 
 \vbox{\hrule height #2 width \squaresize}}%
%
%down edge
 \rlap{\raise -.5\dd!
 \vbox{\hrule height #3 width \squaresize}}%
%
%right edge
 \temp!=\squaresize\advance\temp! by .5\uu!
 \rlap{\kern \squaresize \kern-.5\rr! 
 \vbox{\hrule height \temp! width#4 depth .5\dd!}}%
%
%contents
 \rlap{\kern .5\squaresize\raise .5\squaresize
 \vbox to 0pt{\vss\hbox to 0pt{\hss $#5$\hss}\vss}}%
}%end of \hbox
 \ht0=0pt \dp0=0pt \box0
}%end of \sq!

\def\vsq!#1#2#3#4#5\endvsq!{\vbox to \squaresize
  {\hrule width\squaresize height 0pt%
\vss\sq!{#1}{#2}{#3}{#4}{#5}}}

\newdimen \LL! \newdimen \UU! \newdimen \DD! \newdimen \RR!

\def\vvsq!{\futurelet\next\vvvsq!}
\def\vvvsq!{\relax
  \ifx     \next l\LL!=\Thickness \let\continue!=\skipnexttoken!
  \else\ifx\next u\UU!=\Thickness \let\continue!=\skipnexttoken!
  \else\ifx\next d\DD!=\Thickness \let\continue!=\skipnexttoken!
  \else\ifx\next r\RR!=\Thickness \let\continue!=\skipnexttoken!
  \else\ifx\next P\let\continue!=\place!
  \else\def\continue!{\vsq!\LL!\UU!\DD!\RR!}%
  \fi\fi\fi\fi\fi 
  \continue!}

\def\skipnexttoken!#1{\vvsq!}

\def\place! P#1#2#3{%
\rlap{\kern.5\squaresize\temp!=.5\squaresize\kern#1\temp!
  \temp!=\squaresize \advance\temp! by #2\squaresize
 \temp!=.5\temp!  \raise\temp!\vbox to 0pt%
{\vss\hbox to 0pt{\hss$#3$\hss}\vss}}\vvsq!}

\def\Young#1{\LL!=\thickness \UU!=\thickness \DD! = \thickness \RR! = \thickness
\vbox{\smallskip\offinterlineskip
\halign{&\vvsq! ## \endvsq!\cr #1}}}

\def\blank{\omit\hskip\squaresize}
\catcode `!=12

\squaresize = 40pt
\thickness = 1pt
\Thickness = 3pt
$$ \Young{P{-1.2}{1.2}{\nu'}P{1.2}{1.2}\nu{} \cr
P{-1.3}{1.0}{\mu'}P{1.2}{1.0}\mu P{-1.2}{-1.2}{\lambda'}P{1.2}{-1.2}\lambda {}\cr}$$
\end{figure}

\begin{lemma}\label{lem:revfour}
Suppose we are given a sequence of consecutive shrinks $\nu \rshrink{r} \mu
\rshrink{s} \lambda $ with $r\leq s$ and suppose $\lambda \cover \lambda '$ or
$\lambda \covered \lambda '$.  Working backwards by local rules, suppose
that we do not encounter the case $\times $ or $\bigcirc $.  Then we obtain
$\nu' \rshrink{r'} \mu' \rshrink{s'} \lambda' $ with $r'\leq s'$.
\end{lemma}

\begin{proof}  We take several cases.  First suppose that $\lambda'
\rshrink{t}\lambda $.  If $t\neq s$ then $\mu '\rshrink{t}\mu $ and
$s'=s$.  Now if also $t\neq r$, then $r'=r$ and the lemma follows.
Otherwise, $t=r$ means that $r=t\neq s$, so $r<s $ implying that $r'\leq
r+1\leq s$.  (Note that no applicable local rule allows the row of
shrink to increase by more than one.)  Finally, if $t=s$ then $s'=s+1$,
so in any case $r'\leq r+1\leq s+1=s'$.

Now suppose that $\lambda '\rgrowth{t}\lambda  $. If $t\neq s-1$ or
$t=s-1$ but $\mu /\lambda' $ is not a vertical domino, then $s'=s$ and
$\mu ' \rgrowth{t}\mu $.  Now if also $t\neq r-1$ or $t=r-1$ but $\nu
/\mu' $ is not a vertical domino, then $r'=r$ and the lemma follows.
Otherwise, if $t=s-1$ and $\nu /\mu' $ is a vertical domino, then
$r'=r-1$ implying $r'<r\leq s=s'$.  Finally, if $t=s-1$ and $\mu
/\lambda '$ is a vertical domino, then $s'=s-1$.  But here the $s$th and
$s-1$st parts of $\mu $ must be equal, so {\it a priori} $r\leq s-1$.
We also get $\mu '\rgrowth{s}\mu $, so $r'=r$.  Thus, $r'\leq s-1=s'$.
So the lemma holds in all cases.

\end{proof}

The following lemma is similarly proven by taking cases.
\begin{lemma}\label{lem:revfive}
Suppose we are given a sequence of consecutive growths $\nu
\rgrowth{r} \mu \rgrowth{s}\lambda $ with $r\ge s$ and suppose $\lambda
\cover \lambda '$ or
$\lambda \covered \lambda '$.  Working backwards by local rules, suppose
that we do not encounter the case $\times $.  Then we obtain
$\nu' \rgrowth{r'} \mu' \rgrowth{s'} \nu '$ with $r' \ge s'$.
\end{lemma}

Returning to the proof of Theorem~\ref{thm:revmain},
 first we show that the circles occur in increasing order in the rows $f - 1$
 and $f$ by using Lemma~\ref{lem:revfour}.
Since $\bigcirc$ cannot occupy the same column in both rows (this would place
 a growth and a shrink at the same time on the edge in between),
 it is enough to see that $\bigcirc$ cannot appear earlier (from the right,
 since we are working backwards) in row $f - 1$.
Let $\Bar{ \Bar \kap } \th j$ denote the shape at $( f - 2, j )$
 in the synthesized picture, for any $j$.
Under the current assumption, if we set $\Bar{ \Bar \kap } \th f
 \rshrink{r_{f}} \bar \kap \th f \rshrink{s_{f}} \kap \th f$,
 then we have $r_{f}\leq s_{f}$.
By Lemma~\ref{lem:revfour}, $\Bar{ \Bar \kap } \th j \rshrink{ r_j } \bar \kap \th j
 \rshrink{ s_j } \kap \th j$ with $r_j \le s_j$ holds
 so long as $\bigcirc$ does not appear.
Let $j^o$ denote the coordinate of the column containing $\bigcirc$ in
 row $f - 1$.
This means that both $\Bar{ \Bar \kap } \th{ j^o } \rshrink l
 \bar \kap \th { j^o }$ and $\bar \kap \th{ j^o - 1 } \rshrink l
 \bar \kap \th { j^o }$ occur, and that we have a shrink on
 $( f - 1, j^o )$--$( f, j^o )$.
By what we just saw, the first of these two conditions implies
 $\bar \kap \th { j^o } \rshrink l \kap \th{ j^o }$, since $l$ is
 the bottommost possible row of these shapes.
Consulting Theorem~\ref{ThmLocRules}, we find no local rule which would produce
 the combination $\bar \kap \th{ j^o - 1 } \rshrink l \bar \kap \th{ j^o }
 \rshrink l \kap \th{ j^o }$ at the cell $( f, j^o )$,
 which means that $\bigcirc$ cannot appear earlier in row $f - 1$.
Hence $\bigcirc$ appears earlier in the $f$th row.

Now we show that the $\cross$ occur in increasing order in rows $f - 1$
 and $f$ by using Lemma~\ref{lem:revfive}.
Since $\cross$ cannot occupy the same column in both rows, it is enough to see
 that $\cross$ cannot appear earlier (again from the right) in row $f - 1$.
Continuing to denote by $j^o$ the column containing $\bigcirc$ in row $f - 1$,
 let $j^\cross$ denote the column containing $\cross$ in row $f - 1$,
 and assume that $\cross$ appears more to the left in row $f$.
By definition of $\bigcirc$, we have $\Bar{ \Bar \kap } \th{ j^o - 1 }
 \rgrowth l \bar \kap \th{ j^o - 1 }$.
We also have a growth on $( f - 1, j^o - 1 )$--$( f, j^o - 1 )$,
 since $\bigcirc$ has appeared earlier in row $f$.
Since $l$ is the bottommost possible row of these shapes,
 we are in the situation to start using Lemma~\ref{lem:revfive},
 and we can continue to apply it until we reach
 the segments $( f - 2, j^\cross )$--$( f - 1, j^\cross )$--$( f, j^\cross )$.
The occurrence of $\cross$ at $( f - 1, j^\cross )$ implies that we have both
 $\Bar{ \Bar \kap } \th{ j^\cross } \rgrowth 1 \bar \kap \th{ j^\cross }$ and
 $\bar \kap \th{ j^\cross - 1 } \rgrowth 1 \bar \kap \th{ j^\cross }$.
The consequence of the repetitive appliction of Lemma~\ref{lem:revfive} is that we also have
 $\bar \kap \th{ j^\cross } \rgrowth 1 \kap \th{ j^\cross }$,
 since row $1$ is the uppermost row.
Again consulting Theorem~\ref{ThmLocRules}, we find no local rule which would produce
 the combination $\bar \kap \th{ j^\cross - 1 } \rgrowth 1
 \bar \kap \th{ j^\cross } \rgrowth 1 \kap \th{ j^\cross }$.
Hence we must have $\cross$ in the increasing order, as desired.

\end{proof}

\begin{eg}
We show with an example how we can use Theorem~\ref{thm:revmain} (and some arguments not
 included there) to recover the whole picture from $P$ and $Q$.

Suppose we are given a pair $P$ and $Q$ as in Remark~\ref{rem:stratum}.  
Using the symbols in Theorem~\ref{thm:revmain}, we have $l = \gam = m = 2$.
Since $P$ does not contain $\bar 2$ and $Q$ has a $2$-shrink,
 we need to distinguish whether we are in the case (4) or
 (5).
As prescribed in Theorem~\ref{thm:revmain}, we test the assumption that $( 7, 8 )$--%
 $( 8, 8 )$ is a shrink, so that $P \th 7 = \smallmatrix 1 & 1 & 2 \\ 2 & 2
 \endsmallmatrix$.
Working backwards, we obtain the cells in row $8$ and the shapes on row $7$
 as in Figure~\ref{fig:revcor}.  

By (5), we can determine that this assumption is correct.

To go up to the next level, we again find that $l = \gam = m = 2$
 and that $Q$ has a $2$-shrink.
Theorem~\ref{thm:revmain} says that we again need a test,
 which fails this time producing a row with a nonempty leftmost segment:
$$
\def\nowidth#1{\hbox to0pt{\hss$#1$\hss}}
\newdimen\digitwidth
\setbox0=\hbox{1}
\digitwidth=\wd0
\def\thin#1{{\hbox to\digitwidth{\hss$#1$\hss}}}
\begin{matrix}
\text{row 6} \rightarrow &
 1 && 2 && 3 && 3 && \thin{31} && \thin{32} && \thin{42} && \thin{43} &&
 \thin{42} && \smallmatrix 1 & 1 & 2 & 2 \\ 2 & 2 \endsmallmatrix \hfill \\
&& && && && && && && && & \\
\text{row 7} \rightarrow &
 \varnothing && 1 && 2 && 2 && \thin{21} && \thin{22} && \thin{32} &&
  \thin{33} && \thin{32} && \smallmatrix 1 & 1 & 2 \\ 2 & 2 \endsmallmatrix
  \hfill
\end{matrix}
$$
Hence we are forced to take $P \th 6 = \smallmatrix 1 & 1 \\ 2 & 2
 \endsmallmatrix$ as in Figure~\ref{fig:revcor}.  

In practice we can avoid the test in this case.
If $( 6, 8 )$--$( 7, 8 )$ was again a shrink, then $\bigcirc$ in row $7$
 should occur to the left of column $7$, where we have $\bigcirc$ in row $8$.
This is for same reason as used in the proof of Theorem~\ref{thm:revmain}
 after introducing Lemmas~\ref{lem:revfour} and~\ref{lem:revfive}.     
For this to happen, we must have had a $2$-shrink on grid row $7$
to the left of the point $( 7, 6 )$, but there was no such occurrence.

We did not include this argument in Theorem~\ref{thm:revmain} because we stuck
 to the rules only referring to the information in $P$ and $Q$,
 and that are recursively applicable.
This means that, after completing the grid row $7$, we should only look at
 $P \th 7$ and the up-down tableau on row $7$, as if the initially given
 data was of order $7$.

To go up one more level, Theorem~\ref{thm:revmain} imposes a test once again,
 but if we use the fact that row $6$ was already in the $2$-stratum,
 we immediately decide that $( 5, 8 )$--$( 6, 8 )$ is a growth.
Similarly for $( 4, 8 )$--$( 5, 8 )$.

For $( 3, 8 )$--$( 4, 8 )$ we need a test, which succeeds as shown.
Since there is no $1$-shrink on row $3$ to the left of $\bigcirc$ in row $4$,
 we must switch to the $1$-stratum for the remaining rows.

\begin{figure}
\caption{The correct reverse correspondence}
\label{fig:revcor}
$$
\def\nowidth#1{\hbox to0pt{\hss$#1$\hss}}
\newdimen\digitwidth
\setbox0=\hbox{1}
\digitwidth=\wd0
\def\thin#1{{\hbox to\digitwidth{\hss$#1$\hss}}}
\begin{matrix}
\hfill w \rightarrow && \nowidth 2 && \nowidth 2 && \nowidth{ \bar 2 } &&
 \nowidth 1 && \nowidth 1 && \nowidth 2 && \nowidth{ \bar 1 } && \nowidth 1 &&
 & \displaystyle \genfrac{}{}{0pt}{}{ P \th i }{ \downarrow } \\
\text{row 0} \rightarrow &
 \varnothing && \varnothing && \varnothing && \varnothing && \varnothing &&
  \varnothing && \varnothing && \varnothing && \varnothing && \\
\hfill 1 && && && && \nowidth \cross && && && && & \\
\text{row 1} \rightarrow &
 \varnothing && \varnothing && \varnothing && \varnothing && 1 && 1 && 1 &&
  1 && 1 && \smallmatrix 1 \endsmallmatrix \hfill \\
\hfill 1 && && && && && \nowidth \cross && && && & \\
\text{row 2} \rightarrow &
 \varnothing && \varnothing && \varnothing && \varnothing && 1 && 2 && 2 &&
  2 && 2 && \smallmatrix 1 & 1 \endsmallmatrix \hfill \\
\hfill 1 && && && && && && && && \nowidth \cross & \\
\text{row 3} \rightarrow &
 \varnothing && \varnothing && \varnothing && \varnothing && 1 && 2 && 2 &&
  2 && 3 && \smallmatrix 1 & 1 & 1 \endsmallmatrix \hfill \\
\hfill \bar 1 && && && && && && && \nowidth \cross && \nowidth \bigcirc & \\
\text{row 4} \rightarrow &
 \varnothing && \varnothing && \varnothing && \varnothing && 1 && 2 && 2 &&
  3 && 2 && \smallmatrix 1 & 1 \endsmallmatrix \hfill \\
\hfill 2 && \nowidth \cross && && && && && && && & \\
\text{row 5} \rightarrow &
 \varnothing && 1 && 1 && 1 && \thin{11} && \thin{21} && \thin{21} &&
  \thin{31} && \thin{21} && \smallmatrix 1 & 1 \\ 2 \endsmallmatrix \hfill \\
\hfill 2 && && \nowidth \cross && && && && && && & \\
\text{row 6} \rightarrow &
 \varnothing && 1 && 2 && 2 && \thin{21} && \thin{22} && \thin{22} &&
  \thin{32} && \thin{22} && \smallmatrix 1 & 1 \\ 2 & 2 \endsmallmatrix
  \hfill \\
\hfill 2 && && && && && && \nowidth \cross && && & \\
\text{row 7} \rightarrow &
 \varnothing && 1 && 2 && 2 && \thin{21} && \thin{22} && \thin{32} &&
  \thin{33} && \thin{32} && \smallmatrix 1 & 1 & 2 \\ 2 & 2 \endsmallmatrix
  \hfill \\
\hfill \bar 2 && && && \nowidth \cross && && && && \nowidth \bigcirc && & \\
\text{row 8} \rightarrow &
 \varnothing && 1 && 2 && 3 && \thin{31} && \thin{32} && \thin{33} &&
  \thin{32} && \thin{22} && \smallmatrix 1 & 1 \\ 2 & 2 \endsmallmatrix \hfill
\end{matrix}
$$
\end{figure}
\end{eg}

An example of the full correspondence, which coincides with the example
we give earlier by bumping is shown in Figure~\ref{fig:revfull}.   

\begin{figure}
\caption{An example of the reverse correspondence}
\label{fig:revfull}
$$ \smdfyoung{
&\X $\times$.1.&1&1&1&1&1&1&1&1&1&1&1&1&1&1&1&1&1&1\\
&1&1&1&1&1&\X $\times$.2.&2&2&2&2&2&2&2&2&2&2&2&2&2\\
&1&1&1&1&1&2&2&2&2&2&2&2&2&2&2&2&2&\X $\times$.3.&3\\
&1&1&1&1&\X $\times$.2.&\X $\circ$.1.&1&1&1&1&1&1&1&1&1&1&1&2&2\\
&1&1&1&1&2&1&1&1&\X $\times$.2.&2&2&2&2&2&2&2&2&\X $\circ$.1.&1\\
&1&1&1&1&2&1&1&1&2&2&2&2&2&\X $\times$.3.&3&3&3&2&2\\
&1&1&1&1&2&1&\X $\times$.2.&2&21&21&21&21&21&31&31&31&31&21&21\\
&1&1&1&1&2&1&2&2&21&21&\X $\times$.31.&31&31&32&32&32&32&22&22\\
&1&1&1&1&2&1&2&2&21&21&31&31&31&32&\X $\times$.42.&42&42&32&32\\
&1&1&1&1&2&1&2&2&21&21&31&31&31&32&42&\X $\times$.52.&52&42&42\\
&1&1&1&1&2&1&2&2&21&21&31&31&31&32&42&52&52&42&\X $\times$.52.\\
&1&\X $\times$.2.&2&2&21&11&21&21&\X $\circ$.2.&2&3&3&3&31&41&51&51&41&51\\
&1&2&2&2&21&11&21&21&2&\X $\times$.3.&31&31&31&\X $\circ$.3.&4&5&5&4&5\\
&1&2&2&2&21&11&21&21&2&3&31&31&\X $\times$.41.&4&41&51&51&41&51\\
&1&2&2&\X $\times$.3.&31&21&22&22&21&31&311&311&411&41&411&511&511&411&511\\
&1&2&2&3&31&21&22&22&21&31&311&\X $\times$.411.&421&42&421&521&521&421&521\\
\X $\times$.1.&11&21&21&31&311&211&221&221&211&311&\X $\circ$.31.&41&42&41&411&511&511&411&511\\
1&11&21&\X $\times$.31.&32&321&221&222&222&221&321&32&42&421&411&\X $\circ$.41.&51&51&41&51\\
1&11&21&31&32&321&221&222&\X $\times$.322.&321&331&33&43&431&421&42&52&52&42&52\\
1&11&21&31&32&321&221&222&322&321&331&33&43&431&421&42&52&\X $\times$.62.&52&53\\
} $$
\end{figure}

\section{Open Questions and Remarks}\label{sec:oqr}

Although the current pictorial viewpoint allows some additional insight
into the workings of Berele's correspondence, it is not yet the major
simplification that one could hope for.  In particular the difficulty
of running the algorithm backwards is still less than satisfactory.  It
would also be nice if Berele's correspondence could be seen as a
particular case of some more general correspondence between pairs of
up-down tableau and certain kinds of permutation-like objects.  This is
currently under investigation.

The pictorial version of Schensted's algorithm is connected with a
certain poset invariant due to Greene and Kleitman.  It is a natural
question to try to generalize this to the current case, but all efforts
to date have failed.

\end{document}